\documentclass[12pt]{article}
\usepackage{amssymb}
\usepackage{amsfonts, amsmath, amsthm,graphics,amssymb}

\def\theequation{\thesubsection.\arabic{equation}}
\def\theguess{\thesubsection.\arabic{guess}}
\newtheorem{theorem}{Theorem}

\def\thetheorem{\thesubsection.\arabic{theorem}}
\def\theprop{\thesubsection.\arabic{prop}}

\def\thelemma{\thesubsection.\arabic{lemma}}
\def\thecor{\thesubsection.\arabic{cor}}
\newtheorem{exam}[theorem]{Example}
\def\theexam{\thesubsection.\arabic{exam}}

\def\theremark{\thesubsection.\arabic{remark}}
\newcommand{\eqa}{\begin{eqnarray}}
\newcommand{\eeqa}{\end{eqnarray}}
\newcommand{\beq}{\begin{equation}}
\newcommand{\eeq}{\end{equation}}
\newcommand{\nn}{\nonumber}

\newcommand{\pal}{\partial}

\newcommand{\al}{\alpha}

\newcommand{\tr}{{\rm tr}}
\newcommand{\epf}{$\quad$\hfill
\raisebox{0.11truecm}{\fbox{}}\par\vskip0.4truecm}

\newcommand{\res}{\mathrm {Res}}

\def\otim{\mathop{\otimes}}

\def\ddz{{{\rm d}\over{\rm d} z}}

\def\diag{{\rm diag}}

\newsavebox{\uuunit}
\sbox{\uuunit}
    {\setlength{\unitlength}{0.825em}
     \begin{picture}(0.6,0.7)
        \thinlines
        \put(0,0){\line(1,0){0.5}}
        \put(0.15,0){\line(0,1){0.7}}
        \put(0.35,0){\line(0,1){0.8}}
       \multiput(0.3,0.8)(-0.04,-0.02){12}{\rule{0.5pt}{0.5pt}}
     \end {picture}}
\newcommand{\ID}{\mathord{\!\usebox{\uuunit}}}
\newcommand{\complessi}{\mathbb C}

\begin{document}
\theoremstyle{definition}
\newtheorem{df}[theorem]{Definition}
\newtheorem{ex}[theorem]{Example}
\theoremstyle{theorem}
\newtheorem{lm}[theorem]{Lemma}
\newtheorem{rmk}[theorem]{Remark}
\newtheorem{thm}[theorem]{Theorem}
\newtheorem{cor}[theorem]{Corollary}
\newtheorem{con}[theorem]{Conjecture}
\newtheorem{prop}[theorem]{Proposition}

\title{{\bf Canonical structure and symmetries of the Schlesinger equations.}}
\author{Boris Dubrovin\thanks{SISSA, International School of Advanced
Studies, via Beirut 2-4, 34014 Trieste, Italy} and Marta 
Mazzocco\thanks{School of Mathematics, The University of Manchester, Manchester M60 1QD, 
United Kingdom.}} 
\date{}
\maketitle

\begin{abstract} The Schlesinger equations  $S_{(n,m)}$
describe monodromy preserving deformations of order $m$ Fuchsian systems
with $n+1$ poles. They can be considered as a family of commuting 
time-dependent
Hamiltonian systems on the direct product of $n$ copies of $m\times m$
matrix algebras equipped with the standard linear Poisson bracket.
In this paper we present a new canonical Hamiltonian formulation of the 
general Schlesinger equations $S_{(n,m)}$ for all $n$, $m$ and we
compute the action of the symmetries of the Schlesinger equations in  these coordinates.
\end{abstract}

\tableofcontents

\def\theequation{\thesection.\arabic{equation}}
\def\theguess{\thesection.\arabic{guess}}
\def\thetheorem{\thesection.\arabic{theorem}}
\def\theprop{\thesection.\arabic{prop}}
\def\thelemma{\thesection.\arabic{lm}}
\def\thecor{\thesection.\arabic{cor}}
\def\theexam{\thesection.\arabic{ex}}
\def\theremark{\thesection.\arabic{rmk}}
\setcounter{equation}{0}
\setcounter{theorem}{0}

\section{Introduction.} 

The {\it Schlesinger equations} $S_{(n,m)}$ \cite{Sch} is the following 
system of nonlinear differential equations 
\begin{eqnarray}
&&
{\partial\over\partial u_j} {A}_i= 
{[ {A}_i, {A}_j]\over u_i-u_j},\qquad i\neq j,\nn \\
&&
{\partial\over\partial u_i} {A}_i= 
-\sum_{j\neq i}{[ {A}_i, {A}_j]\over u_i-u_j},\label{sch}
\end{eqnarray} 
for $m\times m$ matrix valued functions ${A}_1=A_1(u),\dots,
{A}_{n}=A_n(u)$,
$u=(u_1,\dots,u_n)$, the independent variables $u_1$, \dots, $u_n$ must be 
pairwise distinct. The first non-trivial case $S_{(3,2)}$
of the Schlesinger equations corresponds to the 
famous sixth Painlev\'e equation \cite{fuchs, Sch, Gar1}, 
the most general of all Painlev\'e
equations. In the case of any number $n>3$ of $2\times2$ 
matrices ${A}_j$, the Schlesinger equations $S_{(n,2)}$ reduce to the
Garnier systems ${\cal G}_n$ (see \cite{Gar1,Gar2,Ok1}).

The Schlesinger equations $S_{(n,m)}$ appeared in the theory of 
{\it isomonodromic deformations}\/
of Fuchsian systems. Namely, the monodromy matrices of the Fuchsian system
\begin{equation}
{d \Phi\over  d{z}} =\sum_{k=1}^{n}{{A}_k(u)\over {z}-u_k}\Phi,
\qquad\qquad {z}\in\complessi\backslash\{u_1,\dots,u_{n}\}
\label{N1}\end{equation}
do not depend on $u=(u_1, \dots, u_n)$ if the matrices $A_i(u)$ satisfy
(\ref{sch}). Conversely, under certain assumptions for the 
matrices $A_1$, \dots, $A_n$ and for the matrix
\begin{equation}\label{a-inf}
A_\infty:= -\left( A_1+\dots +A_n\right),
\end{equation}
{\it all} isomonodromic deformations of the Fuchsian system (\ref{N1})
are given by solutions to the Schlesinger equations (see, e.g.,
\cite{sib})\footnote{Bolibruch constructed non-Schlesinger isomonodromic 
deformations in \cite{Bol}. These can occur when the matrices $A_i$ are 
resonant,
i.e. admit pairs of eigenvalues with positive integer differences.}.

The solutions to the Schlesinger equations can be parameterized by the
{\it monodromy data} of the Fuchsian system (\ref{N1}) (see precise 
definition below 
in Section 2). To reconstruct the solution starting from given monodromy 
data one is to solve the classical {\it Riemann - Hilbert problem} of 
reconstruction of the Fuchsian system
from its monodromy data. The main outcome of this approach says that
the solutions $ {A}_i(u)$ can be continued 
analytically to meromorphic functions on the universal covering of
$$
\left\{(u_1,\dots,u_{n})\in\complessi^{n}\,|\,u_i\neq u_j\,\hbox{for}\, 
i\neq j\right\}
$$
\cite{Mal,Mi}.
This is a generalization of the celebrated {\it Painlev\'e
property}\/ of absence of movable critical singularities (see details in
\cite{ince,its}). In certain cases the technique based on the theory of
Riemann--Hilbert problem gives a possibility to compute the asymptotic 
behavior
of the solutions to Schlesinger equations near the critical locus $u_i=u_j$
for some $i\neq j$, although there are still interesting open problems
\cite{jimbo, DM, guzzetti, costin}. 

It is the Painlev\'e property that was used by Painlev\'e and Gambier
as the basis for their classification scheme of nonlinear differential
equations. Of the list of some 50 second order nonlinear differential equations
possessing Painlev\'e property
the six (nowadays known as {\it Painlev\'e equations}) 
are selected due to the following crucial property: the general solutions
to these six equations cannot be expressed in terms of {\it classical
functions}, i.e., elementary functions, elliptic and other classical
transcendental functions (see \cite{Um1} for a modern approach to this
theory based on a nonlinear version of the differential Galois theory). 
In particular, according to these results
the general solution to the Schlesinger system $S_{(3,2)}$ corresponding 
to the Painlev\'e-VI equation cannot be expressed in terms of classical functions.

A closely related question is the problem of construction and classification
of {\it classical solutions} to the Painlev\'e equations and their generalizations.
This problem remains open even for the case of Painlev\'e-VI although there are
interesting results based on the theory of symmetries of the Painlev\'e equations
\cite{OK, OK3, kitaev} and on the geometric approach to studying the space 
of monodromy
data \cite{DM, Hit, M,M2}.

The above methods do not give any clue to solution of 
the following general problems: are 
solutions of $S_{(n+1,m)}$ or of $S_{(n,m+1)}$ more complicated than those
of $S_{(n,m)}$? Which solutions to $S_{(n+1,m)}$ or $S_{(n,m+1)}$
can be expressed via solutions to $S_{(n,m)}$? Furthermore, which of them can
ultimately be expressed via classical functions?

Interest to these problems was one of the starting points for our work. 
We began to look at the theory of {\it symmetries} of Schlesinger
equations, i.e., of birational transformations acting in the space of Fuchsian
systems that map solutions to solutions. One class of symmetries is wellknown
 \cite{MJ1, MJ2}:
they are the so-called Schlesinger transformations
\begin{equation}\label{sch-tr}
A(z)=\sum_{i=1}^n {A_i\over z-u_i}\mapsto \tilde A(z) = 
{d G(z)\over dz} G^{-1}(z) +G(z) A(z)
G^{-1}(z)= \sum_{i=1}^n {\tilde A_i\over z-u_i} 
\end{equation}
with a rational invertible $m\times m$ matrix valued function $G(z)$
preserving the class of Fuchsian systems. Clearly such transformations
preserve the monodromy of the Fuchsian system.

More general symmetries of the $S_{(3,2)}$ Schlesinger
equations {\it do not} preserve the monodromy. They can be derived
from the theory, due to K.Okamoto \cite{Ok2}, of canonical 
transformations of the 
Painlev\'e-VI equation considered as a time-dependent Hamiltonian system
(see also \cite{manin, arinkin} regarding
an algebro-geometric approach to some of the Okamoto symmetries).
Some of the Okamoto symmetries were later generalized to the Schlesinger 
systems
$S_{(n,2)}$ with arbitrary $n>3$ \cite{OK,Ts} (see also \cite{IKSY}) using the
Hamiltonian formulation of the related Garnier equations.
The generalization of these symmetries to $S_{(n,m)}$ with
arbitrary $n$, $m$ was one of the motivations for our work.

With this problem in mind, in this paper we present a canonical Hamiltonian 
formulation of Schlesinger equations $S_{(n,m)}$ for all $n$, $m$.

Recall \cite{MJU,Man1} that Schlesinger equations can be written as Hamiltonian
systems on the Lie algebra
$$
{\mathfrak g}:= \oplus_{i=1}^n{\mathfrak{gl}}(m) \ni (A_1, \dots, A_n)
$$
with respect to the standard linear Lie - Poisson bracket on 
${\mathfrak g}^\star\sim {\mathfrak g}$
with quadratic time-dependent Hamiltonians of the form
\begin{equation}\label{hjm}
H_k:=\sum_{l\neq k}{{\rm tr}\left(A_k A_l\right)\over u_k-u_l}.
\end{equation}
Because of isomonodromicity they can be restricted onto the symplectic leaves
$$
{\cal O}_1 \times \dots \times {\cal O}_n \in {\mathfrak g}
$$
obtained by fixation of the conjugacy classes ${\cal O}_1$,\dots,
${\cal O}_n$ of the matrices $A_1$, \dots, $A_n$. 
The matrix $A_\infty$ given in (\ref{a-inf}) is a common integral of 
the Schlesinger
equations. Applying the procedure of symplectic reduction \cite{mars-wein}
we obtain the reduced symplectic space 
\begin{eqnarray}\label{level}
&&
\left\{ A_1\in {\cal O}_1 , \dots, A_n\in {\cal O}_n,\,
A_\infty = \hbox{given diagonal matrix}\right\}
\nn\\
&&
\qquad\hbox{modulo simultaneous diagonal conjugations}.
\end{eqnarray}
The dimension of this reduced symplectic leaf in the generic situation
is equal to $2g$ where
$$
g={m(m-1)(n-1)\over2}-(m-1).
$$
Our aim is to introduce ``good'' canonical Darboux coordinates on 
the generic reduced symplectic leaf (\ref{level}).

Actually, there is a natural system of canonical coordinates on (\ref{level}):
it is obtained in \cite{scott, gekhtman} within the general framework of 
algebro-geometric
Darboux coordinates introduced by S.Novikov and A.Veselov \cite{veselov}
(see also \cite{adams, DuDi}). They
are given by the $z$- and $w$-projections of the points of the divisor of 
a suitably normalized
eigenvector of the matrix $A(z)$ considered as a section of a line bundle
on the spectral curve
\begin{equation}\label{s-curve}
\det (A(z) - w \ID) =0.
\end{equation}
However, the symplectic reduction (\ref{level}) is time--dependent. This produces a shift 
in the Hamiltonian functions that can only be computed by knowing the 
explicit parameterization of the matrices $A_1,\dots, A_n$ by the spectral 
coordinates. The same difficulty makes the computation of 
the action of the symmetries on the spectral coordinates for $m>2$ essentially impossible.

Instead, we construct a new system of the so-called {\it isomonodromic Darboux
coordinates} $q_1$, \dots, $q_g$, $p_1$, \dots, $p_g$
on generic symplectic manifolds (\ref{level}) and we give the new Hamiltonians 
in these coordinates. 
Let us explain our construction.

The Fuchsian system (\ref{N1}) can be reduced to a scalar differential
equation of the form
\begin{equation}\label{skal}
y^{(m)} = \sum_{l=1}^{m-1} d_l(z) y^{(l)}.
\end{equation}
For example, one can eliminate last $m-1$ components of 
$\Phi$ to obtain a $m$-th order equation for the first component
$y:=\Phi_1$. (Observe that the reduction procedure depends on the choice of
the component of $\Phi$).
The resulting Fuchsian equation will have regular singularities
at the same points $z=u_1$, \dots, $z=u_n$, $z=\infty$. It will also have
other singularities produced by the reduction procedure. However, they will be
{\it apparent} singularities, i.e., the solutions to (\ref{skal}) will be
analytic in these points. Generically there will be exactly $g$ apparent
singularities (cf. \cite{oht}; a more precise result about the number
of apparent singularities working also in the nongeneric situation was obtained
in \cite{Bol2}); they will be the first part $q_1$, \dots, $q_g$ of the
canonical coordinates. The conjugated momenta are defined by
$$
p_i = \res_{z=q_i} \left( d_{m-2}(z) +{1\over 2} d_{m-1}(z)^2\right), \quad i=1,
\dots, g.
$$
Our first results is

\begin{thm} Let the eigenvalues of the matrices $A_1$, \dots, $A_n$, $A_\infty$
be pairwise distinct. Then the map
\begin{equation}\label{theor}
\left\{\begin{matrix} {\rm Fuchsian} ~ {\rm systems} ~ {\rm with} ~{\rm given}
~{\rm poles}\\
{\rm and}~ {\rm given} ~ {\rm eigenvalues} ~ {\rm of} ~ A_1, \dots, A_n,
A_\infty\\
{\rm modulo} ~ {\rm diagonal} ~ {\rm conjugations} \end{matrix}\right\} \to
(q_1, \dots, q_g, p_1, \dots, p_g)
\end{equation}
gives a system of rational
Darboux coordinates on a large Zariski open set\footnote{A precise characterization of this large open set is given in Theorem \ref{skal-red} and Remark \ref{rmknonzero}.} in the generic reduced symplectic leaf (\ref{level}).
The Schlesinger equations $S_{(n,m)}$ in these coordinates are written in the
canonical Hamiltonian form
\begin{eqnarray}
&&
{\pal q_i\over \pal u_k} = {\pal{\mathcal H}_k\over \pal p_i}
\nn\\
&&
{\pal p_i\over \pal u_k} = -{\pal {\mathcal H}_k\over \pal q_i}\nn
\end{eqnarray}
with the Hamiltonians
$$
{\mathcal H}_k = {\mathcal H}_k(q, p; u) = -\res_{z=u_k} \left( d_{m-2}(z) +{1\over 2} d_{m-1}(z)^2
\right), \quad k=1, \dots, n.
$$
\label{th.intro.1}\end{thm}

Here {\it rational Darboux coordinates} means that the elementary
symmetric functions $\sigma_1(q)$, \dots, $\sigma_g(q)$ and 
$\sigma_1(p)$, \dots, $\sigma_g(p)$ are rational functions of the coefficients
of the system and of the poles $u_1$, \dots, $u_n$. Moreover, there exists
a section of the map (\ref{theor}) given by rational functions
\begin{equation}\label{sect}
A_i = A_i(q, p), \quad i=1, \dots, n,
\end{equation}
symmetric in $(q_1, p_1)$, \dots, $(q_g, p_g)$ with coefficients depending on $u_1,\dots,u_n$ and on the eigenvalues of the matrices $A_i$, $i=1, \dots, n, \infty$.
All other Fuchsian systems with the same poles $u_1,\dots,u_n$, the same exponents and the same
$(p_1,\dots,p_g,q_1,\dots,q_g)$ are obtained by simultaneous diagonal conjugation
$$
A_i(q,p)\mapsto C^{-1}A_i(q,p)C, \quad i=1, \dots, n, 
\quad C={\rm diag}\, (c_1,
\dots, c_m).
$$

In the course of the proof of Theorem \ref{th.intro.1}, we establish that the same
parameters $(p_1,\dots,p_g,q_1,\dots,q_g)$  are rational coordinates in the space 
of what we call {\it special Fuchsian equations}, i.e., $m$-th order Fuchsian 
equations with
$n+1$ regular singularities with given exponents
and $g$ apparent singularities with the exponents\footnote{As it 
was discovered by H.Kimura and  K.Okamoto \cite{KO1} these are the 
exponents at the generic apparent singularities of the scalar
reduction of a Fuchsian system.}
$0$, $1$, \dots, $m-2$, $m$. We then prove that there is a birational map between
special Fuchsian equations and Fuchsian systems. This allows us to conclude that 
 $(p_1,\dots,p_g,q_1,\dots,q_g)$  are rational coordinates in the space of 
Fuchsian systems.

The natural action of the symmetric group $S_n$ on the Schlesinger equations 
is described in the following 

\begin{thm}
The Schlesinger equations $S_{(n,m)}$ written in the canonical form of Theorem 
\ref{th.intro.1} admit  a group of birational canonical 
transformations $\langle S_2,\dots, S_m,S_\infty\rangle$
\begin{equation}
S_k:\quad \left\{\begin{array}{l} 
\tilde q_i= u_1+u_k-q_i, \quad i=1,\dots,g,\\
\tilde p_i=-p_i,\quad i=1,\dots,g,\\
\tilde u_l = u_1+u_k-u_l,\quad l=1,\dots,n, \\
\tilde {\mathcal H}_l=-{\mathcal H}_l,\quad l=1,\dots,n,
\end{array}\right.
\label{isimmk}\end{equation}
\begin{equation}
S_\infty:\quad \left\{\begin{array}{l} 
\tilde q_i={1\over q_i-u_1}, \quad i=1,\dots,g,\\
\tilde p_i= -p_iq_i^2-{2m^2-1\over m}q_i,\quad i=1,\dots,g,\\
\tilde u_l= {1\over u_l-u_1},\quad l=2,\dots,n, \\
u_{1}\mapsto\infty,\\
\infty\mapsto u_{1},\\
\tilde{\mathcal H}_1=\mathcal H_1,\\
\tilde {\mathcal H}_l= -\mathcal H_l (u_l-u_1)^2 + (u_l-u_1) (d^0_{m-1}(u_l-u_1))^2 -\\
\qquad-(u_l-u_1) {(m-1)(m^2-m-1)\over m}d^0_{m-1}(u_l-u_1),\quad l=2,\dots,n\\
\end{array}\right.
\label{simminf}\end{equation}
where 
$$
d^0_{m-1}(u_k) = \sum_{s=1}^g {1\over u_k-q_s}  - {m\, (m-1)\over 2} 
\sum_{l\neq k} {1\over u_k-u_l}.
$$
The transformation $S_k$ acts on the monodromy matrices as follows
\begin{eqnarray}
&&
\tilde  M_1=M_1^{-1}\dots M_{k-1}^{-1} M_k M_{k-1}\dots M_1, \nn\\
&&
\tilde M_j=M_{j},\quad j\neq 1,k,\nn\\
&&
\tilde M_k=M_{k-1} \dots M_2 M_1 M_2^{-1} \dots M_{k-1}^{-1},\quad i=k+1,\dots,n.\nn
\end{eqnarray}
The transformation $S_\infty$ acts on the monodromy matrices as follows
\begin{eqnarray}
&&
\tilde M_\infty=e^{-{2\pi i \over m}}  C_1 M_\infty^{-1}  M_1  M_\infty C_1^{-1},\nn\\
&&
\tilde M_1= e^{2\pi i \over m} C_1 M_\infty C_1^{-1},\nn\\
&&
\tilde M_j=C_1^{-1}M_j C_1,\quad j\neq 1,\infty,\nn
\end{eqnarray}
where $C_1$ is the connection matrix defined in Section 2.
\end{thm}

Our approach to the construction of Darboux coordinates seems not to work for
nongeneric reduced symplectic leaves. The problem is that, for a nongeneric
orbit the number of apparent singularities of the scalar reduction is bigger
than half of the dimension of the symplectic leaf. The most striking is the
example of rigid Fuchsian systems. They correspond to the extreme nongeneric
case where the reduced symplectic leaf is zero dimensional. One could expect
to have no apparent singularities of the scalar reduction; by no means this is
the case (see \cite{katz}). The study of reductions of Schlesinger systems
on nongeneric reduced symplectic leaves can reveal  new interesting systems
of nonlinear differential equations. Just one example of a nongeneric
situation is basic for the analytic theory of semisimple Frobenius manifolds.
In this case $m=n$, the monodromy group of the Fuchsian system is generated
by $n$ reflections \cite{Dub7}. The dimension of the symplectic leaves
equals
$$
{n(n-1)\over 2} -\left[{n\over 2}\right].
$$
We do not know yet how to construct isomonodromic
Darboux coordinates for this case if the dimension $n$ of the Frobenius manifold is
greater than 3.

Our next result is the comparison of the isomonodromic Darboux coordinates
with those obtained in the framework of the theory of algebro-geometrically
integrable systems (dubbed here {\it spectral Darboux coordinates}). The spectral
Darboux coordinates are constructed as follows. In the generic situation under
consideration the genus of the spectral curve (\ref{s-curve}) equals $g$. The
affine part of the divisor of the eigenvector 
$$
A(z)\psi = w \psi, \quad \psi=(\psi_1, \dots, \psi_m)^T
$$
of the matrix $A(z)$ has degree $g$. Denote $\gamma_1$, \dots, $\gamma_g$
the $z$-projections of the points of the divisor and $\mu_1$, \dots, $\mu_g$
their $w$-projections. These are the spectral Darboux coordinates in the case
under consideration. 

\begin{thm} Let us consider a family of Fuchsian systems 
$$
\epsilon {d \Phi\over dz} = A(z) \Phi, \quad A(z) = \sum_{i=1}^n {A_i\over
z-u_i}
$$
depending on a small parameter $\epsilon$. The matrices $A_1$, \dots,
$A_n$, $A_\infty$ are assumed to be independent of $\epsilon$. Then the
isomonodromic Darboux coordinates of this Fuchsian system have the following
expansion as $\epsilon \to 0$
$$
q_k = \gamma_k +O(\epsilon), \quad p_k = \epsilon^{-1} \mu_k + O(1), \quad k=1,
\dots, g.
$$
Here $\gamma_k, \mu_k$ are the spectral Darboux coordinates of the matrix $A(z)$
defined above.
\end{thm}

We do not even attempt in this paper to discuss physical applications of
our results. However one of them looks particularly attractive. According to an
idea of N.Reshetikhin \cite{res} (see also \cite{harnad}) the wellknown in 
conformal field theory 
Knizhnik--Zamolodchikov equations can be considered as quantization of 
Schlesinger equations. We believe that our isomonodromic canonical 
coordinates could play an important role in the analysis of the quantization 
procedure, somewhat similar to the role played by the spectral
canonical coordinates in the Sklyanin scheme of quantization of integrable
systems \cite{sklyanin}. We plan to address this problem in subsequent
publications.

The paper is organized as follows. In Section 2 we recall the
relationship between Schlesinger equations and isomonodromic deformations of
Fuchsian systems. In Section 3 we discuss the Hamiltonian formulation of 
Schlesinger equations. A formula for the symplectic structure of Schlesinger 
equations recently found by I.Krichever \cite{krichever} proved to be useful for
subsequent calculations with isomonodromic coordinates; we prove that 
this formula is equivalent to the standard one. In Section 4 we construct 
the isomonodromic
Darboux coordinates and establish a birational isomorphism between the space of
Fuchsian systems considered modulo conjugations and the space of special
Fuchsian differential equations. In Section 5 we express the semiclassical
asymptotics of the isomonodromic Darboux coordinates via spectral Darboux
coordinates and  we apply the above results
to constructing the nontrivial symmetries of Schlesinger equations.
The necessary facts from the theory of spectral Darboux coordinates
are collected in the Appendix.




\setcounter{equation}{0}
\setcounter{theorem}{0}

\section{Schlesinger equations as monodromy preserving deformations of Fuchsian
systems.}

In this section we establish our notations, remind few basic definitions and
prove some technical lemmata that will be useful throughout this paper.

The Schlesinger equations ${\cal} S_{(n,m)}$ describe monodromy preserving 
deformations of 
Fuchsian systems (\ref{N1}) with $n+1$ regular singularities at 
$u_1,\dots,u_{n}$, $u_{n+1}=\infty$:
\begin{equation}\label{f-ur}
{{\rm d}\over{\rm d}{z}}  \Phi =\sum_{k=1}^{n}{{A}_k\over {z}-u_k}\Phi,
\qquad\qquad {z}\in\complessi\backslash\{u_1,\dots,u_{n}\}
\end{equation}
${A}_k$ being $m\times m$ matrices independent on ${z}$, 
and $u_k\neq u_l$ for $k\neq l$, $k,l=1,\dots,n+1$. Let us explain 
the precise meaning of this claim.

\subsection{Levelt basis near a logarithmic singularity and local monodromy
data}

A system
\begin{equation} \label{loga}
{d\Phi\over dz} ={A(z)\over z-z_0} \Phi
\end{equation}
is said to have a {\it logarithmic}, or {\it Fuchsian} singularity at $z=z_0$ if
the $m\times m$ matrix valued function $A(z)$ is analytic in some neighborhood of $z=z_0$.
By definition the {\it local monodromy data}\/ of the system is the class of
equivalence of such systems w.r.t. local  gauge transformations
\begin{equation}\label{gauge0}
A(z) \mapsto G^{-1}(z) A(z)\, G(z) + (z-z_0) G^{-1}(z)\partial_z G(z)
\end{equation}
analytic near $z=z_0$ satisfying
$$
\det G(z_0)\neq 0.
$$
The parameters of the local monodromy can be obtained by choosing a suitable
fundamental matrix solution of the system (\ref{loga}). The most general
construction of such a fundamental matrix was given by Levelt \cite{Levelt}.
We will briefly recall this construction in the form suggested in \cite{Dub7}.

Without loss of generality one can assume that $z_0=0$. Expanding the system
near $z=0$ one obtains
\begin{equation}\label{sys-exp}
{d \Phi\over d z} = \left( {A_0\over z} + A_1 + z\, A_2 + \dots
\right)  \Phi.
\end{equation}
Let us now describe the structure of local monodromy data. 

Two linear operators $\Lambda$, $R$ acting in the complex $m$-dimensional space $V$
$$
\Lambda, \, R : V \to V
$$
are said
to form an {\it admissible pair} if the following conditions are fulfilled.

\noindent 1. The operator $\Lambda$ is semisimple and the operator $R$ is
nilpotent.

\noindent 2. $R$ commutes with $e^{2\pi i \Lambda}$,
\begin{equation}\label{lev1}
e^{2\pi i \Lambda} R = R\, e^{2 \pi i \Lambda}.
\end{equation}
Observe that, due to the last condition the operator $R$ satisfies
\begin{equation}\label{lev2}
R(V_\lambda) \subset \oplus_{k\in {\mathbb Z}} V_{\lambda+k} \quad {\rm for}
\quad {\rm any} \quad \lambda \in {\rm Spec} \, \Lambda,
\end{equation}
where $V_\lambda \subset V$ is the subspace of all eigenvectors of 
$\Lambda$ with
the eigenvalue $\lambda$. The last condition says that

\noindent 3. The sum in the r.h.s. of (\ref{lev2}) contains only 
non-negative values of
$k$.

A decomposition
\begin{equation}\label{lev3}
R=R_0 + R_1 +R_2 + \dots
\end{equation}
is defined where
\begin{equation}\label{lev4}
R_k(V_\lambda) \subset V_{\lambda+k} \quad {\rm for}
\quad {\rm any} \quad \lambda \in {\rm Spec} \, \Lambda.
\end{equation}
Clearly this decomposition contains only finite number of terms.
Observe the useful identity
\begin{equation}\label{lev5}
z^\Lambda R \, z^{-\Lambda} = R_0 + z\, R_1 + z^2 R_2 + \dots.
\end{equation}

\begin{thm} For a system (\ref{sys-exp}) with a logarithmic singularity at
$z=0$ there exists a fundamental matrix solution of the form
\begin{equation}\label{lev6}
\Phi(z) = \Psi(z) z^\Lambda z^R
\end{equation}
where $\Psi(z)$ is a matrix valued function analytic near $z=0$ satisfying
$$
\det \Psi(0)\neq 0
$$
and $\Lambda$, $R$ is an admissible pair.
\end {thm}

The formula (\ref{lev6}) makes sense after fixing a branch of logarithm 
$\log z$ near $z=0$.
Note that $z^R$ is a polynomial in $\log z$ due to nilpotency of $R$.

The proof can be found in \cite{Levelt} (cf. \cite{Dub7}). Clearly $\Lambda$
is the semisimple part of the matrix $A_0$; $R_0$ coincides with its nilpotent
part. The remaining terms of the expansion appear only in the {\it resonant
case}, i.e., if the difference between some eigenvalues of $\Lambda$ is a
positive integer. In the important particular case of a diagonalizable matrix
$A_0$,
$$
T^{-1}A_0 T=\Lambda ={\rm diag}\, (\lambda_1, \dots, \lambda_m)
$$
with some nondegenerate matrix $T$, the matrix function $\Psi(z)$ in the 
fundamental matrix solution (\ref{lev6}) can be obtained in the form
$$
\Psi(z) = T \left( \ID + z \,\Psi_1 + z^2 \Psi_2 + \dots \right).
$$
The matrix coefficients $\Psi_1$, $\Psi_2$, \dots of the expansion as well as the 
matrix components $R_1$, $R_2$, \dots of the matrix $R$ (see (\ref{lev3}))
can be found recursively from the equations
$$
[\Lambda, \Psi_k]-k\, \Psi_k =B_k-R_k +\sum_{i=1}^{k-1} \Psi_{k-i} B_i - R_i \Psi_{k-i},
\quad k\geq 1.
$$
Here
$$
B_k:= T^{-1} A_k T, \quad k\geq 1.
$$
If $k_{\rm max}$ is the maximal integer among the differences
$\lambda_i-\lambda_j$ then 
$$
R_k =0 \quad {\rm for} ~ k>k_{\rm max}.
$$
Observe that vanishing of the logarithmic terms in the fundamental matrix
solution (\ref{lev6}) is a constraint imposed only on the first $k_{\rm max}$
coefficients $A_1$, \dots, $A_{k_{\rm max}}$ of the expansion (\ref{sys-exp}). 

It is not difficult to describe the ambiguity in the choice of the admissible
pair of matrices
$\Lambda$, $R$ describing the local monodromy data 
of the system (\ref{sys-exp}). Namely, the diagonal matrix $\Lambda$ is defined up to
permutations of diagonal entries. Assuming the order fixed, the ambiguity in the
choice of $R$ can be described as follows \cite{Dub7}.
Denote ${\mathcal C}_0(\Lambda)\subset GL(V)$ the subgroup consisting of invertible linear operators $G:V\to V$ satisfying
\beq\label{gr-c0}
 z^\Lambda G \,z^{-\Lambda} = G_0+z\, G_1 +z^2 G_2+\dots .
\eeq
The definition of the subgroup can be reformulated \cite{Dub7} in terms of invariance of certain flag in $V$ naturally associated with the semisimple operator $\Lambda$.
The matrix $\tilde R$ obtained from $R$ by the conjugation of the form
\beq\label{eq-r}
\tilde R= G^{-1} R \, G
\eeq
will be called {\it equivalent} to $R$. Multiplying (\ref{lev6}) on the right
by $G$ one obtains another  fundamental matrix solution to the same system of the same
structure
$$
\tilde\Phi(z):=\Psi(z)z^\Lambda z^RG =\tilde \Psi(z) z^\Lambda z^{\tilde R}
$$ 
i.e., $\tilde \Psi(z)$ is analytic at $z=0$ with $\det \tilde\Psi(0)\neq 0$.

The columns of the fundamental matrix (\ref{lev6}) form a distinguished basis
in the space of solutions to (\ref{sys-exp}).

\begin{df} The basis given by the columns of the matrix (\ref{lev6}) is called
{\it Levelt basis} in the space of solutions to (\ref{sys-exp}). 
The fundamental
matrix (\ref{lev6}) is called {\it Levelt fundamental matrix solution}.
\end{df}

The monodromy transformation of the Levelt fundamental matrix solution
reads
\begin{equation}\label{lev7}
\Phi\left(z\, e^{2\pi i}\right) = \Phi(z) M, \quad M= 
e^{2\pi i \Lambda} e^{2\pi i R}.
\end{equation}

To conclude this Section let us denote ${\cal C}(\Lambda,R)$ the subgroup of invertible transformations of the form
\begin{equation}\label{lev9}
{\cal C}(\Lambda,R)=\{ \, G\in GL(V)\, |\, z^\Lambda  G \,  z^{-\Lambda} = \sum_{k\in {\mathbb Z}} G_k z^k \, \mbox{and} \, [G,R]=0\}.
\end{equation}
The subgroups ${\cal C}(\Lambda,R)$ and ${\cal C}(\Lambda,\tilde R)$
associated with equivalent matrices $R$ and $\tilde R$ are conjugated.
It is easy to see that this subgroup coincides with the centralizer
of the monodromy matrix (\ref{lev7})
\begin{equation}\label{lev8}
G\in {\cal C}(\Lambda, R) \quad {\rm iff} \quad G \, 
e^{2\pi i \Lambda} e^{2\pi i R}  = e^{2\pi i \Lambda} e^{2\pi i R} G, \quad \det
G\neq 0.
\end{equation}

Denote 
\begin{equation}\label{lev10}
{\cal C}_0(\Lambda, R) \subset {\cal C}(\Lambda, R)
\end{equation}
the subgroup consisting of matrices $G$ such that the expansion (\ref{lev9})
contains only non-negative powers of $z$. Multiplying the Levelt fundamental matrix
(\ref{lev6}) by a matrix $G\in {\cal C}_0(\Lambda, R)$ one obtains another
Levelt solution to (\ref{sys-exp})
\begin{equation}\label{lev11}
\Psi(z) z^{\Lambda} z^R G = \tilde \Psi(z) z^{\Lambda} z^R.
\end{equation}

In the next Section we will see that the quotient ${\cal C}(\Lambda, R)/
{\cal C}_0(\Lambda, R)$ plays an important role in the theory of monodromy
preserving deformations.

\begin{exam} For 
$$
\Lambda=\left(\begin{matrix} -1 & 0 & 0 \\ 0 & 0 & 0 \\ 0 & 0 &
1\end{matrix}\right), \quad R=\left( \begin{matrix}
0 & 0 & 0\\ a & 0 & 0 \\ c & b & 0\end{matrix}\right)
$$
the quotient ${\cal C}(\Lambda, R)/
{\cal C}_0(\Lambda, R)$ is trivial {\rm iff} $a\, b \neq 0$.
\end{exam}

\subsection{Monodromy data and isomonodromic deformations of a Fuchsian system}

Denote $\lambda^{(k)}_j$, $j=1,\dots,m$, 
the eigenvalues of the matrix ${A}_k$,
$k=1,\dots,n,\infty$ where the matrix ${A}_\infty$ is defined as
$$
{A}_\infty:= -\sum_{k=1}^{n}{A}_k.
$$
For the sake of technical simplicity let us assume that
\begin{equation}
\lambda^{(k)}_i\neq\lambda^{(k)}_j \,\hbox{for}\quad i\neq j,
\qquad k=1, \dots, n, \infty.
\label{N1.3}
\end{equation}
Moreover, it will be assumed that ${A}_\infty$ is a constant diagonal 
$m\times m$ matrix
with eigenvalues $\lambda^{(\infty)}_j$, $j=1,\dots,m$.

Denote $\Lambda^{(k)}$, $R^{(k)}$ the local monodromy data of the Fuchsian
system near the points $z=u_k$, $k=1$, \dots, $n$, $\infty$. 
The matrices $\Lambda^{(k)}$ are all diagonal
\begin{equation}\label{l-r}
\Lambda^{(k)} = {\rm diag}\, (\lambda_1^{(k)}, \dots, \lambda_m^{(k)}), \quad
k=1, \dots, n, \infty.
\end{equation}
and, under our assumptions
$$
\Lambda^{(\infty)} = {A}_\infty.
$$
Recall that the matrix $G\in GL(m,{\mathbb C})$ belongs to the group ${\mathcal C}_0(\Lambda^{(\infty)})$ {\it iff}
\begin{equation}\label{cinf0}
z^{-\Lambda^{(\infty)}}  G \,  z^{\Lambda^{(\infty)}}
=G_0 +\frac{G_1}{z} +\frac{G_2}{z^2}+\dots .
\end{equation}
It is easy to see that our assumptions about eigenvalues of $A_\infty$ imply diagonality of the matrix $G_0$.

Let us also remind that the matrices $\Lambda^{(k)}$ satisfy
\begin{equation}\label{l-sum}
\tr\, \Lambda^{(1)} + \dots + \tr\, \Lambda^{(\infty)} = 0.
\end{equation}
The numbers 
$\lambda_1^{(k)},\dots,\lambda_m^{(k)}$ are called the {\it exponents}\/ of 
the system (\ref{N1}) at the singular point $u_k$.

Let us fix a fundamental matrix solutions of the form (\ref{lev6}) near all singular points $u_1$, \dots, $u_n$, $\infty$. To this end
we are to fix branch cuts on the complex plane and choose the branches of
logarithms $\log(z-u_1)$, \dots, $\log (z-u_n)$, $\log z^{-1}$. We will do it 
in the following way: perform parallel branch cuts $\pi_k$ 
between $\infty$ and each of the $u_k$, $k=1,\dots,n$ along a given (generic)
direction. After this we can fix Levelt fundamental matrices analytic on
\begin{equation}\label{domain}
z\in {\mathbb C} \setminus \cup_{k=1}^n \pi_k,
\end{equation}
\begin{equation}
\Phi_k({z})=T_k\left(\ID +{\cal O}({z}-u_k)\right) 
({z}-u_k)^{\Lambda^{(k)}}({z}-u_k)^{R^{(k)}}, \quad z\to u_k, \quad 
k=1, \dots, n
\label{N6.1}
\end{equation}
and 
\begin{equation}
\Phi({z})\equiv
\Phi_\infty({z})=\left(\ID+{\cal O}({1\over z})\right)
{z}^{-{A}_\infty}{z}^{-{R}^{(\infty)}},\quad\hbox{as}\quad 
{z}\rightarrow\infty,
\label{M3}
\end{equation}
Define the {\it connection matrices}  by 
\begin{equation}
\Phi_\infty({z})=\Phi_k({z}){ C}_k,\label{N3}
\end{equation}
where $\Phi_\infty (z)$ is to be analytically continued in a vicinity of the
pole $u_k$ along the positive side of the branch cut $\pi_k$.

The monodromy matrices ${M}_k$, $k=1,\dots,n,\infty$ are defined
with respect to a basis $l_1,\dots,l_{n}$ of loops in the
fundamental group
$$
\pi_1\left(\complessi\backslash\{u_1,\dots u_{n}\},
\infty\right).
$$ 
Choose the basis in the following way. The loop $l_k$ arrives 
from infinity in a vicinity of 
$u_k$ along one side of the branch cut $\pi_k$ that will be called 
{\it positive}, then it encircles $u_k$ going in anti-clock-wise direction 
leaving all other poles outside and, finally it returns to infinity along 
the opposite side of the branch cut $\pi_k$ called {\it negative}.

Denote $l_j^* \Phi_\infty(z)$ the result of  analytic
continuation of the fundamental matrix $\Phi_\infty(z)$ along the loop
$l_j$. The monodromy matrix ${M}_j$ is defined by
\begin{equation}\label{momo}
l_j^* \Phi_\infty (z) = \Phi_\infty (z) {M}_j, ~~j=1, \dots, n.
\end{equation}
The monodromy matrices satisfy
\begin{equation}
{M}_\infty {M}_{n} \cdots {M}_1=\ID,\qquad
{M}_\infty=\exp\left(2\pi i{A}_\infty\right)
\exp\left(2\pi iR^{(\infty)}\right)
\label{N6}
\end{equation}
if the branch cuts $\pi_1$, \dots, $\pi_{n}$
enter the infinite point according to the order of their labels, i.e., the
positive side of $\pi_{k+1}$  looks at the negative side of $\pi_k$, $k=1,
\dots, n-1$.

Clearly one has
\begin{equation}
{M}_k={ C}_k^{-1} \exp\left(2\pi i \Lambda^{(k)}\right)
\exp\left(2\pi i R^{(k)}\right) 
{ C}_k,\qquad k=1,\dots,n.
\label{N4}
\end{equation}

The collection of the local monodromy data $\Lambda^{(k)}$, $R^{(k)}$ together
with the central connection matrices $C_k$
will be used in order to uniquely fix the Fuchsian system with given poles.
They will be defined up to an equivalence that we now describe. The eigenvalues 
of the diagonal matrices
$\Lambda^{(k)}$
are defined up to permutations. Fixing the order of the eigenvalues, we define
the class of equivalence of the nilpotent part $R^{(k)}$ and of the connection matrices $C_k$ by factoring out the transformations of the form
\begin{eqnarray}\label{class-mon}
&&
R_k \mapsto G_k^{-1} R_k G_k,
\quad
C_k \mapsto G_k ^{-1}C_k G_\infty, \quad k=1, \dots, n,
\nn\\
&&
\nn\\
&& 
G_k \in {\mathcal C}_0 (\Lambda^{(k)}),
\quad G_\infty \in {\cal C}_0(\Lambda^{(\infty)}).
\end{eqnarray}
Observe that the monodromy matrices (\ref{N4}) will transform by a simultaneous conjugation
$$
M_k \mapsto G_\infty^{-1} M_k G_\infty, \quad k=1, 2, \dots, n, \infty.
$$
\begin{df} The class of equivalence (\ref{class-mon}) of the collection 
\begin{equation}\label{m-data}
\Lambda^{(1)}, R^{(1)}, \dots, \Lambda^{(\infty)}, R^{(\infty)},  C_1, \dots, C_n
\end{equation}
is called {\it monodromy data} of the Fuchsian system with respect to a fixed ordering of the eigenvalues of the matrices $A_1$, \dots, $A_n$ and a given
choice of the branch cuts.
\end{df}

\begin{lm}
Two Fuchsian systems of the form (\ref{N1}) with the same poles
$u_1,\dots,u_{n},\infty$ and the same matrix $A_\infty$ coincide, modulo diagonal conjugations if and only if they have 
the same monodromy data 
with respect to the same system of branch cuts $\pi_1,\dots,\pi_{n}$. 
\label{lm2.8}
\end{lm}

\noindent {\bf Proof.} 
Let 
$$
\Phi_\infty^{(1)}({z})=\left(\ID +O({1\over z})\right)z^{-\Lambda^{(\infty)}}
z^{-R^{(\infty)}}, \quad \Phi_\infty^{(2)}({z})=\left(\ID +O({1\over z})\right)z^{-\tilde\Lambda^{(\infty)}}
z^{-\tilde R^{(\infty)}}
$$ 
be
the fundamental matrices of the form (\ref{M3}) of the two Fuchsian 
systems. Using assumption about $A_\infty$ we derive that $\tilde \Lambda^{(\infty)}=\Lambda^{(\infty)}$.
Multiplying $\Phi_\infty^{(2)}({z})$
if necessary on the right by a matrix 
$G\in {\cal C}_0(\Lambda^{(\infty)})$, we can obtain another fundamental
matrix of the second system with 
$$
\tilde R^{(\infty)} = R^{(\infty)}.
$$
Consider the following matrix:
\begin{equation}\label{tYPhi}
Y({z}):= \Phi_\infty^{(2)}({z})[\Phi_\infty^{(1)}({z})]^{-1}.
\end{equation}
$Y({z})$ is an analytic function around infinity: 
\begin{equation}\label{YPhi1}
Y({z})=G_0+{\cal O}\left({1\over{z}}\right),\quad\hbox{as}\, 
{z}\rightarrow\infty
\end{equation}
where $G_0$ is a diagonal matrix.
Since the monodromy matrices coincide, $Y({z})$ is a single valued 
function on the punctured Riemann sphere 
$\overline\complessi\backslash\{u_1,\dots,u_{n}\}$. 
Let us prove that $Y({z})$ is analytic also at the points $u_k$. Indeed, 
having fixed the monodromy data, we can choose the fundamental matrices 
$\Phi_k^{(1)}({z})$ and $\Phi_k^{(2)}({z})$ of the form
(\ref{N6.1}) with the same connection matrices ${ C}_k$
and the same matrices $\Lambda^{(k)}$, $R^{(k)}$. Then near the point
$u_k$, $Y({z})$ is analytic:
\begin{equation}\label{YPhi2}
Y({z})=T_k^{(2)}\left(\ID+{\cal O}({z}-u_k) \right)
\left[T_k^{(1)}\left(\ID+{\cal O}({z}-u_k) \right)\right]^{-1}.
\end{equation}
This proves that $Y({z})$ is an analytic function on all $\overline\complessi$ 
and then, by the Liouville theorem
$Y({z})=G_0$, which is constant. So the two Fuchsian systems coincide, after conjugation by the diagonal matrix $G_0$. 
\epf

\begin{rmk} \label{rem1} The connection matrices are determined, within their 
equivalence
classes by the monodromy matrices if the quotients ${\cal C}(\Lambda^{(k)},
R^{(k)})/ {\cal C}_0(\Lambda^{(k)},
R^{(k)})$ are trivial for all $k=1$, \dots, $n$. In particular this is the case
when all the characteristic exponents at the poles $u_1$, \dots, $u_n$ are
non-resonant.
\end{rmk}

From the above Lemma the following result readily follows.

\begin{thm} \label{iso1} If the matrices ${A}_k(u_1, \dots, u_{n})$ satisfy
Schlesinger equations (\ref{sch}) and the matrix 
$$
{A}_\infty = -({A}_1 + \dots + {A}_n)
$$
is diagonal then all the characteristic exponents
do not depend on $u_1$, \dots, $u_{n}$. The fundamental matrix 
$\Phi_\infty(z)$ can be chosen in such a way that the nilpotent matrix 
${R}^{(\infty)}$ and also all the monodromy matrices are constant in 
$u_1$, \dots, $u_{n}$. 
Moreover, the Levelt fundamental matrices $\Phi_k(z)$ can be chosen in 
such a way that
all the nilpotent matrices $R^{(k)}$ and also all the connection matrices
${\cal C}_k$ are constant. Viceversa, if the deformation 
${A}_k={A}_k(u_1, \dots, u_{n})$ is such that the monodromy data do
not depend on $u_1$, \dots, $u_n$ then the matrices 
${A}_k(u_1, \dots, u_{n})$, $k=1$, \dots, $n$  satisfy Schlesinger
equations.
\end{thm}

At the end of this Section we give a criterion that ensures that the ``naive''
definition of monodromy preserving deformations still gives rise to the
Schlesinger equations.

\begin{thm} Let $A_k={A}_k(u_1, \dots, u_{n})$, $k=1$, \dots, $n$ be a
deformation of the Fuchsian system (\ref{N1}) such that the following conditions
hold true.

\noindent 1. The matrix $A_\infty=-A_1-\dots -A_n$ is constant and diagonal.

\noindent 2. The Fuchsian system admits a fundamental matrix solution
of the form (\ref{M3}) with the $u$-independent matrix $R^{(\infty)}$.
Denote $\Phi_\infty(z;u)$ the fundamental matrix solution of the family of
Fuchsian systems of the form (\ref{M3}).

\noindent 3. The monodromy matrices $M_1$, \dots, $M_n$ defined as in
(\ref{momo}) with respect to the fundamental matrix $\Phi_\infty(z;u)$ do not
depend on $u_1$, \dots, $u_n$. Note that this implies constancy of the diagonal
matrices $\Lambda^{(k)}$ of exponents, $k=1, \dots, n$.

\noindent 4. The (class of equivalence of) local monodromy data $(\Lambda^{(k)},
R^{(k)})$ does not depend on $u_1, \dots, u_n$.

\noindent 5. The quotients 
${\cal C}(\Lambda^{(k)}, R^{(k)})/{\cal C}_0(\Lambda^{(k)},
R^{(k)})$ are zero dimensional for all $k=1, \dots, n$.

\noindent Then the deformation satisfies the Schlesinger equations. Moreover, under
the assumption of validity of 1 - 4, if the condition 5 does not hold true,
then there exist non-Schlesinger deformations 
preserving the monodromy matrices.
\end{thm}

\noindent{\bf Proof.} The first statement of the Theorem easily follows from Remark \ref{rem1}
and Theorem \ref{iso1}. To
prove the second part, let us assume that the dimension of the quotient
${\cal C}(\Lambda^{(k)}, R^{(k)})/{\cal C}_0(\Lambda^{(k)},
R^{(k)})$  is positive for some $k$. Here $\Lambda^{(k)}$, $R^{(k)}$ are local
monodromy data of the Fuchsian system (\ref{N1}) with some poles $u_1$, \dots,
$u_n$, $\infty$.
Let us choose a nontrivial family of matrices
$G(s) \in {\cal C}(\Lambda^{(k)}, R^{(k)})/{\cal C}_0(\Lambda^{(k)},
R^{(k)})$ for sufficiently small $s$, $G(0)=\ID$. We will now obtain a
deformation of the Fuchsian system (\ref{N1}) in the following way. Let us
deform the $k$-th connection matrix $C_k$ by putting
$$
C_k(s):= G(s) C_k.
$$
To reconstruct the deformation of the Fuchsian system, we are to solve the
suitable Riemann - Hilbert problem. It will be solvable for sufficiently small
$s$ because of solvability for $s=0$. At this point one can also deform the
poles $u_i(s)$, $u_i(0)=u_i$. This deformation is obviously isomonodromic but
not of the Schlesinger type. The Theorem is proved.\epf



\setcounter{equation}{0}
\setcounter{theorem}{0}

\section{Hamiltonian structure of the Schlesinger system}

\subsection{Lie-Poisson brackets for Schlesinger system}

The Hamiltonian description of the Schlesinger system can be derived
\cite{Hit2}
from the general construction of a Poisson bracket on the
space of flat connections in a principal $G$-bundle over a surface with 
boundary
using Atiyah - Bott symplectic structure (see \cite{audin}). Explicitly this
approach yields the following well known formalism representing
the Schlesinger system ${S}_{(n,m)}$ in
Hamiltonian form with $n$ time variables $u_1,\dots,u_{n}$ and
$n$ commuting time--dependent Hamiltonian flows on the dual space to the 
direct sum of $n$ copies
of the Lie algebra ${\mathfrak{sl}}(m)$
\begin{equation}\label{lie}
{\mathfrak g}:=\oplus_{n} {\mathfrak{sl}}(m)\ni 
\left( {A}_1, {A}_2, \dots, {A}_{n} \right).
\end{equation}
The standard Lie-Poisson bracket on ${\mathfrak g}^*$
reads
\begin{equation}
\left\{\left({A}_p\right)^i_j,\left({A}_q\right)_l^k\right\}
= \delta_{pq}\left(\delta_l^i \left({A}_p\right)_j^k- 
\delta_j^k \left({A}_q\right)^i_l\right).
\label{ham2}\end{equation}
(we identify ${\mathfrak{sl}}(m)$ with its dual by using the Killing form
$(A,B)={\rm Tr}\, AB, \quad A, B \in {\mathfrak{sl}}(m)$.) The 
following statement is well-known
(see \cite{MJ1,Man1}) and can be checked by a straightforward computation.

\begin{thm} \label{thm1.1}
The dependence of the solutions ${A}_k$, $k=1,\dots,n$, of the
Schlesinger system ${S}_{(n,m)}$ upon the
variables $u_1,\dots,u_{n}$ is determined by  Hamiltonian systems 
on (\ref{lie})
with time-dependent quadratic Hamiltonians 
\begin{eqnarray}
&&
H_k= \sum_{l\neq k}
{{\rm Tr}\left({A}_k{A}_l\right)\over u_k-u_l},
\label{ham0}\\
&&
{\partial\over\partial u_k}{A}_l = \{{A}_l,H_k\}.\label{ham000}
\end{eqnarray}
\end{thm}

To arrive from the Hamiltonian systems (\ref{ham0}), (\ref{ham000}) to the isomonodromic
deformations one is to impose an additional constraint. Define
\begin{equation}\label{ainf}
{A}_\infty: = - ({A}_1 + \dots +{A}_n).
\end{equation}
It can be easily seen that 
\begin{equation}\label{int-inf}
{\partial \over \partial u_i} {A}_\infty = \{ {A}_\infty, H_i\} =0,
\quad i=1, \dots, n.
\end{equation} 
So, the matrix entries of ${A}_\infty$ are integrals of the Schlesinger
equations. They generate the action of the group ${\mathfrak{sl}}(m)$ on ${\mathfrak g}$
by the diagonal conjugations
\begin{equation}\label{act-diag}
{A}_k \mapsto G^{-1} {A}_k \, G, \quad k=1, \dots, n, \quad G\in
{\mathfrak{sl}}(m).
\end{equation}
To identify the Hamiltonian equations (\ref{ham000}) with isomonodromic
deformations one is to apply the Marsden - Weinstein procedure of symplectic
reduction \cite{mars-wein}. In our setting this procedure works as follows.
Let us
choose a particular level surface of the moment map
corresponding to the first integrals (\ref{int-inf}). We will mainly deal with
the level surfaces of the form
\begin{equation}\label{moment}
{A}_\infty = {\rm diag}\, \left(\lambda_1^{(\infty)}, \dots,
\lambda_m^{(\infty)}\right), \quad \lambda_i^{(\infty)} \neq
\lambda_j^{(\infty)}, \quad i\neq j
\end{equation}
for some pairwise distinct numbers $\lambda_1^{(\infty)}, \dots,
\lambda_m^{(\infty)}$.

After restricting the Hamiltonian systems onto the level surface (\ref{moment})
there remains a residual symmetry with respect to conjugations by diagonal
matrices
\begin{equation}\label{diag}
{A}_k \mapsto D^{-1} {A}_k \, D, \quad k=1, \dots, n, \quad D={\rm
diag}\, (d_1, \dots, d_n).
\end{equation}

Denote 
$$
Diag\simeq \left[{\mathbb C}^*\right]^{m-1} \subset SL(m)
$$
the subgroup of diagonal matrices acting on ${\mathfrak g}$ by simultaneous
conjugations (\ref{diag}).

\begin{df} The quotient of the restriction of the Hamiltonian system
(\ref{ham0}), (\ref{ham000}) onto the level surface of the form (\ref{moment}) 
w.r.t. the transformations (\ref{diag}) will be called
{\it reduced Schlesinger system}.
\end{df}

From the above results it readily follows that, the reduced Schlesinger
system describes all nontrivial monodromy preserving deformations of the
Fuchsian system (\ref{N1}).


\subsection{Symplectic structure of the isomonodromic deformations.}

The symplectic leaves of the Poisson bracket (\ref{ham2})
are products of adjoint orbits
\begin{equation}\label{leaves}
({A}_1,\dots,{A}_{n})\in
{\cal O}_1\times\dots\times{\cal O}_{n}\subset {\mathfrak g}
\end{equation}
where ${\cal O}_k$ is the adjoint orbit of ${A}_k$. The
symplectic structure $\omega_H$ induced by (\ref{ham2}) on the orbits 
can be represented in 
the following form \cite{Hit2}.
Given two tangent vectors
\begin{equation}\label{tangent}
\delta_1{A}=\sum_k{\delta_1{A}_k\over z-u_k}, \quad
\delta_2{A}=\sum_k{\delta_2{A}_k\over z-u_k}
\end{equation} 
at the given
point ${A}({z})$  
the value of the symplectic form can be computed
by 
\begin{equation}
\omega_H(\delta_1 {A},\delta_2{A})=-\sum_{k}{\rm Tr}(U_k^{(1)} \delta_2
{A}_k),
\label{omegah}\end{equation}
where the matrices $U_k^{(i)}$ are such that 
$\delta_i{A}_k=[U_k^{(i)},{A}_k]$. 

Actually, the symplectic structure (\ref{omegah}) was obtained in \cite{Hit2} from the general Atiyah - Bott Poisson structure \cite{AB} on the moduli space of flat $SL(m)$ connections on the surface $\Sigma$. It is obtained by projecting the gauge invariant symplectic form
\begin{equation}\label{ab}
\omega_{AB}(\delta_1 A, \delta_2 A) =\frac1{2\pi i} \int_\Sigma {\rm Tr}\, (\delta_1 A\wedge \delta_2 A)
\end{equation}
onto the moduli space.
In our case $\Sigma$ is the Riemann sphere without small discs around the poles
of $A(z)$. 

The eigenvalues of the matrices ${A}_k$ are 
Casimirs of the Poisson bracket (\ref{ham2}). We will mainly consider
the generic case where these eigenvalues are distinct 
for all $k=1$, \dots, $n$.
Then the level sets of the Casimirs coincide with the symplectic leaves
(\ref{leaves}). 

Denote 
$$
\left( {\rm Spec}\, A_1, \dots, {\rm Spec}\, A_n\right) 
$$
the collection of the eigenvalues of the matrices $(A_1, \dots,
A_n)\in{\mathfrak g}$. Generically these are the parameters of the symplectic
leaves. Fixing the level surface of the moment map (\ref{moment}) and taking the
quotient over the action (\ref{diag}) of the group $Diag\subset SL(m)$ of
diagonal matrices one obtains a manifold that we denote
\begin{equation}\label{red-leaf}
{M}\left( {\rm Spec}\, A_1, \dots, {\rm Spec}\, A_n; A_\infty\right).
\end{equation}
The dimension of this manifold is equal to $2g$ where
\begin{equation}\label{rod}
g={m(m-1)(n-1)\over 2} -(m-1).
\end{equation}
Indeed, the dimension of a generic adjoint orbit ${\mathcal O}_i$ is equal
to $m^2 -m$. Choosing a level surface (\ref{moment}) of the momentum map
$$
(A_1, \dots, A_n)\mapsto A_\infty:=-(A_1+\dots +A_n)
$$  
we
impose only $m^2-1$ independent equations, since the trace of the matrix
$A_\infty$ is equal to the sum of traces of $A_1$, \dots, $A_n$. Finally,
subtracting the dimension $m-1$ of the group $Diag$ we arrive at (\ref{rod}).
The manifold still carries a symplectic structure since the action of the Abelian group
$Diag$ preserves the Poisson bracket (\ref{ham2}). One of the aims of this paper
is to construct a system of canonical Darboux coordinates on generic manifolds
of the form (\ref{red-leaf}). 

For our aims it will be useful the following approach to the Hamiltonian theory
of monodromy preserving deformations developed
recently by Krichever \cite{Kr}. He has obtained a general formula for
the symplectic structure on the space of isomonodromic deformations of 
a generic linear system of ODEs of the form
$$
\ddz\Phi={A}({z})\Phi,
$$
where ${z}$ is a variable on a punctured genus $g$ algebraic curve and 
$A({z})$ is any meromorphic matrix function of ${z}$ with poles 
of any order at $P_1,\dots,P_n$.
In the case of genus $g=0$ the formula reads
\begin{equation}
\omega_K=-{1\over2}\sum_k \res_{z=P_k}{\rm Tr}
(\delta{A}\wedge\delta\Phi\,\Phi^{-1})
\label{omegak}\end{equation}
where the variations $\delta{A}$ and $\delta\Phi$ are
independent. The corresponding Hamiltonian function describing the 
isomonodromic deformations in the parameter $P_j$ is
\begin{equation}
H_{K_j}=-{1\over2} \res_{z=P_j}{\rm Tr}
({A}({z})^2).
\label{hamk}
\end{equation}
Due to gauge invariants \cite{Kr} the symplectic form admits a restriction
 onto the manifold (\ref{leaves}). It defines therefore a symplectic
 structure $\omega_K$ on the product of adjoint orbits. Let us prove that
 the symplectic structures $\omega_H$ and $\omega_K$ coincide, up to a
 sign.

\begin{thm} In the setting of isomonodromic deformations of the Fuchsian system
(\ref{N1}) the two symplectic forms (\ref{omegah}) and (\ref{omegak}) 
coincide (up to a sign)
\begin{equation}\label{om-k-h}
\omega_K=-\omega_H.
\end{equation}
\label{krichit}\end{thm}

\noindent{\bf Proof.} The variations $\delta{A}$ tangent to the orbit
are obtained by means of 
an infinitesimal gauge transform with $G=\ID+\phi$, i.e.
$$
\delta{A} = G AG^{-1} + {{\rm d}G\over{\rm d}z}G^{-1}-{A}=-[{A},\phi]+\ddz\phi +O(\phi^2).
$$
So, representing the tangent vectors (\ref{tangent}) as
$$
\delta_i{A} =- [{A}, \phi_i] +\ddz \phi_i, \quad i=1, \, 2
$$
we obtain
$$
\delta_i {A}_k =- [{A}_k, U_k^{(i)}]
$$
where $U_k^{(i)}$ is given by the first term of the expansion of 
$\phi_{i}$ at $u_k$ 
$$
\phi_{i}({z})=U_k^{(i)}+{\cal O}({z}-u_k).
$$
Because of this
$$
\omega_H(\delta_1{A}, \delta_2{A}) = - \sum_k {\rm Tr} (
{A}_k [U_k^{(1)}, U_k^{(2)}]).
$$
In the formula (\ref{omegak})
we can take 
$$
\delta_i\Phi\,\Phi^{-1}=\phi_i, \quad i=1, \, 2.
$$ 
Indeed, the matrices $\phi_i$ and $\delta_i \Phi \, \Phi^{-1}$
satisfy the same equation. This follows from
$$
\ddz\delta\Phi=\delta{A}\,\Phi+
{A}\delta\Phi.
$$
Thus
\begin{eqnarray}
&&
\omega_K(\delta_1{A}, \delta_2{A})-{1\over2}\sum_k \res_{z=u_k}{\rm Tr}
(\delta_1 {A}\, \phi_2 -\phi_1 \delta_2 {A})
\nonumber\\
&&
= \sum_k {\rm Tr} ({A}_k \, [U_k^{(1)}, U_k^{(2)}]) 
-\omega_H(\delta_1{A}, \delta_2{A}).\nonumber
\end{eqnarray}
The Theorem is proved.\epf

One can also prove that the Hamiltonians (\ref{hamk}) correctly reproduce 
the formula (\ref{ham0}), up to a minus sign (this is not a problem because also
the signs of the two symplectic structures are opposite).
We shall use the Krichever formula in the next Section in order to compute
Poisson brackets for 
a new set of Darboux coordinates that we call {\it isomonodromic coordinates.}

\subsection{Multi-time dependent Hamiltonian systems and their canonical
transformations}

Let ${\cal P}$ be a manifold equipped with a Poisson bracket $\{~,~\}$. A
function $H=H(x;t)$ depending explicitly on time
defines a time dependent Hamiltonian system of the form
\beq\label{t-dep1}
\dot x = \{ x,H\}.
\eeq
The total energy
\beq\label{t-dep2}
E(t) :=H(x(t);t)
\eeq
computed on an arbitrary solution $x=x(t)$ to (\ref{t-dep1}) is not conserved.
However, the following well known identity describes its dependence on time
\beq\label{t-dep3}
\dot E = {\pal H(x;t)\over \pal t}|_{x=x(t)}.
\eeq
One can recast the equations (\ref{t-dep1}), (\ref{t-dep3}) into an (autonomous)
Hamiltonian
form using the following standard trick of introducing extended phase space
$$
\hat{\cal P} := {\cal P}\times {\mathbb R}^2_{t, E}
$$
with a Poisson bracket $\{~,~\}\hat{}$ such that
\beq\label{t-dep4}
\{~,~\}\hat{} _{_{\cal P}} = \{~,~\}, \quad \{x,t\}\hat{} =\{x,E\}\hat{}=0, \quad
\{E,t\}\hat{}=1
\eeq
The Hamiltonian
\beq\label{t-dep5}
\hat H:= H(x;t)-E
\eeq
yields the dynamics
(\ref{t-dep1}), (\ref{t-dep3}) along with the trivial equation
$$
\dot t =1.
$$
The new Hamiltonian $\hat H$ is a conserved quantity. One returns to the
original setting considering dynamics on the zero level surface
$$
\hat H(x; t,E) =0.
$$

Let us now consider $n$ functions on ${\cal P}$ depending on $n$ times
$H_1=H_1(x;{\bf t})$, \dots, $H_n=H_n(x;{\bf t})$ where
$$
{\bf t}:= (t_1, \dots, t_n).
$$
Assume that the time dependent Hamiltonian systems
\beq\label{m-t1}
{\pal x\over \pal t_i} = \{ x, H_i\}, \quad i=1, \dots, n
\eeq
commute pairwise, i.e.
\beq\label{m-t2}
{\pal\over \pal t_j} {\pal x\over \pal t_i} = 
{\pal\over \pal t_i} {\pal x\over \pal t_j} \quad {\rm for} ~ {\rm all}~i\neq j.
\eeq
Because of commutativity there exist common solutions $x=x({\bf t})$ of the
family (\ref{m-t1}) of differential equations.
We want to introduce an analogue of the extended phase space for these multi-time
dependent systems. We begin with the following simple

\begin{lm} The Hamiltonian systems (\ref{m-t1}) commute
{\rm iff} the functions
\beq\label{m-t3}
c_{ij}(x; {\bf t}):={\pal H_i\over \pal t_j} -{\pal H_j\over \pal t_i} + \{ H_i, H_j\} , \quad
i\neq j
\eeq
are Casimirs of the Poisson bracket, i.e.
$$
\{ x, c_{ij}\} =0.
$$
The energy functions
\beq\label{en}
E_i=E_i({\bf t}) :=  H_i(x({\bf t}); {\bf t}) , \quad i =1, \dots, n
\eeq
satisfy
\beq\label{m-t4}
{\pal E_i\over \pal t_j} = \left[{\pal H_j\over \pal t_i}+ c_{ij}(x; {\bf
t})\right]_{x=x({\bf t})}.
\eeq
In these equations 
\beq\label{pal}
{\pal H_i\over \pal t_j}:= {\pal H_i(x;{\bf t})\over \pal t_j}
\eeq
are partial derivatives.
\end{lm}

\noindent{\bf Proof.} Spelling out the left hand side of (\ref{m-t2}) gives
\eqa
&&
{\pal\over \pal t_j} {\pal x\over \pal t_i}= Lie_{\pal\over \pal t_j} \{ x,
H_i\}
\nn\\
&&
=\left\{ Lie_{\pal\over \pal t_j}x, H_i\right\} + \left\{ x, 
Lie_{\pal\over \pal t_j} H_i\right\}
\nn\\
&&
=\left\{ \{ x, H_j\}, H_i\right\} + \left\{ x, {\pal H_i\over \pal t_j} + \{
H_i, H_j\} \right\}.
\nn
\eeqa
In this calculation we have used that the Hamiltonian vector fields are
infinitesimal symmetries of the Poisson bracket. Substituting this expression
into (\ref{m-t2}) and using Jacobi identity we arrive at
$$
\left\{ x, {\pal H_i\over \pal t_j} -{\pal H_j\over \pal t_i} + \{ H_i,
H_j\}\right\} =0.
$$
This proves the first part of the Lemma. The second part easily follows from
(\ref{m-t3}):
$$
{\pal E_i\over \pal t_j} =\{ H_i, H_j\} +{\pal H_i\over \pal t_j} 
={\pal H_j\over \pal t_i} +c_{ij}.
$$
The Lemma is proved.\epf

\begin{df} The functions $H_1(x; {\bf t})$, \dots, $H_n(x; {\bf t})$ on 
${\cal P}\times {\mathbb R}^n$ define $n$ multi-time dependent commuting
Hamiltonian systems if they satisfy equations
\beq\label{m-t-def}
{\pal H_i\over \pal t_j} -{\pal H_j\over \pal t_i} + \{ H_i,
H_j\} =0, \quad i, \, j=1, \dots, n.
\eeq
\end{df}

The energy functions $E_i=E_i({\bf t})$ of a multi-time dependent commuting
family satisfy
\beq\label{m-t-en}
{\pal E_i\over \pal t_j} = {\pal H_j\over \pal t_i}{\large |}_{x=x({\bf t})}.
\eeq

\begin{rmk} The evolution equations for the energy functions take more
``natural'' form
\beq\label{m-t-en1}
{\pal E_i\over \pal t_j} = {\pal H_i\over \pal t_j}{\large |}_{x=x({\bf t})}.
\eeq
similar to (\ref{t-dep3}) under an additional assumption of commutativity of the
Hamiltonians 
$$
\{ H_i, H_j\} =0, \quad i\neq j
$$
as functions on the phase space ${\cal P}$.
Observe that the one-form
\beq\label{1-form}
\varpi = H_1(x({\bf t}))dt_1 + \dots + H_n(x({\bf t}))dt_n
\eeq
is closed for any solution $x({\bf t})$ if the Hamiltonians commute. 
The commutativity holds true in the case of Schlesinger equations (see below).
The closeness of the one-form (\ref{1-form}) is crucial in the definition of
the isomonodromic tau-function.
\end{rmk}

Let us now define the extended phase space
\beq\label{ex-sp}
\hat{\cal P} :={\cal P}\times {\mathbb R}^{2n}
\eeq
where the coordinates on the second factor will be denoted $t_1$, \dots, $t_n$,
$E_1$, \dots, $E_n$. The Poisson bracket $\{~,~\}\hat{}$ on the extended phase
space is defined in a way similar to (\ref{t-dep4})
\beq\label{pb-ext}
\{~,~\}\hat{} _{_{\cal P}} = \{~,~\}, \quad \{x,t_i\}\hat{} =\{x,E_i\}\hat{}=0, \quad
\{E_i,t_j\}\hat{}=\delta_{ij}.
\eeq 
The Hamiltonians on the extended phase space are given by
\beq\label{ham-ext}
\hat H_i = H_i(x;{\bf t}) - E_i, \quad i=1, \dots, n.
\eeq
On the extended phase space the multi-time dependent commuting Hamiltonian
equations can be put into a form of autonomous commuting Hamiltonian systems.
Namely, the following statement holds true.

\begin{lm} For a multi-time dependent commuting family of Hamiltonian systems
the Hamiltonians (\ref{ham-ext}) commute pairwise. The corresponding Hamiltonian
equations on the extended phase space (\ref{ex-sp}) read
\eqa
&&
{\pal x\over \pal t_j} =\{ x, \hat H_j\}\hat{} = \{ x, H_j\}
\nn\\
&&
{\pal E_i\over \pal t_j} =\{ E_i, \hat H_j\}\hat{} = {\pal H_j\over \pal t_i}
\nn\\
&&
{\pal t_i\over \pal t_j} =\{t_i, \hat H_j\}\hat{} = \delta_{ij}.
\label{eq-ext}
\eeqa
\end{lm}

\noindent Proof is straightforward.

On the common level surface
$$
\hat H_1=0, \dots, \hat H_n =0.
$$ 
in $\hat{\cal P}$ one recovers the original multi-time dependent dynamics.

Let us now consider the particular case of a symplectic phase space ${\cal P}$
of the dimension $2g$.
Introduce canonical Darboux coordinates $q_1$, \dots, $q_g$, $p_1$, \dots,
$p_g$, such that the symplectic form $\omega$ on ${\cal P}$ becomes
$$
\omega=\sum_{i=1}^g dp_i \wedge dq_i.
$$
Then the extended phase space carries a natural symplectic structure
\beq\label{sym-ext}
\hat\omega = \omega-\sum_{i=1}^n d E_i \wedge dt_i.
\eeq
The canonical transformations of a multi-time dependent commuting Hamiltonian
family are defined as symplectomorphisms of the extended phase space
$\hat{\cal P}$ equipped with the symplectic structure (\ref{sym-ext}).

\begin{exam} A multi-time dependent 
generating function $S=S({\bf q}, \tilde{\bf q}, {\bf t})$
satisfying
$$
\det {\pal^2 S\over \pal q_i \pal \tilde q_j}\neq 0
$$
defines a canonical transformation of the form
$$
\tilde p_i ={\pal S\over \pal \tilde q_i}, \quad p_i =-{\pal S\over \pal q_i},
\quad \tilde E_k = E_k-{\pal S\over \pal t_k}.
$$
Usually in textbooks the last equation is written as the transformation law
of the Hamiltonians, i.e., the new Hamiltonians $\tilde H_k$ are given by
$$
\tilde H_k = H_k-{\pal S\over \pal t_k}, \quad k=1, \dots, n.
$$
We will stick to this tradition.
\end{exam}

Let us come back to Schlesinger equations. In this case the position of the
poles $u_1$, \dots, $u_n$ play the role of the (complexified) time variables.
It is straightforward to prove that the Schlesinger equations on ${\mathfrak g}^*$ can be
considered as a multi-time dependent commuting Hamiltonian family.

\begin{thm} The Hamiltonians $H_k$ of the form (\ref{ham0}) on ${\mathfrak g}^*$
Poisson commute
$$
\{H_k,H_l\}=0,\quad\forall k,l=1,\dots n.
$$
They also satisfy
$$
{\pal H_k\over\pal u_l} = {\pal H_l\over\pal u_k}.
$$
\end{thm}

We end this Section with the following simple observations about the Hamiltonians
(\ref{ham0}). First, these Hamiltonians are not independent. Indeed, 
\beq\label{sum-ham}
\sum_{i=1}^n H_i =0.
\eeq
Therefore the solutions to the Schlesinger equations depend only on the
differences $u_i-u_j$. Moreover, 
\beq\label{sum-ham1}
\sum_{i=1}^n u_i H_i = \sum_{i<j} \tr \,A_iA_j =\frac12 \tr \left( A_\infty^2
-\sum_{i=1}^n A_i^2\right).
\eeq
So, for a solution to the Schlesinger equations
\beq\label{scal-inv}
\sum_{i=1}^n u_i{\pal\over\pal u_i} A_k =[A_k, A_\infty].
\eeq
Thus the Hamiltonian (\ref{sum-ham1}) generates trivial dynamics on the reduced
symplectic leaves. This implies that the solutions to the Schlesinger equations
depend only on $n-2$ combinations of the variables $u_1$, \dots, $u_n$
invariant w.r.t. the action of one-dimensional affine group
$$
u_i\mapsto a\, u_i +b, \quad i=1, \dots, n, \quad a\neq 0.
$$
Due to this invariance it is sometimes convenient to normalize
the position of the poles of the Fuchsian systems by
\beq\label{normalize}
u_1=0, \quad u_2=1.
\eeq



\setcounter{equation}{0}
\setcounter{theorem}{0}

\section{Scalar reductions of Fuchsian systems.}

In this section we establish a birational transformation, that we
call {\it scalar reduction}, between the space of all
$m\times m$ Fuchsian systems of the form (\ref{N1}) considered
modulo diagonal conjugations and the space
of {\it special 
Fuchsian differential equations}, that we describe in the next sub-section.

\subsection{Special Fuchsian differential equations}

Recall \cite{codd} that a scalar linear differential equation of order $m$ with
rational coefficients is called Fuchsian if it has only regular singularities.
Writing the differential equation in the form
\begin{equation}\label{ur}
y^{(m)} = a_1(z) y^{(m-1)} + \dots + a_m(z) y
\end{equation}
one spells out the condition of regularity of a point $z=z_0$ in the 
form of existence of the limits
$$
b_k(z_0):=-\lim_{z\to z_0}(z-z_0)^k a_k(z), \quad k=1, \dots, m.
$$
The infinite point $z=\infty$ is regular if there exist the limits
$$
b_k(\infty):=-\lim_{z\to\infty} z^k a_k(z), \quad k=1, \dots, m.
$$

All the solutions to the equation (\ref{ur}) are analytic at the points of
analyticity of the coefficients.
Let $z=z_0$ be a pole of the coefficients of the Fuchsian equation.
The {\it indicial equation} at the point $z=z_0$ reads
\begin{eqnarray}
&&
\lambda(\lambda-1) \cdots (\lambda-m+1) + 
b_1(z_0) \lambda(\lambda-1)\cdots(\lambda-m+2) + \dots \nn\\
&&
\qquad \dots+ b_{m-1}(z_0)\lambda +b_m(z_0)=0.\label{ind-0}
\end{eqnarray}
If the roots $\lambda_1=\lambda_1(z_0)$, \dots, $\lambda_m=\lambda_m(z_0)$ of 
this equation are non-resonant,
i.e. none of the differences $\lambda_i-\lambda_j$ is a positive integer,
then there exists a fundamental system of solutions of the form
\begin{equation}\label{fund-ur}
y_j(z) = (z-z_0)^{\lambda_j} \left(1+\sum_{l>0} a_{jl} (z-z_0)^l \right), \quad
j=1, \dots, m, \quad z\to z_0.
\end{equation}
Therefore the roots of the indicial equation coincide with the exponents at the
regular singularity $z=z_0$.
Similarly, the indicial equation at $z=\infty$ reads
\begin{eqnarray}
&&
\lambda(-\lambda-1) \cdots (-\lambda-m+1) + 
b_1(\infty) \lambda(-\lambda-1) \cdots (-\lambda-m+2) +\dots\nn\\
&&
\qquad\dots + b_{m-1}(\infty)\lambda - b_m(\infty)=0.\label{ind-inf}
\end{eqnarray}
The roots of this equation give the exponents $\lambda_i=\lambda_i(\infty)$
at $z=\infty$.
The corresponding fundamental system of solutions, in the non resonant case is
given by
\begin{equation}\label{fund-ur-inf}
y_j(z) =z^{-\lambda_j}\left( 1+\sum_{l>0} a_{jl} z^{-l }\right), \quad j=1, \dots,
m. \quad, z\to \infty.
\end{equation}

\begin{rmk} \label{can-mon} For a Fuchsian differential equation with
non resonant roots of (\ref{ind-inf}) there exists a unique canonical, up to a
permutation of the roots, basis of solutions  (\ref{fund-ur-inf}). Therefore there exists a unique 
canonical normalization of the monodromy matrices, for a given choice of a basis in the fundamental
group of the punctured plane.
\end{rmk}

If the Fuchsian equation has $N+1$ poles including infinity then the exponents
satisfy the following {\it Fuchs relation}
\begin{equation}\label{f-rel}
\sum_{{\rm all}~{\rm poles}~ z_0} \sum_{i=1}^m \lambda_i(z_0) (N-1){m(m-1)\over 2}.
\end{equation}

In the resonant case logarithmic terms are to be added. They can be obtained
by a method similar to the one described above for the case of Fuchsian systems.
See \cite{codd} for the details.

\begin{rmk} \label{big} If $\lambda$ is a root of the indicial equation (\ref{ind-0}) and
$\lambda + n$ is not a root for any positive integer $n$ then there exists a solution
$$
y(z) =(z-z_0)^\lambda \left( 1+\sum_{l>0} c_l (z-z_0)^l\right), \quad z\to z_0.
$$
\end{rmk}

\begin{df} A pole $z=z_0$ of the coefficients of the Fuchsian equation is called
{\it apparent singularity} if all solutions $y(z)$ are analytic at $z=z_0$.
\end{df}

A necessary condition for the pole $z=z_0$ to be an apparent
singularity is that all the roots $\lambda_1$, \dots, $\lambda_m$
of the indicial equation (\ref{ind-0}) must be
non-negative integers. Absence of logarithmic terms impose additional 
constraints onto the coefficients of the Fuchsian equation. 
Let $M$ be the maximum of the roots. Then the full set of 
constraints can be obtained by plugging into the equation the expansions 
(\ref{fund-ur}) truncated at the term of order $M$ 
and requiring compatibility of the resulting linear system for the coefficients
$a_{jl}$, $j=1, \dots, m$, $l=1$, \dots, $M$. 

\begin{df} A Fuchsian differential equation of order $m$ is called {\it special}
if it has $n+1$ regular singularities at the points $z=u_1$, \dots, $z=u_n$,
$z=\infty$ and also 
$g$
apparent singularities $z=q_1$, \dots, $z=q_g$, where $g$ is given by the
formula (\ref{rod}) with the indices $0$, $1$, \dots,
$m-2$, $m$.
\end{df}

Observe that, due to Fuchs relation \cite{AB} for a special Fuchsian
equation the sum of the indices
at the points $z=u_1$, \dots, $z=u_m$,
$z=\infty$ is equal to $m-1$. For this reason we will denote, as above
$$
\lambda_1^{(i)}, \dots, \lambda_m^{(i)}
$$
the indices
at $z=u_i$ and
$$
\lambda_1^{(\infty)}, \lambda_2^{(\infty)}+1, \dots, \lambda_m^{(\infty)}+1
$$
the indices at infinity. These numbers have zero sum (cf. (\ref{l-sum})
\begin{equation}\label{sviaz}
\sum_{i=1}^n \sum_{j=1}^m \lambda_j^{(i)} +\sum_{j=1}^m \lambda_j^{(\infty)} =0.
\end{equation}

We now describe in more details the behavior of the coefficients of the
special Fuchsian equation near an apparent singularity.

\begin{lm} \label{basisapp} Near an apparent singularity $z=q_i$ 
there exist $m$ linear independent solutions ${y}_1(z),\dots,{y}_m(z)$ 
to the special Fuchsian equation (\ref{ur}) having expansions at ${z}=q_i$ of the following 
form
\begin{eqnarray}\label{basa}
&&
{y}_1 = 1 +{\alpha_1^{(i)}\over(m-1)!}({z}-q_i)^{m-1}+
{\mathcal O}({z}-q_i)^{m+1},\nn\\
&&
{y}_2 =  ({z}-q_i) + {\alpha_2^{(i)}\over(m-1)!}({z}-q_i)^{m-1}+
{\mathcal O}({z}-q_i)^{m+1},\nn\\
&&
\dots\dots\dots \\
&&
{y}_{m-1} = {1\over(m-2)!}({z}-q_i)^{m-2} +
{\alpha_{m-1}^{(i)}\over(m-1)!}({z}-q_i)^{m-1}+
{\mathcal O}({z}-q_i)^{m+1},\nn\\
&&
{y}_m = {1\over m!} ({z}-q_i)^m + 
{\alpha_m^{(i)}\over(m+1)!}({z}-q_i)^{m+1}+
{\mathcal O}({z}-q_i)^{m+2},
\nn\end{eqnarray}
where $\alpha_1^{(i)}$, \dots, $\alpha_m^{(i)}$ are some constant coefficients.
\end{lm}

\noindent The proof is obvious.
\vskip 0.1 cm

Denote
$$
{\cal L}:= -{d^m\over dz^m} + a_1(z) {d^{m-1}\over dz^{m-1}} + \dots + a_m(z)
$$
the differential operator in the l.h.s. of (\ref{ur}).

\begin{lm} \label{krit} The point $z=q_i$ is an apparent singularity of the special Fuchsian
equation ${\cal L}y=0$ if and only if the following conditions are satisfied.

\noindent 1. The coefficients $a_1(z)$, $a_2(z)$, \dots, $a_n(z)$ have at most
simple poles at $z=q_i$ and
\begin{equation}\label{vychet}
\res_{z=q_i} a_1(z) =-1.
\end{equation}

\noindent 2. There exist coefficients $\alpha_1^{(i)}$, \dots, 
$\alpha_{m-1}^{(i)}$ such that, after the substitution of the expansions
(\ref{basa}) into the differential equation one obtains
\begin{equation}\label{basa1}
{\cal L} y_j=O(z-q_i), \quad j=1, \dots, m-1, \quad z\to q_i.
\end{equation}
\end{lm}

\noindent{\bf Proof}. Suppose that $z=q_i$ is an apparent singularity. The indicial equation (\ref{ind-0}) 
for $z_0=q_i$ by assumption
must have the roots $0$, $1$, \dots, $m-2$, $m$. Because of it $b_1(z_0)=-1$
and $b_k(z_0)=0$ for $k>1$. So
$$
a_k(z) ={c_k\over (z-q_i)^{k-1}} \left( 1+O(z-q_i)\right), 
\quad z\to q_i, \quad
k=2, \dots, m,
$$
for some constants $c_1,\dots,c_m$.
Let us prove that $c_3= \dots = c_m=0$ 
(only the non trivial case $m\geq 3$ is to be studied). 
Indeed substituting the solution $y_1(z)$
from (\ref{basa}) into the equation one obtains that the l.h.s. behaves as
$$
{\cal L}y_1 \sim {c_m\over (z-q_i)^{m-1}}, \quad z\to q_i.
$$
Hence $c_m=0$. Similarly, substituting $y_2(z)$ one proves that, for $m\geq 4$
$c_{m-1}=0$. Continuing this procedure we prove that all the poles are at most
simple.

Validity of (\ref{basa1}) means that the solutions $y_1(z)$, \dots,
$y_{m-1}(z)$ corresponding to the roots $0$, $1$, \dots, $m-2$ of the indicial
equation at $z=q_i$ contain no logarithmic terms up to the order $O(z-q_i)^m$.
As it was explained in the previous page, this implies absence of logarithmic terms also in the
higher orders since the order of resonance by assumption is equal to $m$.
It remains to observe that the holomorphic solution $y_m(z)$ corresponding to
the maximal root $m$ always exists (see Remark \ref{big}). The vice--versa is obvious.\epf

The main result of this Section is a coordinate description of the
space of all special Fuchsian equations with given indices. Denote
\begin{eqnarray}\label{mom}
&&
p_s:= -\alpha_{m-1}^{(s)}+\delta_1^{(s)}, \quad s=1, \dots g.
\nn\\
&&
\delta_1^{(s)}=\sum_{t\neq s} {1\over q_s-q_t} + \sum_{i=1}^n
{1\over q_s-u_i} \left[ \sigma_1^{(i)} -{m(m-1)\over
2}\right]
\\
&&
\sigma_1^{(i)} = \sum_{j=1}^m \lambda_j^{(i)}.
\nn
\end{eqnarray}

\begin{thm} \label{special} Any special Fuchsian equation of order $m$ with given indices
$\lambda_j^{(i)}$, 
$j=1$, \dots, $m$, $i=1$, \dots, $n$,
$\infty$ satisfying (\ref{sviaz})
must have the form
\begin{equation}\label{sfe}
y^{(m)}= a_1(z) y^{(m-1)} + \dots + a_m(z) y
\end{equation}
where the coefficients are given by
\begin{eqnarray}\label{papper}
&&
a_1(z) = \sum_{s=1}^g {1\over z-q_s} + \sum_{i=1}^n {1\over z-u_i} \left[
\sum_{j=1}^m \lambda_j^{(i)} - {m\, (m-1)\over 2}\right]
\nn\\
&&
a_k(z) = \left[-\sum_{s=1}^g {c_{m-k+1}^{(s)}R(q_s)^{k-1}\over z-q_s}
+ (-1)^{k-1} \sum_{i=1}^n {\beta_k^{(i)}\over z-u_i} [R'(u_i)]^{k-1}
\right.+\\
&&
\left.+\beta_k^{(\infty)} z^{k\, n - n - k} 
+P_{k\, n - n - k -1}(z)\right] {1\over R(z)^{k-1}},\nn
\\
&&
k=2, \dots, m.
\nn
\end{eqnarray}
Here $R(z)=\prod_{k=1}^n (z-u_k)$,  $P_{n-3}(z)$, $P_{2n-4}(z)$, 
\dots, $P_{m\, n-m-n-1}(z)$ are some polynomials labeled by their 
degrees, $c_1^{(s)}$, \dots, $c_{m-1}^{(s)}$ are some numbers.
The coefficients $\beta^{(k)}_1,\dots,\beta^{(k)}_m$, $k=1,\dots,n,\infty$ 
depend only on the indices. They are determined from the identities
\begin{eqnarray}\label{betas}
&&
\lambda (\lambda-1) \dots (\lambda-m+1) - \beta_1^{(i)} \lambda(\lambda-1) \dots
(\lambda -m+2) +
\nn\\
&&\qquad +\beta_2^{(i)} \lambda(\lambda-1) \dots (\lambda-m+3) - \dots
-\beta_m^{(i)}
= \prod_{j=1}^m (\lambda-\lambda_j^{(i)})
\\\nn\end{eqnarray}
\begin{eqnarray}\label{betas1}
&&
\lambda (\lambda+1) \dots (\lambda+m-1) -\left[ {m(m+1)\over 2} -1 +
\sum_{j=1}^m \lambda_j^{(\infty)}\right] 
\lambda(\lambda+1) \dots
(\lambda +m-2) -
\nn\\
&&
\qquad -\beta_2^{(\infty)} \lambda(\lambda+1) \dots (\lambda+m-3) + \dots
+(-1)^{m-1}\beta_m^{(\infty)}\\
&&
\qquad = (\lambda-\lambda_1^{(\infty)})\prod_{j=2}^m 
(\lambda-\lambda_j^{(\infty)}-1).
\nn
\end{eqnarray}
The coefficients of the polynomials $P_{n-3}(z)$, $P_{2n-4}(z)$,  \dots, $P_{m\, n-m-n-1}(z)$ 
and the parameters $c_1^{(s)}$, \dots, $c_{m-1}^{(s)}$ are rational functions of 
$q_1$, \dots, $q_g$, $p_1$, \dots, $p_g$ and $u_1,\dots,u_n$. 
\end{thm}

\noindent{\bf Proof of the Theorem.}  The ansatz (\ref{papper}) follows 
from the definition of a Fuchsian equation and from the first of the 
claims of Lemma \ref{krit}. The expressions (\ref{betas}), (\ref{betas1})  via indices 
is nothing but the spelling of the indicial
equations (\ref{ind-0}), (\ref{ind-inf}). Let us now use the second statement
of Lemma \ref{krit} in order to show that all the remaining coefficients are uniquely determined
by $q_s$ and $p_s$.

Denote $\delta_k^{(s)}$ the constant term in the Laurent expansion of $a_k(z)$
near the apparent singularity $z=q_s$, $s=1$, \dots, $g$
\begin{eqnarray}\label{del-k}
&&
a_1(z) = {1\over z-q_s} + \delta_1^{(s)} +O(z-q_s)
\nn\\
&&
a_k(z) = -{c_{m-k+1}^{(s)}\over z-q_s} + \delta_k^{(s)} +O(z-q_s), \quad
k=2, \dots, m.
\end{eqnarray}

\begin{lm} \label{lmur-ja} 1. The coefficients $c_j^{(s)}$ in (\ref{papper}) coincide with
$\alpha_j^{(s)}$ in (\ref{basa}).

\noindent 2. The equation (\ref{ur}) with coefficients given by (\ref{papper}) is special Fuchsian {\rm iff} the following equations hold valid for all $s=1$, \dots, $g$
\begin{eqnarray}\label{ur-ja}
&&
p_s \alpha_1^{(s)} + \delta_m^{(s)}=0
\nn\\
&&
p_s \alpha_2^{(s)}+ \delta_{m-1}^{(s)}-\alpha_1^{(s)}=0
\nn\\
&&
\dots \quad \dots \quad \dots
\nn\\
&&
p_s \alpha_{m-2}^{(s)} + \delta_3^{(s)}-\alpha_{m-3}^{(s)}=0
\nn\\
&&
p_s \alpha_{m-1}^{(s)} + \delta_2^{(s)} -\alpha_{m-2}^{(s)}=0,
\end{eqnarray}
where 
\begin{equation}\label{rho}
p_s:= -\alpha_{m-1}^{(s)}+\delta_1^{(s)}, \quad s=1, \dots, g.
\end{equation}
\end{lm}

\noindent{\bf Proof.} Substituting the solution $y_l(z)$ for $1\leq l \leq m-2$
into (\ref{ur}), by using (\ref{papper}), one obtains, modulo
terms of order $O(z-q_s)$ the nontrivial contributions only from the terms
$$
a_1(z) y_l^{(m-1)} + a_2(z) y_l^{(m-2)} + a_{m-l+1}(z) y_l^{(l-1)} +
a_{m-l+2}(z) y_l^{(l-2)}.
$$
According to Lemma \ref{krit} this expression must be of the order $O(z-q_s)$ for  $z\to q_s$.
Spelling this out gives
\begin{eqnarray}\label{krit1}
&&
\left[ {1\over z-q_s} +\delta_1^{(s)}\right]\, \alpha_l^{(s)}
+\left[ -{c_{m-1}^{(s)}\over z-q_s} +\delta_2^{(s)}\right] \, \alpha_l^{(s)}
(z-q_s) 
\nn\\
&&
+ \left[ -{c_{l}^{(s)}\over z-q_s} +\delta_{m-l+1}^{(s)}\right]
\left[1+\alpha_l^{(s)} {(z-q_s)^{m-l}\over (m-l)!}\right] 
+ \left[ -{c_{l-1}^{(s)}\over z-q_s} + \delta_{m-l+2}^{(s)}\right] \, (z-q_s)
\nn\\
&&
=O(z-q_s).
\end{eqnarray}
Expanding these equations for $l=1, \dots, m-2$ yields
$$
c_l^{(s)}=\alpha_l^{(s)}
$$
and also the first $m-2$ equations of (\ref{ur-ja}).  Analogously for $l=m-1$ the nontrivial 
contributions in  (\ref{papper}) arise only in the terms
$$
a_1(z) y_{m-1}^{(m-1)} + a_2(z) y_{m-1}^{(m-2)} +
a_{3}(z) y_{m-1}^{(m-3)}.
$$
Again, imposing that this expression must be of the order $O(z-q_s)$ , one obtains 
$$
c_{m-1}^{(s)}=\alpha_{m-1}^{(s)}
$$
and also the last equation of (\ref{ur-ja}).  

Due to Lemma \ref{krit} the equations (\ref{krit1})
for $l=1, \dots, m-1$
are necessary and sufficient for the points $z=q_s$ to be apparent singularities
of a special Fuchsian equation. The Lemma is proved. \epf

\noindent{\bf End of the proof of the Theorem \ref{special}.}  We have derived a system 
of linear equations (\ref{ur-ja}) for the parameters $\alpha_1^{(s)}$, 
\dots, $\alpha_{m-1}^{(s)}$, $s=1, \dots, g$ and for the coefficients 
of the polynomials $P_{n-3}(z)$, $P_{2n-4}(z)$, \dots, $P_{m\, n-m-n-1}(z)$. 
It is easy to see that the number of equations is equal to the number of 
unknowns. It remains to prove that the determinant of this linear
system is not an identical zero.

Let us first eliminate the $\alpha$'s.  To this end we need to spell out the terms $\delta^{(s)}_2,\dots,\delta^{(s)}_{m}$ of order zero in the expansions of $a_2,\dots,a_{m}$ at $z=q_s$:
$$
\delta^{(s)}_k=(k-1)\frac{R'(q_s)}{R(q_s)}\alpha_{m-k+1}^{(s)}-
\sum_{t\neq s}\frac{\alpha_{m-k+1}^{(t)}R(q_t)^{k-1}}{(q_s-q_t)R(q_s)^{k-1}}+
f_k(q_s)+\frac{P_{kn-n-k-1}(q_s)}{R(q_s)^{k-1}},
$$
where 
$$
f_k(z)={(-1)^{k-1} \sum_{i=1}^n {\beta_k^{(i)}\over z-u_i} [R'(u_i)]^{k-1}+
\beta_k^{(\infty)} z^{k\, n - n - k}  \over R(z)^{k-1}}.
$$
The rational functions $f_1(z)$, \dots, $f_m(z)$ depend only on the positions of the poles $u_i$ 
and the indices. Using these notations we rewrite the equations (\ref{ur-ja})  as follows
\begin{eqnarray}\label{ur-ja1}
&&
p_s \alpha_1^{(s)} +\frac{(m-1)R'(q_s)}{R(q_s)}\alpha_{1}^{(s)}-
\sum_{t\neq s} {\alpha_1^{(t)}R(q_t)^{m-1}\over(q_s-q_t)R(q_s)^{m-1}} +
f_m(q_s) +{P_{m\, n - m - n -1}(q_s)\over R(q_s)^{m-1}}=0,
\nn\\
&&
p_s \alpha_2^{(s)} +\frac{(m-2)R'(q_s)}{R(q_s)}\alpha_{2}^{(s)}-
\sum_{t\neq s} {\alpha_2^{(t)}R(q_t)^{m-2}\over(q_s-q_t)R(q_s)^{m-2}} 
+f_{m-1}(q_s) +{P_{m n - m - 2n}(q_s)\over R(q_s)^{m-2}}=\alpha_1^{(s)}
\nn\\
&&
\dots \quad \dots \quad \dots
\nn\\
&&
p_s \alpha_{m-1}^{(s)} +\frac{R'(q_s)}{R(q_s)}\alpha_{m-1}^{(s)}-
\sum_{t\neq s}\frac{\alpha_{m-1}^{(t)}R(q_t)}{(q_s-q_t)R(q_s)}+
f_{2}(q_s) +{P_{n-3}(q_s)\over R(q_s)}\alpha_{m-2}^{(s)}.
\end{eqnarray}
Let us introduce the following vector notations. Denote
$$
{\bf q} = (q_1, \dots, q_g), \quad {\bf p}=(p_1, \dots, p_g).
$$
For any function $f(z)$ introduce vector
$$
f({\bf q}):= (f(q_1), \dots, f(q_g)).
$$
Similar notations will be used for functions of ${\bf p}$.
For example,
$$
{\bf p}^2 = (p_1^2, \dots, p_g^2).
$$
We also introduce $g$-component vectors $\alpha_j$ with the coordinates 
$\alpha_j^{(s)}$ and $\delta_1$ with the coordinates $\delta_1^{(s)}$. 
The last ingredient will be the $g\times g$ matrices $M^{(l)}=(M^{(l)}_{ij})$, $l=1,\dots,m-1$ with the matrix entries
\begin{equation}\label{cal-mos}
M^{(l)}_{ij} =\left(p_i +(m-l)\frac{R'(q_i)}{R(q_i)}\right)\delta_{ij} -
{R(q_j)^{m-l}(1-\delta_{ij})\over R(q_i)^{m-l}(q_i -q_j)}.
\end{equation}
Using these notations we can rewrite the equations (\ref{ur-ja1}) as follows

\begin{eqnarray}\label{ur-ja2}
&&
\alpha_{m-2} =  M^{(m-1)} \left(\delta_1-\bf{p}\right) + f_{2}({\bf q}) 
+{P_{n-3}({\bf q})\over
R({\bf q})}
\nn\\
&&
\alpha_{m-3} = M^{(m-2)}\alpha_{m-2} +f_3({\bf q})  +{P_{2n-4}({\bf q})\over
R({\bf q})^2}
\nn\\
&&
\dots \quad \dots \quad \dots 
\nn\\
&&
\alpha_1=M^{(2)}\alpha_2+f_{m-1}({\bf q})+{P_{(m-1)(n-1)-n-1}({\bf q})\over R({\bf q})^{m-2}}\nn\\
&&
0=M^{(1)}\alpha_{1} +f_m({\bf q})  +{P_{m(n-1)-n-1}({\bf q})\over
R({\bf q})^{m-1}}.
\end{eqnarray}
Substituting the first equation into the second equation we obtain
$$
\alpha_{m-3} = M^{(m-2)}\left(M^{(m-1)} \alpha_{m-1} + f_{2}({\bf q}) 
+{P_{n-3}({\bf q})\over
R({\bf q})}\right)+f_3({\bf q})  +{P_{2n-4}({\bf q})\over
R({\bf q})^2}.
$$
Continuing this process we express all $\alpha$'s via the known functions
and the coefficients of the polynomials $P_{n-3}(z)$, $P_{2n-4}(z)$, \dots,
$P_{m\, n-m-n-1}(z)$. On the last step we arrive at a linear equation for these
coefficients:
\begin{eqnarray}\label{main-ur}
&&
\hat P_{m(n-1)-n-1}({\bf q}) + M^{(1)} \hat P_{(m-1)(n-1)-n-1}({\bf q})+ \dots +
 M^{(1)} M^{(2)}\cdots M^{(m-2)} \hat P_{n-3}({\bf q})\nn\\
&&
\qquad= M^{(1)} M^{(2)}\cdots M^{(m-1)} [{\bf p}-\delta_1] -
M^{(1)} M^{(2)}\cdots M^{(m-2)} f_2({\bf q}) -\dots -\nn\\
&&
\qquad-M^{(1)} f_{m-1}({\bf q})-f_{m}({\bf q})
\end{eqnarray}
where we denote
$$
\hat P_{k\, n - k -n -1}(z):= {P_{k\, n - k -n -1}(z)\over R(z)^{k-1}}, \quad
k=2, \dots, m.
$$
It remains to prove that the determinant of the linear operator 
in the left hand
side of this system does not identically vanish.

Indeed, this determinant is a polynomial in $p_1$, \dots, $p_g$. Let us
compute the terms of highest degree in these variables. It is easy to see that
those terms can be written down explicitly 
$$
{W_{m,n}(q_1, \dots, q_g, \hat p_1, \dots,\hat p_g )\over [R(q_1) \dots
R(q_g)]^{m-1}}
$$
where the polynomial $W_{m,n}$ in $2g$ variables is defined in (\ref{det}),
$$
\hat p_s := p_s R(q_s).
$$
Clearly it is not an identical zero. This proves the Theorem.\epf

\begin{cor} \label{corallo}The positions $q_1$, \dots, $q_g$ of the apparent singularities
along with the auxiliary parameters $p_1$, \dots, $p_g$ are coordinates
on a Zariski open subset in the space of all special Fuchsian equations
with given indices and given Fuchsian singularities $u_1,\dots,u_n,\infty$.
\end{cor}

Observe that, in terms of the special Fuchsian equation the auxiliary parameters
$p_i$ are defined by
\begin{equation}\label{def-r}
p_i = \res_{z=q_i} \left[a_2(z)+{1\over 2} a_1(z)^2\right], \quad i=1, \dots, g.
\end{equation}
They are related to $\rho_s:=\alpha_{m-1}^{(s)}$ by the shift $\delta_1^{(s)}$ (see (\ref{mom})) 
of the form
\begin{eqnarray}\label{can-shift}
&&
\rho_s = -p_s +{\partial S({\bf q}, {\bf u})\over \pal q_s}
\\
&&
S({\bf q}, {\bf u})={1\over 2} \log\prod_{s\neq t} + \log\prod_{i=1}^n
\prod_{s=1}^g (q_s-u_i)^{\sigma_1^{(i)} -{m(m-1)\over 2}}
\nn
\end{eqnarray}
We will see below that the coordinates $q_s$ and $p_s$ are canonically
conjugated variables for the Schlesinger equations. The shift (\ref{can-shift})
is a canonical transformation. So, the variables $q_s$ and $-\rho_s$ are also
canonically conjugated.

\begin{rmk} The linear system (\ref{main-ur}) to be solved in order to
reconstruct the special Fuchsian equation with given indices and poles
and $g$ given pairs $(q_i, p_i)$ is very similar to the linear equation
(\ref{ansatz1}) used in the Appendix below in order to reconstruct the spectral
curve (\ref{comden4}) starting from its behavior over $z=u_1$, \dots, $z=u_n$,  
$z=\infty$ and a given divisor $(\gamma_1, \mu_1)+ \dots +(\gamma_g,\mu_g)$
of the degree $g$. The essential
difference is that, in matrix notations in (\ref{ansatz1}) enter the powers
of the diagonal matrix $(\mu_1, \dots, \mu_g)$ while in (\ref{main-ur}) this is
to be replaced by the matrix $M$ of the form (\ref{cal-mos}). It is a surprise
that the matrix $M^{(l)}$ for $l=m$ coincides with the Lax matrix for the Calogero - Moser system
\cite{moser}! At the moment we do not have an explanation of this coincidence.
\end{rmk}

\subsection{Transformation of Fuchsian systems into special Fuchsian
differential equations.}

In this Section we will assign to a Fuchsian system (\ref{N1}) a special
Fuchsian differential equation of the form
\begin{equation}\label{scal}
y^{(m)}= \sum_{l=0}^{m-1} d_{l} ({z}) y^{(l)}.
\end{equation}
Note a change of notations with respect to (\ref{ur}):
$$
d_l(z) = a_{m-l}(z), \quad l=1, \dots, m.
$$

The reduction of a system of differential equations to a scalar equation is
given by the following wellknown classical construction.
Denote by $\phi_1$, $\phi_2$, \dots, $\phi_m$ the components of the vector
function $\Phi$. The $m$-th order linear
differential equation for the scalar function $y:=\phi_1$ can be written
in the determinant form
$$
\det\,\left( \begin{matrix} y && y_1 && \dots && y_m \\
y' && y'_1 && \dots && y_m' \\
\dots && \dots && \dots && \dots \\
y^{(m)} && y_1^{(m)} && \dots && y_m^{(m)}\end{matrix}\right)=0.
$$ 
Here $y_1$, \dots, $y_m$ is the first row of a fundamental matrix solution
for the system (\ref{N1}). Expanding the determinant one obtains the needed
differential equation in the form
\begin{equation}\label{1-ass}
W(z) y^{(m)} =W'(z) y^{(m-1)} +W_{m-2}(z) y^{(m-2)} + \dots +W_0(z) y
\end{equation}
where
\begin{equation}\label{wro}
W(z) = \det\, \left( \begin{matrix}  y_1 && \dots && y_m \\
 y'_1 && \dots && y_m'\\
 \dots && \dots && \dots \\
 y_1^{(m-1)} && \dots && y_m^{(m-1)}\end{matrix}\right)
\end{equation}
is the Wronskian of the functions $y_1(z)$, \dots, $y_m(z)$, the functions
$W_l$ are certain determinants with the rows constructed from $(y_1, \dots,
y_m)$ and their derivatives. Therefore
\beq\label{wro-der}
d_{m-1}(z) = {W'(z)\over W(z)}, \quad d_l(z)={W_l(z)\over W(z)}, \quad l\leq
m-2.
\eeq
It readily follows that the coefficients of the scalar differential equation can
have poles only at zeroes of the Wronskian and at the points $u_1$, \dots,
$u_n$. Let us call the Fuchsian equation(\ref{1-ass}) the {\it 1-associated}
with the Fuchsian system (\ref{N1}). In a similar way one can obtain for any
$j=1, \dots, m$
the
$j$-associated Fuchsian equation for the $j$-th component of $\Phi$.

We will now give a more precise description of the poles of the scalar
Fuchsian equation and the corresponding exponents. Let us begin with the poles
that come from zeroes of the Wronskian. They are so-called {\it apparent
singularities} of the Fuchsian differential equation (\ref{1-ass}). That is,
the coefficients of the differential equation have poles at zeroes of $W(z)$ but
all solutions are analytic at these poles. Let us first compute exponents
at the apparent singularities.

\begin{lm} Let $z=q$ be a zero of the Wronskian $W(z)$ of the multiplicity $k$.
Then, if $k\leq m-1$ then the exponents of solutions to the Fuchsian equation
are
\beq\label{expo-k}
0,\, 1, \, \dots, m-k-1, \, m-k+j_1, \, m-k+1+j_1+j_2, \dots, m-1+j_1+\dots+j_k
\eeq
where $j_1, \dots, j_k$ are nonnegative integers satisfying
\beq\label{sum-expo-k}
k\, j_1+(k-1) j_2 + \dots + j_k =k;
\eeq
if $k\geq m$ then the exponents are
\beq\label{expo-m}
0,\, 1+j_1, \, 2+j_1+j_2, \dots, m-1+j_1+\dots+j_{m-1}
\eeq
where $j_1, \dots, j_{m-1}$ are nonnegative integers satisfying
\beq\label{sum-expo-m}
(m-1)j_1+(m-2)j_2 + \dots + j_{m-1}=k.
\eeq
\end{lm}

Proof. All the exponents at an apparent singularity must be nonnegative integers
$0\leq n_1 \leq \dots \leq n_m$. The corresponding basis of solutions
must have the form 
$$
y_1(z) =(z-q)^{n_1} (1+O(z-q)), \dots, y_m(z)= (z-q)^{n_m} (1+O(z-q)).
$$
Therefore all the exponents are pairwise distinct: in the opposite case the
difference of two basic solutions would have a higher exponent. Besides,
necessarily $n_1=0$. Indeed, otherwise all elements of the first line of the
fundamental matrix $\Phi$ would vanish at $z=q$. This contradicts
to non-degeneracy of the determinant $\det \Phi(z)$ of the fundamental matrix
solution to the Fuchsian system.

Let us now spell out the indicial equation (\ref{ind-0}) at $z=q$,
$$
\lambda (\lambda-1) \dots (\lambda-m+1) - b_{1} 
\lambda (\lambda-1) \dots (\lambda-m+2)- \dots -b_{m}=0
$$
where
$$
b_{m-i}=\lim_{z\to q} (z-q)^{m-i} d_i(z), \quad i=0, 1, \dots, m-1.
$$
Because of (\ref{wro-der}) we always have
$$
b_{1}=k
$$
where $k$ is the multiplicity of $z=q$ as a zero of the Wronskian.
Besides, if $k\leq m-1$ then 
$$
b_{k+1}= \dots = b_m=0.
$$
So the indicial equation factorizes
\eqa
&&
\lambda (\lambda-1) \dots (\lambda-m+k+1)
\nn\\
&&
\times \left[ 
(\lambda-m+k) \dots (\lambda-m+1)-k \,(\lambda-m+k) \dots (\lambda-m+2)- \dots
-b_{k}\right] =0.
\nn
\eeqa
The sum of the $k$ roots of the second factor is equal to
$$
(m-k) + (m-k+1) + \dots + (m-1) + k.
$$
As these roots must be pairwise distinct positive integers 
different from zeroes
of $\lambda (\lambda-1) \dots (\lambda-m+k+1)$, they can be represented in the
form (\ref{expo-k}), (\ref{sum-expo-k}). 

Let us now consider the case of multiplicity $k\geq m$. Since $\lambda=0$ must
be a root of the indicial equation, one has $b_0=0$. Hence the indicial equation
reads
$$
\lambda\,\left[ (\lambda-1) \dots (\lambda-m+1) -k \,(\lambda-1) \dots
(\lambda-m+2)  -\dots -b_{m-1}\right]=0.
$$
Again, the $(m-1)$ roots of the second factor are pairwise distinct positive integers
with the sum equal to 
$$
1+2+\dots+(m-1) +k.
$$
So, they can be represented in the form (\ref{expo-m}), (\ref{sum-expo-m}).
The Lemma is proved.

\begin{cor} \label{wro-cor} If $z=q$ is a simple root of the Wronskian $W(z)$, then the
exponents of the Fuchsian differential equation (\ref{1-ass}) are 
$$
0, 1, \dots, m-2, \, m.
$$
\end{cor}

\begin{rmk} In \cite{KO1} Kimura and Okamoto claimed that, for an apparent
singularity of multiplicity $k$ the exponents are
$$
0,\, 1, \dots, m-2, \, m+k-1.
$$
We were unable to reproduce the proof of this statement for $k>1$.
\end{rmk}

\begin{thm}
Take a Fuchsian system of the form (\ref{N1}) with 
pairwise distinct nonresonant exponents at $\infty$
$\lambda^{(\infty)}_i$, $i=1,\dots,m$.  Denote
$\lambda^{(k)}_1,\dots,\lambda^{(k)}_{m}$ the exponents at $u_k$, 
$k=1,\dots,n$. 
Suppose that for some $j=1,\dots,m$ 
\begin{equation}\label{nonzero1}
\sum_{k=1}^n A_{k_{ji}} u_k\neq 0,\,\forall\, i=1,\dots,m,\, i\neq j.
\end{equation}
Then the scalar differential equation for the $j$-th component
$$
y:=\phi_j
$$
of the solution $\Phi$ of the deformed system (\ref{N1}) possesses the following properties.

1. The coefficients of the equation depend rationally on the matrix elements of ${A}_k$,
$k=1,\dots,n$. They can be represented in the form
$$
d_{l}({z}) ={f_{l}({z})\over\prod_{k=1}^{n}({z}-u_k)^{m-l}
\prod_{i=1}^{g}({z}-q_i)},
$$
where 
$$
g={(n-1)m(m-1)\over 2}-(m-1),
$$
the functions $f_{l}({z})$ are polynomials of degree 
$$
{\rm deg}\, f_l(z)=(n-1)(m-l)+g,
$$
and $q_1$, \dots, $q_g$ are zeroes of the Wronskian (\ref{wro}).
They are apparent singularities of the Fuchsian equation.

2. This differential equation has regular singularities at 
$u_1,\dots,u_{n}$ with the 
exponents $\lambda^{(k)}_i$, $i=1,\dots,m$, at $u_k$ and also at $\infty$
with the exponents
$\lambda^{(\infty)}_1,\lambda^{(\infty)}_2+1,\dots,\lambda^{(\infty)}_{m}+1$ . 

3. If the numbers 
$q_1$, \dots, $q_g$
are pairwise distinct then the
exponents at each apparent singularity are $0,1,\dots,m-2,m$. 
\label{skal-red}\end{thm}

\noindent{\bf Proof.} Without loss of generality we may assume $j=1$. The above
construction gives a Fuchsian equation with regular singularities at $u_1$, \dots,
$u_n$, $\infty$ and apparent singularities at the zeroes of the Wronskian.
Denote $\tilde \lambda_1^{(i)}$, \dots, $\tilde \lambda_m^{(i)}$
the exponents of the solutions to the scalar equation
at the point $z=u_i$ for every $i=1, \dots, n$. From the construction it immediately follows that
for $j=1, \dots, m$
\begin{equation}\label{ner-i}
\tilde \lambda_j^{(i)} =\lambda_j^{(i)}+n^{(i)}_j, \quad\hbox{for some}\quad 
n ^{(i)}_j\in{\mathbb N}, 
\end{equation}
after a suitable labelling of the exponents. 
Let us now consider the exponents
$\tilde \lambda_1^{(\infty)}$, \dots, $\tilde \lambda_m^{(\infty)}$. They
satisfy modified equalities 
\begin{equation}\label{ner-inf}
\tilde \lambda_1^{(\infty)} =\lambda_1^{(\infty)}, \quad
\tilde \lambda_j^{(\infty)}=\lambda_j^{(\infty)}+1+n_j^{(\infty)}, 
\quad\hbox{for some}\quad 
n ^{(\infty)}_j\in{\mathbb N}, 
\quad j=2, \dots, m.
\end{equation}
Indeed, because of diagonality of the matrix $A_\infty$ there exists a
fundamental matrix solution of the Fuchsian system of the form
$$
\Phi= \left( \ID +O\left({1\over z}\right)\right)z^{-A_\infty}z^{-R^{(\infty)}}.
$$
Looking at the first row of this matrix yields the above estimates.

We will prove below that, doing if necessary a small monodromy preserving
deformation the above equalities are satisfied with $n^{(i)}_j=0$ for all $i=1,\dots,n,\infty$, $j=1,\dots,m$. To this end
we now proceed to considering the apparent singularities.

We already know from Corollary \ref{wro-cor} that, if $z=q$ is 
a simple zero of the Wronskian then the exponents at the apparent singularity are
$0$, $1$, \dots, $m-2$, $m$.
Let us now prove that, under the assumption (\ref{nonzero1}) the Wronskian
has exactly $g$ zeroes (counted with their multiplicities).

Denote 
$$
R(z):=\prod_{k=1}^n(z-u_k).
$$
Differentiating the linear system
\begin{equation}
\phi_i'({z} )=\sum_{j=1}^m {A}_{ij}({z} ) \phi_j({z} )
\label{r1}\end{equation}
where 
$$
{A}_{ij}(z)=\sum_{k=1}^{n}{{A}_{k_{ij}}\over{z}-u_k}
$$ 
and using Leibnitz rule we have for each $l=1,\dots,m$
\begin{equation}
\phi_i^{(l)}({z} )=\sum_{j=1}^m\sum_{k=0}^{l-1}
\left(\begin{array}{c}l-1\\ k\\ \end{array}\right) 
{A}_{ij}^{(l-1-k)}({z} ) \phi_j^{(k)}({z} ).
\label{r2}\end{equation}
Denoting $y=\phi_1$, we rewrite the system in the form
\begin{equation}
y^{(l)}({z} )=\sum_{j=2}^m{\cal P}_{l+1,j}({z} ) \phi_j({z} ) +
\sum_{k=0}^{l-1}{\cal Q}_{l+1,l-k}({z} )y^{(l-1-k)}({z} )
\label{r3}\end{equation}
for $l=1,\dots,m$. The coefficients ${\cal P}_{l+1,j}({z} )$ and 
${\cal Q}_{l+1,l-k}({z} )$ are rational functions of ${z} $ 
such that $R({z} )^l {\cal P}_{l+1,j}({z})$ is a polynomial
in ${z}$ of degree ${n l-l -1}$ and
$R({z})^{k+1}{\cal Q}_{l+1,l-k}({z})$ is a polynomial
in ${z}$ of degree $(k+1)(n-1)$. They can be computed from the following
recursion relations
starting with
$$
{\cal Q}_{2,1}({z})={A}_{11}({z}), \quad
{\cal P}_{2,j}({z})={A}_{1j}({z}).
$$
The recursion reads
\begin{eqnarray}
&&
{\cal P}_{l+2,j}(z)={\cal P}'_{l+1,j}(z)+\sum_{s=2}^m 
{\cal P}_{l+1,s}(z){A}_{s j}(z)\qquad j=2,\dots,m,\nn\\
&&
{\cal Q}_{l+2,l+1-k}(z)={\cal Q}'_{l+1,l-(k-1)}(z)+{\cal Q}_{l+1,l-k}(z),
\qquad k=1,\dots,l-1,\nn\\
&&
\nn\\
&&
{\cal Q}_{l+2,l+1}(z)={\cal Q}_{l+1,l}(z),\nn\\
&&
{\cal Q}_{l+2,1}(z)={\cal Q}'_{l+1,1}(z)+\sum_{s=2}^m 
{\cal P}_{l+1,s}(z){A}_{s 1}(z).
\nn\end{eqnarray}
Observe that  
${\cal Q}_{l+1,l}(z)={A}_{11}(z)$ for all $l$.

The system (\ref{r3}) can be written in the form
\begin{equation}
{\cal P} \left(\begin{array}{c}
\phi_1\\\phi_2\\ \dots \\ \phi_m \\
\end{array}\right)={\cal Q}\left(\begin{array}{c}
y\\ y'\\ \dots \\ y^{(m-1)} \\
\end{array}\right)
\label{r5}\end{equation}
where
\begin{equation}
\begin{array}{c}
{\cal P}= \left(\begin{array}{cccc}
1&0&\dots&0\\
0&{\cal P}_{22}(z)&\dots&{\cal P}_{2m}(z)\\
\dots&\dots&\dots\\
0&{\cal P}_{m,2}(z)&\dots&{\cal P}_{m,m}(z)\\
\end{array}\right)\\
\\
{\cal Q}= \left(\begin{array}{ccccc}
1&0 &\dots &\dots & 0\\
-{\cal Q}_{21}(z) & 1 & 0 &\dots & 0\\
-{\cal Q}_{31}(z) & -{\cal Q}_{32}(z) & 1 &\dots & 0\\
\dots&\dots&\dots&\dots&\dots\\
-{\cal Q}_{m,1}(z) &\dots &\dots & -{\cal Q}_{m,m-1} (z)& 1 \\
\end{array}\right)\\
\end{array}
\label{r5.5}\end{equation}
Let us prove invertibility of the $m\times m$ matrix ${\cal P}$. 

\begin{lm} 
\begin{equation}
{\rm det}\left({\cal P}\right){\Delta({z} )\over R({z} )^{{m(m-1)\over 2}}}.
\label{detl}\end{equation}
where $\Delta({z})$ is a polynomial of degree 
$g={(n-1)m(m-1)\over 2}-(m-1)$ with leading coefficient
\begin{equation}\label{det1}
-\underset{\scriptstyle 2\leq i<j\leq m}
{\prod}\left(\lambda^{(\infty)}_i-\lambda^{(\infty)}_j\right)
\left[\prod_{j=2}^m\sum_{k=1}^{n}u_k{A}_{k_{1j}}
\right].
\end{equation}
\end{lm}

Proof. First of all we prove that, if 
$p_{l+1,j}{z}^{nl-l-1}$ is the leading term in 
$R(z)^l {\cal P}_{l+1,j}$ then the 
leading term in ${R}(z)^{l+1}{\cal P}_{l+2,j}$ is 
$p_{l+2,j}{z}^{n(l+1)-(l+1)-1}$ with 
$p_{l+2,j}=(-\lambda^{(\infty)}_j-(l+1))p_{l+1,j}$. In fact, from the above
recursion relations one obtains
\begin{eqnarray}
&&
R(z)^{l+1}{\cal P}_{l+2,j}(z)=R(z)^{l+1}{\cal P}'_{l+1,j}(z)+
R(z)^{l+1}\sum_{s=2}^m {\cal P}_{l+1,s}(z){A}_{s j}(z)\sim\nn\\
&&
\qquad\sim((n l-l -1)p_{l+1,j} - n l p_{l+1,j}){z}^{nl-l+n}+
R(z)^{l+1}{\cal P}_{l+1,j}(z){A}_{jj}(z)\sim\nn\\
&&
\qquad\sim\left(\sum_{k=1}^{n}{A}_{k_{jj}}(z)-(l+1)\right) 
p_{l+1,j} {z} ^{n(l+1)-(l+1)-1}.
\nn\end{eqnarray}
Thus we have 
$$
p_{l+1,j}=(-\lambda^{(\infty)}_j-l)p_{l,j}\dots= (-\lambda^{(\infty)}_j-l)\cdots (-\lambda^{(\infty)}_j-2)
p_{2,j},
$$
with $p_{2,j}=\sum_{k=1}^{n}{A}_{k_{1j}} u_k$ 
because ${\cal P}_{2,j}(z)={A}_{1j}(z)$. Substituting these leading 
terms in the entries ${\cal P}_{ij}(z)$ we obtain that the leading term of
${\cal P}$ is ${\bf R}{\bf M}{\bf P}$ where
$$
{\bf P}={\rm diagonal}\left(1,p_{2,2},\dots,p_{2,m}\right)
$$

$$
{\bf R}={\rm diagonal}\left(1,{{z}^{n-2}\over R(z)},
\dots,{{z}^{(n-1)(i-1)-1}\over R(z)^{i-1}},\dots,
{{z}^{(n-1)(m-1)-1}\over R(z)^{m-1}}\right),
$$

$$
{\bf M}=\left(\begin{array}{cccc}
1&0&\dots &0\\0&1&\dots &1\\
0&(-\lambda^{(\infty)}_2-2) &\dots &(-\lambda^{(\infty)}_m-2)\\
\dots &\dots &\dots &\dots\\
0&\prod_{j=2}^{m-1}(-\lambda^{(\infty)}_2-j)&
\dots &\prod_{j=2}^{m-1}(-\lambda^{(\infty)}_m-j)\\
\end{array}\right),
$$
and computing the determinant we 
obtain (\ref{detl}), (\ref{det1}). The Lemma is proved. \epf

Observe that the leading coefficient of the polynomial $\Delta(z)$ 
cannot be zero, thanks to our hypothesis (\ref{nonzero1}).

Continuing the procedure used in the proof of the Lemma it is easy to obtain explicit 
formulae for the coefficients of the scalar equation.
Let ${\cal I}$ be the inverse matrix of ${\cal P}$. Its leading term equals
$$
{\cal I}_{ij}(z)\sim\left({\bf M}^{-1}\right)_{ij} 
{R({z})^{j-1}\over p_{2,i}{z}^{(n-1)(j-1)-1}}
$$
that is,
$$
{\cal I}_{ij}(z)={D_{ij}({z} )R({z})^{j-1}
\over\Delta({z})}
$$
where $D_{ij}({z} )$ is a polynomial in ${z}$ of degree 
${(n-1)m(m-1)\over 2}-(m-1)-((j-1)(n-1)-1)$ with coefficients depending 
on $i,j$. Solving the system (\ref{r5}) we obtain for $i>1$
\begin{equation}
\phi_{i}(z) = \sum_{j=1}^{m}
\left(-\sum_{2\leq s<j}{\cal I}_{is}(z){\cal Q}_{sj}(z)+{\cal I}_{ij}(z)\right)
y^{(j-1)}.
\label{r6}\end{equation}
Substituting (\ref{r6}) in (\ref{r2}) with $l=m$, we obtain 
$$
y^{(m)}(z)=\sum_{j=1}^m\left[\sum_{i=2}^m{\cal P}_{m+1,i}(z)\left(
{\cal I}_{ij}(z)-\sum_{s=2}^{j-1}{\cal I}_{is}(z){\cal Q}_{sj}(z)\right)+
{\cal Q}_{m+1,j}(z)\right] y^{(j-1)}(z)
$$
that is the requested differential equation 
$$
y^{(m)}(z)= \sum_{l=0}^{m-1} d_{l} ({z}) y^{(l)}(z),
$$
with 
$$
d_{l}=\sum_{i=2}^m{\cal P}_{m+1,i}(z)\left(
{\cal I}_{i,l+1}(z)-\sum_{s=2}^{l+1}{\cal I}_{is}(z)
{\cal Q}_{s,l+1}(z)\right)+{\cal Q}_{m+1,l+1}(z){f_{1,l}({z})\over
\Delta({z}){R}({z})^{m-l}},
$$
where $f_{1,l}({z})$ are polynomials of degree $(n-1)(m-l)+g$.

It remains to prove the statement about the exponents of the scalar Fuchsian
equation at the poles $z=u_i$, $z=\infty$. Let us assume for simplicity
that all the apparent singularities are pairwise distinct (the general case
can be considered in a similar way). Taking the sum of the equalities
(\ref{ner-i}), (\ref{ner-inf}) we obtain the following estimate
$$
\sum_{j=1}^m\left( \sum_{i=1}^n \tilde\lambda_j^{(i)} +\tilde
\lambda_j^{(\infty)}\right) =m-1+\sum_{i,j} n_j^{(i)}.
$$
The sum of exponents at the apparent singularity $z=q_s$ is equal to
$$
1+2+\dots+(m-2)+m={m(m-1)\over 2} + m.
$$
So the total sum of exponents of the Fuchsian equation satisfies
$$
\sum \, {\rm exponents} ~= m-1 + g \left( {m(m-1)\over 2 }+ m\right)+\sum_{i,j} n_j^{(i)}(g+n-1) {m(m-1)\over 2}+\sum_{i,j} n_j^{(i)}.
$$
But, according to the Fuchs relation (\ref{f-rel}) the total sum of exponents
over all $g+n+1$ regular singularities must be equal just to $(g+n-1) {m(m-1)\over 2}$. Therefore all non-negative integers $n^{(i)}_j$ appearing in the equalities  
(\ref{ner-i}), (\ref{ner-inf}) must be zero. 
The Theorem is proved.\epf
\vskip 0.2 cm

\begin{rmk} The above construction of the Fuchsian equation (\ref{scal}) for a given 
Fuchsian
system is clearly invariant w.r.t. simultaneous diagonal conjugations
of the coefficients of the latter.
\end{rmk}

To make sure that (\ref{scal}) is a special Fuchsian
equation for a generic Fuchsian system, we are to prove that, in the 
generic case all the roots of the polynomial $\Delta(z)$ are pairwise distinct. This will
follow from Theorem \ref{special} claiming that, in the space of all 
special Fuchsian equations,
the positions of the apparent singularities are independent variables and from the
result of the next Section that says that, under certain genericity assumption
the Fuchsian system can be reconstructed from the special Fuchsian equation
uniquely up to a conjugation by constant diagonal matrices.
  



\subsection{Inverse transformation.}

\begin{thm} Consider an $m$-th order special Fuchsian 
equation of the form
\begin{equation}
y^{(m)}({z})= \sum_{l=0}^{m-1} {f_{l}({z})\over
\Delta({z}){ R}({z})^{m-l}} y^{(l)}({z}),
\label{eqscal}\end{equation}
where
${R}({z})=\prod_{k=1}^{n}({z}-u_k)$, $\Delta({z})=\prod_{i=1}^{g}({z}-q_i)$, 
$g={(n-1)m(m-1)\over 2}-(m-1)$
and $f_{l}({z})$ are polynomials of degree $(n-1)(m-l)+g$. Let the 
exponents of the pole $u_k$, $k=1,\dots,n$,
be $\lambda^{(k)}_i$, $i=1,\dots,m$, and the ones of $\infty$ be
$\lambda^{(\infty)}_1,\lambda^{(\infty)}_2+1,\dots,\lambda^{(\infty)}_{m}+1$ 
and let $q_1,\dots,q_g$ be pairwise distinct apparent singularities 
of exponents 
$0,1,\dots,m-2,m$. If the monodromy group of the Fuchsian equation
(\ref{eqscal}) is irreducible, then there exists 
a $m\times m$ Fuchsian system of the form 
$$
\ddz\Phi=\sum_{k=1}^{n}{{A}_k\over{z}-u_k}\Phi
$$
with exponents 
$\lambda^{(\infty)}_1,\dots,\lambda^{(\infty)}_{m}$ at $\infty$ and
$\lambda^{(k)}_1,\dots,\lambda^{(k)}_{m}$ at $u_k$, and no apparent 
singularities, such that the first row of its fundamental matrix 
satisfies the given $m$-th order Fuchsian equation. The 
matrix entries of the matrices ${A}_k$, $k=1,\dots,n$,
depend rationally on the coefficients of the polynomials $f_{l}$
and on $q_1,\dots,q_g$, $u_1,\dots,u_{n}$.
Moreover, if 
$$
\ddz\tilde\Phi=\sum_{k=1}^{n}{\tilde{A}_k\over{z}-u_k}\tilde\Phi
$$
is another Fuchsian system corresponding to the given special 
Fuchsian equation, then there exists a diagonal matrix $D$ such that
$$
\tilde A_k = D^{-1} A_k \, D, \quad k=1, \dots, n.
$$
\label{eqtosystem}\end{thm}

\noindent{\bf Proof.} 
This proof follows essentially the proof due to Bolibruch of reconstruction 
of a Fuchsian system from a given Fuchsian equation \cite{AB,Bol2}. We 
need some extra machinery in our case to eliminate the apparent 
singularities $q_1,\dots,q_g$.  

\begin{lm}
The system 
$$
{{\rm d}Y\over{\rm d}{z}}=F({z}) Y
$$
constructed from (\ref{eqscal}) by the transformation
$$
Y^j=[\Delta({z}) {R}({z})]^{j-1}
{{\rm d}^{j-1}y\over{\rm d}{z}^{j-1}},\qquad j=1,\dots,m,
$$
is Fuchsian at $ u_1,\dots, u_{n}, q_1,\dots,q_g$ with the same exponents 
of (\ref{eqscal}) and has a regular singularity at $\infty$.
\label{bol1}\end{lm}

\noindent The proof is straightforward and can be found in \cite{AB,Bol2}.

First of all we want to eliminate the apparent singularities 
$ q_1,\dots, q_g$.  By Lemma \ref{basisapp}, near the point $z=q_i$ we can 
choose a basis of solutions $y_1,\dots,y_m$ such that
\begin{eqnarray}
&& 
{{\rm d}^{l-1}y_k\over{\rm d}{z}^{l-1}}={1\over(k-l)!}({z}-q_i)^{k-l}+{\cal O}({z}-q_i)^{m-l},
\qquad l\leq k<m,
\nn\\
&&
{{\rm d}^{l-1}y_k\over{\rm d}{z}^{l-1}}
={\alpha^{(i)}_k\over(m-l)!}({z}- q_i)^{m-l}+
{\cal O}({z}- q_i)^{m+2-l},\qquad l>k, \, k<m,
\nn\\
&&
{{\rm d}^{l-1}y_m\over{\rm d}{z}^{l-1}}
={1\over(m-l+1)!}({z}- q_i)^{m-l+1}+{\cal O}({z}- q_i)^{m+2-l},
\qquad k=m.
\nn
\end{eqnarray}
To eliminate all apparent singularities $ q_1,\dots, q_g$, we apply the following 
gauge transformation
$$
\hat Y=\Gamma(z) \Delta(z)^{-M} Y,
$$
where $M=\diag\left(0,1,\dots,m-2,m\right)$ and $\Gamma(z)$ is a lower triangular matrix
with all diagonal elements equal to 1 and all off-diagonal elements equal to zero apart from
the last row which is given by
$$
\Gamma(z)_{ml} =-\frac{R(z)^{m-1}}{\Delta(z)} g_l(z), \quad l=1,\dots,m-1.
$$
where $g_l(z)$ is a degree $g$ polynomial in $\frac{1}{z}$ such that $g_l(q_i)=\alpha_l^{(i)}$ and
$g_l(z)\sim z^{-g}$ as $z\to \infty$.
Let us show that the new matrix $\hat Y$ is holomorphic and invertible at $z=q_i$. 

In fact near $z=q_i$, we have
$$
\Gamma^{(i)}_{ml} =- {\alpha^{(i)}_l R(q_i)^{m-1}\over (z-q_i) \Delta'(q_i)}
+{\cal O}(1),
$$
and $\Delta(z)^{-M} Y=\diag\left(1,\dots,1,{1\over \Delta(z)}\right) 
\diag\left(1,R(z),\dots,R(z)^{m-1}\right)G(z)$, where
$$
\begin{array}{l}
G(z)_{lk} = {\cal O}(z-q_i),\quad l\neq k, m,\\
G(z)_{ll} =1+{\cal O}(z-q_i),\quad l\neq m\\
G(z)_{mk}=\alpha_k^{(i)} +{\cal O}(z-q_i),\quad k\neq m,\\
G(z)_{mm} = (z-q_i)+{\cal O}(z-q_i)^2.
\end{array}
$$
This gives
$$
\hat Y(z)= \hat Y_0 +{\cal O}(z-q_i),
$$
where ${\rm det}(\hat Y_0)= {R(q_i)^{\frac{m(m-1)}{2}}\over \Delta'(q_i)}\neq0$, as we wanted to prove.

We now need to study infinity. First, in the non-resonant case, we can 
choose a basis $y_1,\dots,y_m$  of solutions for the differential equation 
of the form
$$
\begin{array}{l}
y_1(z)= a_1 z^{-\lambda_1^{\infty}} \left(1+{\cal O}({1\over z})\right)\\
y_k(z)= a_k z^{-\lambda_k^{\infty}-1} \left(1+{\cal O}({1\over z})\right),
\qquad k=2,\dots,m,\\
\end{array}
$$
for some arbitrary non-zero coefficients $a_1,\dots,a_m$.
As a consequence we obtain that
$$
\begin{array}{l}
Y_{l1}=(-1)^{l-1} a_1 \lambda^{(\infty)}_1(\lambda^{(\infty)}_1+1)\cdots
(\lambda^{(\infty)}_1+l-2)
z^{(g+n)(l-1)} z^{-\lambda_1^{\infty}-l+1}
\left(1+{\cal O}({1\over z})\right),\\
Y_{lk}=(-1)^{l-1} a_k (\lambda^{(\infty)}_k+1)(\lambda^{(\infty)}_k+2)\cdots
(\lambda^{(\infty)}_k+l-1)z^{(g+n)(l-1)} z^{-\lambda_k^{\infty}-l}
\left(1+{\cal O}({1\over z})\right),\\
\qquad k=2,\dots,m\\
\end{array}
$$
that gives
$$
Y(z) = z^C G^{(\infty)}(z) z^{-\Theta^{(\infty)}}
$$
where $\Theta^{(\infty)}={\rm diagonal}(\lambda_1^{\infty},
\lambda_2^{\infty}+1,\dots,\lambda_m^{\infty}+1)$, $C ={\rm diagonal}(
0,(g+n-1),\dots,(g+n-1)(m-1))$ and $G^{(\infty)}(z)$ is holomorphicly 
invertible at infinity such that $G^{(\infty)}(\infty)$ has all minors not 
equal to zero. In particular the arbitrary choice of the parameters
$a_1,\dots,a_m$ implies freedom of multiplication of $Y$ by a diagonal 
matrix, ${\rm diagonal}(a_1,\dots,a_m)$, from the right.

In the resonant case, in similar way one obtains
$$
Y(z) = z^C G^{(\infty)}(z) z^{-\Theta^{(\infty)}} {z}^{-R^{(\infty)}}.
$$
To estimate the indices at infinity of our $\hat Y({z})= \Gamma(z)
\Delta(z)^{-M} z^C G^{(\infty)}(z) z^{-\Theta^{(\infty)}} {z}^{-R^{(\infty)}}$ we want to use the 
following lemma proved in \cite{AB,Bol2} (see Lemma 4.1.2 in both references).

\begin{lm}
Let $U({z})$ be a matrix holomorphicly invertible at $\infty$ and 
let all the principal minors of $U(\infty)$ be non zero. Then for any integers 
$k_1\leq k_2\leq\dots\leq k_m$ there exists an lower triangular matrix 
$\Gamma^{(\infty)}(z)$ with elements on the principal diagonal 
equal to 1, $\Gamma^{(\infty)}(z)$ polynomial in $z$, and a 
matrix $V^{(\infty)}({z})$ holomorphicly invertible in a neighborhood of 
$\infty$ such that
$$
\Gamma^{(\infty)}(z)
z^{K} U({z}) 
=V^{(\infty)}({z}){z}^{K}.
$$
where $K=\diag\left(k_1,k_2,\dots,k_m\right)$.
\label{bol3}\end{lm}

We add that $V^{(\infty)}({z})$ and $\Gamma^{(i)}({z})$
depend only on the first $S:=k_m-k_1$ terms of the series expansion of 
$U({z})=\sum_{s=1}^\infty U_s z^{-s}$ near $\infty$. 

To apply the above lemma, we first observe that 
$$
\Gamma(z)z^{-M g+C} =z^{-M g+C} \tilde\Gamma(z),
$$
where $\tilde\Gamma(z)$ is a lower triangular matrix
with all diagonal elements equal to 1 and all off-diagonal elements equal to zero apart from
the last row which is given by
$$
\tilde\Gamma(z)_{ml} \sim 
z^{(n-1)(l-m)-g}, \qquad l=1,\dots,m-1.
$$

To  apply Lemma \ref{bol3} we need to introduce a permutation $P$ such that
$P z^{-M g+C} P^{-1}=z^K$ where, in the case $m\geq 3$
$$
K={\rm diagonal}\left((m-1)(n-1)-g,0,(n-1),2(n-1),\dots, (m-2)(n-1)\right),
$$ 
and in the case $m=2$, $K= {\rm diagonal}(0,1)$ and $P=\ID$. Moreover, in the case $m>2$,
 we need to show that $P\tilde\Gamma(\infty)  G^{(\infty)}(\infty) P^{-1} $ has all principal 
 minors different from zero.  This is a straightforward  consequence of the fact that 
$\tilde\Gamma(\infty) =\ID$ and $\tilde G=P G^{(\infty)}(\infty) P^{-1} $ is given by
$$
\begin{array}{l}
\tilde G_{11}=(-1)^{m-1} a_1\lambda_1^{(\infty)}(\lambda_1^{(\infty)}+1)\dots
(\lambda_1^{(\infty)}+m-2),\\
\tilde G_{1k} =(-1)^{m-1} a_k(\lambda_k^{(\infty)}+1)\dots
(\lambda_k^{(\infty)}+m-1), \quad k\neq 1\\
\tilde G_{l1}=(-1)^{l-1} a_1\lambda_1^{(\infty)}(\lambda_1^{(\infty)}+1)\dots
(\lambda_1^{(\infty)}+l-2),\quad l\neq 1,m,\\
\tilde G_{lk} =(-1)^{l-1} a_k(\lambda_k^{(\infty)}+1)\dots
(\lambda_k^{(\infty)}+l-1), \quad l\neq 1,m,\quad k\neq 1\\
\tilde G_{mk}=a_k,\quad k=1,\dots,m.\\
\end{array}
$$
We can then apply Lemma  \ref{bol3} to 
$$
P \hat Y= P z^{-M g+C}\tilde\Gamma(z) G^{(\infty)}(z) z^{-\Theta^{(\infty)}} {z}^{-R^{(\infty)}}=z^{K}
P \tilde\Gamma(z) G^{(\infty)}(z) P^{-1} P z^{-\Theta^{(\infty)}} {z}^{-R^{(\infty)}}.
$$
We obtain a gauge transformation with the matrix 
$\Gamma^{(\infty)}(z)$ polynomial in $z$,
such that the new fundamental matrix 
$$
\tilde Y= \Gamma^{(\infty)}(z) P \hat Y= \Gamma^{(\infty)}(z) z^K P
 \tilde\Gamma(z) 
G^{(\infty)}(z) P^{-1} P z^{-\Theta^{(\infty})}{z}^{-R^{(\infty)}},
$$
factors as
$$
\tilde Y({z})= V^{(\infty)}({z}) z^{K} P
{z}^{-\Theta^{(\infty)}} {z}^{-R^{(\infty)}}= V^{(\infty)}({z}) P^{-1}z^{-M g+C}
{z}^{-\Theta^{(\infty)}} {z}^{-R^{(\infty)}},
$$
with the matrix $V^{(\infty)}({z})$ holomorphicly invertible in a 
neighborhood of $\infty$.
The new exponents at $\infty$ are $\lambda^{(\infty)}_1$,
$\hat\lambda_m^{(\infty)}=\lambda^{(\infty)}_{m}+1+m g-(m-1)(n+g-1)$, 
and, for $j=2,\dots,m-1$,
$\hat\lambda^{(\infty)}_j=\lambda^{(\infty)}_j+1-(j-1)(n-1)$. 
Their sum is zero,
therefore $\infty$ is a Fuchsian singularity (see [AB]). 

So we have constructed a $m\times m$ Fuchsian system of the form 
$$
\ddz\tilde Y=\sum_{k=1}^{n}{\tilde{A}_k\over{z}-u_k}\tilde Y
$$
with exponents $\lambda^{(\infty)}_1$, $\hat\lambda_2^{(\infty)},\dots,
\hat\lambda^{(\infty)}_m$ at $\infty$ and
$\lambda^{(k)}_1,\dots,\lambda^{(k)}_{m}$ at $u_k$, and no apparent 
singularities. We now want to map this system to
a $m\times m$ Fuchsian system 
with exponents 
$\lambda^{(\infty)}_1,\dots,\lambda^{(\infty)}_{m}$ at $\infty$ and
$\lambda^{(k)}_1,\dots,\lambda^{(k)}_{m}$ at $u_k$. We need the following:

\begin{lm}
Given a Fuchsian system of the form (\ref{N1}). Let 
$\lambda^{(k)}_1,\dots \lambda^{(k)}_m$ be the eigenvalues of the matrix 
${A}_k$ for $k=1,\dots,n,\infty$ and let ${\cal G}_k$ be its
diagonalizing matrix,
$$
{\cal G}_k^{-1} A_k {\cal G}_k ={\rm diag}\, 
( \lambda^{(k)}_1,\dots \lambda^{(k)}_m ).
$$
Assume that there are
two eigenvalues, say $\lambda^{(k)}_1$ and $\lambda^{(k)}_m$ such that
$\lambda^{(k)}_m\neq \lambda^{(k)}_1+ 1$, $\lambda^{(k)}_m\neq 
\lambda^{(k)}_1+ 2$, 
and for all $l\neq 1,m$, $\lambda^{(k)}_{1}\neq \lambda^{(k)}_l-1$ and
$\lambda^{(k)}_m\neq \lambda^{(k)}_l+1$. If not all entries in position $m1$
of the matrices ${\cal G}_k^{-1} {A}_l {\cal G}_k$,
$l=1,\dots,n$, are zero, then there exists a gauge transformation 
$G_k({z};{A}_1\dots,{A}_{n},u_1,\dots,u_{n})$,
rational in all arguments, such that the new matrices 
$\tilde {A}_l$, $l=1,\dots,n$, $l\neq k$ have the same 
eigenvalues as the old ones ${A}_l$ and the new matrix
$\tilde {A}_k$ has eigenvalues 
$\lambda^{(k)}_1+1,\lambda^{(k)}_2\dots,\lambda^{(k)}_{m-1},
\lambda^{(k)}_m-1$.
Moreover the gauge transformation $G_k({z};{A}_1\dots,{A}_{n},u_1,\dots,u_{n})$ preserves the Schlesinger equations.
\label{lmb0}\end{lm}

\noindent{\bf Proof.} We give here the gauge transformation $G_\infty(z)$ 
giving
rise to the change $\lambda^{(\infty)}_1\to\lambda^{(\infty)}_1+1$,
$\lambda^{(\infty)}_m\to\lambda^{(\infty)}_m-1$. So we assume that $\lambda^{(\infty)}_1$ and $\lambda^{(\infty)}_m$ are such that
$\lambda^{(\infty)}_m\neq \lambda^{(\infty)}_1+ 1$, $\lambda^{(\infty)}_m\neq 
\lambda^{(\infty)}_1+ 2$, 
and for all $l\neq 1,m$, $\lambda^{(\infty)}_{1}\neq \lambda^{(\infty)}_l-1$ and
$\lambda^{(\infty)}_m\neq \lambda^{(\infty)}_l+1$, and if not all entries in position $m1$
of the matrices ${A}_l $, $l=1,\dots,n$, are zero.

Let us fix a fundamental matrix $\Phi$ normalized at infinity
$$
\Phi_\infty= \left(\ID +\frac{\Psi_1}{z} +\frac{\Psi_2}{z^2}+
{\mathcal O}\left(\frac{1}{z^3}\right)\right) z^{- A_\infty} z^{-R^{(\infty)}},
$$
where
$$
\Lambda=A_\infty, \quad R^{(\infty)}=R_1 + R_2 +\dots,
$$
\begin{eqnarray}\label{primer}
&&
\left(R_1\right)_{i\, j} =\left\{\begin{array}{cc} \left( B_1\right)_{i\, j}, & \lambda_i=\lambda_j+1\\ 0, & \mbox{\rm otherwise}\end{array}\right.
\nn\\
&&
\nn\\
&&
B_1 =-\sum_k A_k u_k
\nn\\
&&
\left( \Psi_1\right)_{i\, j} =\left\{ \begin{array}{cc} -\frac{\left( B_1\right)_{i\, j}}{\lambda_i-\lambda_j-1}, & 
\lambda_i\neq \lambda_j+1\\  & \\
\mbox{\rm arbitrary}, & \mbox{\rm otherwise}\end{array}\right.
\nn\\
&&
\\
&&
\left(R_2\right)_{i\, j} =\left\{\begin{array}{cc} \left( B_2 -\Psi_1 R_1 +B_1\Psi_1\right)_{i\, j} , & \lambda_i=\lambda_j+2\\ 0, & \mbox{\rm otherwise}\end{array}\right.
\nn\\
&&
\nn\\
&&
B_2 =-\sum_k A_k u_k^2\nn\\
&&
\left( \Psi_2\right)_{i\, j} =\left\{ \begin{array}{cc} \frac{\left( -B_2+\Psi_1 R_1 -B_1\Psi_1\right)_{i\, j}}{\lambda_i-\lambda_j-2}, & 
\lambda_i\neq \lambda_j+2\\  & \\
\mbox{\rm arbitrary}, & \mbox{\rm otherwise}\end{array}\right.
\nn
\end{eqnarray}
Consider the following gauge transformation $\Phi(z) = (I(z)+G)\tilde\Phi(z)$  where 
$$
I(z):={\rm Diagonal}\left(z,0,\dots,0\right),
$$
and
\begin{eqnarray}\label{newgauge}
&&
G_{m1}= \Psi_{1_{m1}},\qquad G_{{1m}}=-\frac{1}{G_{m1}},\nn\\
&&
\hbox{if}\quad p\neq 1,m,\quad
G_{{pp}}=1,\quad
G_{{1p}}=\Psi_{1_{mp}} G_{1m},\qquad
G_{{p1}}=\Psi_{1_{p1}},\\
&&
\hbox{if}\quad p,q \neq 1,\, p\neq q,\quad G_{{pq}}=0,\nn\\
&&
G_{{11}}= G_{1m} \Psi_{2_{m1}}+\Psi_{1_{11}}, \quad\hbox{and}\qquad G_{mm}=0.
\nn\end{eqnarray}
In order to see that this gauge transformation is always well defined  is enough to observe that $\Psi_{1_{m1}}(u)$ is never identically equal to zero if at least one of the $(m, 1)$ matrix entries of the matrices $A_1(u)$, \dots, $A_n(u)$ is different from identical zero. Indeed, this follows from the linear equations
\begin{eqnarray}
&&
\label{ddpsi1}
\partial_i\Psi_1 =-A_i
\\
&&\label{ddpsi2}
\partial_i \Psi_2 = -A_i\Psi_1 -u_i A_i,
\end{eqnarray}
which are a straightforward consequence of the equation
$$
\partial_i \Phi_\infty(z;u)=-\frac{A_i}{z-u_i}\, \Phi_\infty(z;u), \quad i=1, \dots, n,
$$
describing the  $u$-dependence of the fundamental matrix $\Phi_\infty(z;u)$. 

Let us prove that this transformation maps the matrices $A_1,\dots,A_n$ to new matrices 
$\tilde A_1,\dots,\tilde A_n$ given by
$$
\tilde A_k:= (I(u_k)+G)^{-1} A_k (I(u_k)+G),
$$
such that
\begin{equation}\label{con-gauge}
\tilde A_\infty=-\sum_{k=1}^n \tilde A_k={\rm diagonal}\left(
\lambda^{(\infty)}_1+1,\lambda^{(\infty)}_2,\dots,
\lambda^{(\infty)}_{m-1},\lambda^{(\infty)}_m-1\right).
\end{equation}
In fact $(I(z)+G)^{-1}=J(z)+G^{-1}$ where
$$
J(z):={\rm Diagonal}\left(0,\dots,0,z\right),
$$
therefore
$$
\tilde A_k:=G^{-1} A_k I(u_k) + G^{-1} A_k G + J(u_k) A_k I(u_k)+J(u_k) A_k G.
$$
Multiplying by $G$ from the left and summing on all $k$ we get that the condition (\ref{con-gauge}) is satisfied if and only if
\begin{eqnarray}\label{gauge1}\nn
&&
\left(\begin{array}{ccc}
-g_{11}&(\lambda^{(\infty)}_1-\lambda^{(\infty)}_2) g_{12}&\dots\\
(\lambda^{(\infty)}_2-\lambda^{(\infty)}_1-1) g_{21}&0&\dots\\
\dots&0&\dots\\
(\lambda^{(\infty)}_m-\lambda^{(\infty)}_1-1) g_{m1}&0&\dots\\
\end{array}
 \right.\qquad\qquad\nn\\
 &&
\qquad\qquad\qquad\left.\begin{array}{ccc}
\dots&(\lambda^{(\infty)}_1-\lambda^{(\infty)}_{m-1}) g_{1\,m-1}&
(\lambda^{(\infty)}_1-\lambda^{(\infty)}_m+1) g_{1m}\\
\dots&\dots&0\\
\dots &\dots&0\\
\dots &\dots&0\\
 \end{array}
 \right)=\nn\\
 &&\nn\\
&&\nn
=\left(\begin{array}{cccc}
\sum_{k} A_{k_{11}}u_k&0&\dots&0\\
\dots &0&\dots&0\\
\sum_{k} A_{k_{m1}}u_k&0&\dots&0\\
\end{array}
 \right)+
 \left(\begin{array}{cccc}
g_{1m} \sum_{k} A_{k_{m1}}u_k^2&0&\dots&0\\
0 &0&\dots&0\\
\dots&0&\dots&0\\
\end{array}
 \right)+\\
 &&
 \\
 &&\nn
\qquad\quad +\left(\begin{array}{cccc}
g_{1m}\sum_{s}\sum_{k} A_{k_{ms}}u_k  g_{s1}&\dots&
g_{1m}\sum_{s}\sum_{k} A_{k_{ms}}u_k  g_{sm}\\
0&\dots&0\\
\dots &\dots&\dots\\0&\dots&0
\end{array}
 \right).
 \end{eqnarray}
Observe that in our assumptions on the eigenvalues $\lambda^{(\infty)}_1$ and $\lambda^{(\infty)}_m$,
 these formulae are clearly satisfied thanks to the fact that $\Psi_1$, $\Psi_2$ and $R^{(\infty)}$ are given by formulae (\ref{primer}). 

Let us prove that this gauge transformation preserves the Schlesinger equations. Differentiating $\tilde A_k$ w.r.t. $u_j$, with $j\neq k$ and using the Schlesinger equations for $A_1,\dots,A_n$ we get:
\begin{eqnarray}\nn
&&
\frac{\partial \tilde A_k}{\partial u_j}=\left[\tilde A_k, (I(u_k)+G)^{-1} \frac{\partial G}{\partial u_j}
+\frac{ (I(u_k)+G)^{-1} A_j(I(u_k)+G)}{u_k-u_j}
\right]=\nn\\
&&
\qquad\qquad=\frac{\left[\tilde A_k, \tilde A_j\right] }{u_k-u_j}+\nn\\
&&
+
 \left[\tilde A_k,(I(u_k)+G)^{-1}\left(\frac{\partial G}{\partial u_j}
+\frac{A_j(I(u_k)-I(u_j))-B_{kj} A_j(I(u_j)+G)
}{u_k-u_j}\right)\right],
\nn\end{eqnarray}
where 
$$
B_{kj}=\left(\begin{array}{cccc}
0&\dots&0&\frac{u_k-u_j}{g_{m1}}\\
0&\dots&0&0\\
\dots&\dots&\dots&\dots\\
0&\dots&0&0
\end{array}\right).
$$
Given the formulae (\ref{newgauge}), it is straightforward to prove that the equation
$$
\frac{\partial G}{\partial u_j}
+\frac{A_j(I(u_k)-I(u_j))-B_{kj} A_j(I(u_j)+G)}{u_k-u_j}=0,
$$
is equivalent to the equations (\ref{ddpsi1}), (\ref{ddpsi2}).
This proves that also $\tilde A_1,\dots,\tilde A_n$ satisfy the Schlesinger equations.

Analogous formulae 
can be derived for the transformation $\lambda^{(k)}_1\to\lambda^{(k)}_1+1$,
$\lambda^{(k)}_m\to\lambda^{(k)}_m-1$, for $k=1,\dots,n$. In fact, 
suppose $u_1=0$ 
and $k\neq1$. We can simply apply the conformal transformation 
$\tilde{z}={1\over u_k}- {1\over{z}}$. The new residue matrices are 
$\tilde{A}_{l}={A}_{l}$ for $l\neq 1,\infty$, 
$\tilde{A}_{1}=-\sum_l{A}_{l}$, 
$\tilde{A}_{\infty}={A}_{k}$. We then need to diagonalize 
$\tilde{A}_{\infty}$ and apply the above gauge transformation to
the new system.\epf

We show that it is possible to make a finite sequence of gauge transformations
described in Lemma \ref{lmb0}  in such a way that the final 
Fuchsian system has exponents $\lambda^{(\infty)}_1,\dots,\lambda^{(\infty)}_m$
at infinity and $\lambda^{(\infty)}_k,\dots,\lambda^{(\infty)}_k$ at $u_k$, 
$k=1,\dots n$. 

By means of a permutation, we 
choose the following ordering of the parameters $\lambda_k^{(\infty)}$
$$
\Re\left(\lambda_m^{(\infty)}\right)\geq\Re\left(\lambda_1^{(\infty)}\right)\geq
\Re\left(\lambda_2^{(\infty)}\right)
\dots\geq
\Re\left(\lambda_{m-1}^{(\infty)}\right).
$$
We start with $\hat\lambda^{(\infty)}_2\to\hat\lambda^{(\infty)}_2+1$ and
$\hat\lambda^{(\infty)}_m\to\hat\lambda^{(\infty)}_m-1$. 
We want to apply such
gauge $s=\lambda^{(\infty)}_2-\hat\lambda^{(\infty)}_2$ times. To do this 
we need to check that for all $p=0,1,\dots,s-1$ and for all $l=1,\dots,m$ 
the following conditions are satisfied:
\begin{eqnarray}
&&\hat\lambda^{(\infty)}_m-\hat\lambda^{(\infty)}_2\neq 2p+1,2p +2\label{cond1}\\
&&\hat\lambda^{(\infty)}_l-\hat\lambda^{(\infty)}_2\neq p+1,
\qquad \forall\,l\neq 2,m,
\label{cond2} \\
&&\hat\lambda^{(\infty)}_m-\hat\lambda^{(\infty)}_l\neq p+1,
\qquad \forall\,l\neq 2,m.\label{cond3}
\end{eqnarray}
To prove (\ref{cond1}) we observe that since 
$\Re\left(\lambda^{(\infty)}_m\right)>\Re\left(\lambda^{(\infty)}_2\right)$,
$\Re\left(\hat\lambda^{(\infty)}_m-\hat\lambda^{(\infty)}_2\right)>\max(2p+2)=2s$.
To prove (\ref{cond2}) we observe that 
$\Re\left(\hat\lambda^{(\infty)}_l-\hat\lambda^{(\infty)}_2\right)$ is a negative number.
To prove (\ref{cond3}) we observe that 
$\Re\left(\hat\lambda^{(\infty)}_m-\hat\lambda^{(\infty)}_1\right)>s$.
Therefore all conditions (\ref{cond1}), (\ref{cond2}), (\ref{cond3})
are satisfied and thanks to the hypothesis that the monodromy group of the Fuchsian equation
(\ref{eqscal}) is irreducible (which implies that at each step at least one residue matrix has $m1$ entry non-identically $0$) there exists a  gauge transformation
$G_2(z)$ such that the new Fuchsian system has exponents
$\hat\lambda^{(\infty)}_2+s$, $\hat\lambda^{(\infty)}_m-s$, and
$\hat\lambda^{(\infty)}_j$ for all $j=1,3,\dots,m-1$.

At the $j$-th step of this procedure the parameters are
$\lambda^{(\infty)}_1,\dots,\lambda^{(\infty)}_j$, 
$\hat\lambda^{(\infty)}_{j+1},$ \dots, $\hat\lambda^{(\infty)}_{m-1}$,
$\tilde\lambda^{(\infty)}_m=\lambda^{(\infty)}_m-{m-j-1\over2} 
[(m+j-2)(n-1)-2]$. We want to apply a  gauge transformation 
$G_{j+1}(z)$ that maps 
$\hat\lambda^{(\infty)}_{j+1}\to\hat\lambda^{(\infty)}_{j+1}+1$,
$\tilde\lambda^{(\infty)}_m\to\tilde\lambda^{(\infty)}_m-1$, a number
$j (n-1)-1=\hat\lambda^{(\infty)}_{j+1}-\lambda^{(\infty)}_{j+1}$ of times.
As above we need to verify that for all $p=0,1,\dots,j ( n-1)-2$ 
for all $l=1,\dots,j,j+2,\dots,m-1$ the following
conditions are satisfied:
\begin{eqnarray}
&&\tilde\lambda^{(\infty)}_m-\hat\lambda^{(\infty)}_{j+1}
\neq 2p+2,2p+1,\label{cond4}\\
&&\lambda^{(\infty)}_l-\hat\lambda^{(\infty)}_{j+1}
\neq p+1,\label{cond5} \\
&&\tilde\lambda^{(\infty)}_m-\lambda^{(\infty)}_l\neq p+1,
\label{cond6}\\
&&\tilde\lambda^{(\infty)}_m-\hat\lambda^{(\infty)}_s\neq p+1,
\label{cond7}
\end{eqnarray}
The proof that these conditions are fulfilled at each step is 
straightforward. Again we can use the  hypothesis that the monodromy group of the Fuchsian equation
(\ref{eqscal}) is irreducible to prove that at each step at least one residue matrix has $mj$ entry non-identically $0$.

Therefore we have obtained a gauge transformation 
$G_m(z) G_{m-1}(z) \cdots G_2(z)$ such that the new Fuchsian system
has exponents $\lambda^{(\infty)}_1,\dots,\lambda^{(\infty)}_m$
at infinity and $\lambda^{(\infty)}_k,\dots,\lambda^{(\infty)}_k$ at $u_k$, 
$k=1,\dots n$. The new fundamental matrix at infinity is
$$
\tilde\Phi_\infty:=\prod_{j=1}^{m-1}G_{j+1}(z)
V^{(\infty)}({z}) z^{-M g +C}
{z}^{\Theta^{(\infty)}} {z}^{R}.
$$
In order to normalize it at infinity  we need to perform one last 
gauge transform:
$$
\Phi = \left(\prod_{j=1}^{m-1}G_{j+1}(\infty)
V^{(\infty)}({\infty})\right)^{-1} \tilde\Phi.
$$
The final new Fuchsian system
$$
\ddz\Phi=\sum_{k=1}^{n}{{A}_k\over{z}-u_k}\Phi
$$
has exponents 
$\lambda^{(\infty)}_1,\dots,\lambda^{(\infty)}_{m}$ at $\infty$ and
$\lambda^{(k)}_1,\dots,\lambda^{(k)}_{m}$ at $u_k$, and no apparent 
singularities.
The matrix entries of the matrices ${A}_k$, 
$k=1,\dots,n$, depend rationally on the coefficients of the polynomials 
$f_{l}$ and on $q_1,\dots,q_g$. This concludes the proof of existence. 
Due to the ambiguity $\prod_{j=1}^{m-1}G_{j+1}(z)\to D
\prod_{j=1}^{m-1}G_{j+1}(z)$ where $D$ is any constant diagonal matrix with 
non-zero entries, we have that from the differential 
equation \ref{eqscal} we have constructed not one Fuchsian system, but a family
of them, all related by diagonal conjugation.

We now prove the last statement of the theorem. Let us start from another
Fuchsian system
\begin{equation}\label{tils}
\ddz\check\Phi=\sum_{k=1}^{n}{\check{A}_k\over{z}-u_k}\check\Phi
\end{equation}
with exponents 
$\lambda^{(\infty)}_1,\dots,\lambda^{(\infty)}_{m}$ at $\infty$ and
$\lambda^{(k)}_1,\dots,\lambda^{(k)}_{m}$ at $u_k$, and no apparent 
singularities. Let us normalize its fundamental matrix
$\check\Psi$ as usual 
$$
\check\Phi^{(\infty)}:=\left(\ID+
{\cal O}({1\over z})\right)z^{-A^{(\infty)}}z^{-{R}^{(\infty)}}. 
$$
Let us apply the reduction procedure described 
in theorem \ref{skal-red}. This means that we construct a gauge transformation
$G(z)={\cal Q}^{-1}{\cal P}$,
$$
{\cal Q}^{-1}{\cal P}\check\Phi=Y
$$
where $Y$ is the Wronskian matrix of the differential equations (\ref{eqscal}).
Now from such equation we constructed a Fuchsian system 
$$
\ddz\Phi=\sum_{k=1}^{n}{{A}_k\over{z}-u_k}\Phi
$$
with the same exponents as (\ref{tils}). By the above construction, $\Phi$ 
is also related to $Y$ by a gauge, $\Phi=\check G(z) Y$. Therefore
$\Phi=\check G(z) G(z)\check\Phi$. This gauge transformation preserves the
normalization at infinity by construction. All monodromy data 
are preserved and by the uniqueness lemma \ref{lm2.8} we 
conclude that $\check G(z) G(z)$ must be diagonal and constant in $z$.

In particular this proves that the first row of the fundamental matrix $\Phi$
satisfies the given $m$-th order Fuchsian equation. \epf


\subsection{Darboux coordinates for Schlesinger system.}

According to  Corollary \ref{corallo}, the parameters $q_i$, $p_i$,
$i=1, \dots, g$ are coordinates on a Zariski open
subset in the space of all special Fuchsian equations with given indices and given $u_1,\dots,u_n$.

Due to Theorems \ref{skal-red} and \ref{eqscal}, the parameters $q_i$, $p_i$,
$i=1, \dots, g$ can be used as coordinates on a Zariski open
subset in the space of all Fuchsian systems with given $u_1,\dots,u_n$ and given exponents, considered modulo diagonal conjugations.  Indeed, for fixed $u_1,\dots, u_n$, the condition (\ref{nonzero1}) defines a Zariski open set in the space of all Fuchsian systems with given exponents. 

In this Section we will prove that these coordinates
are canonically conjugated with respect to the isomonodromic symplectic structure $\omega_K$ (see (\ref{omegak})) on (\ref{red-leaf}).

\begin{rmk}\label{rmknonzero} In order to apply our coordinates to the description of solutions to the Schlesinger equations, one has to make sure that the Zariski closed subset where the map
$$
\left\{\begin{matrix} {\rm Fuchsian} ~ {\rm systems} ~ {\rm with} ~{\rm given}
~{\rm poles}\\
{\rm and}~ {\rm given} ~ {\rm eigenvalues} ~ {\rm of} ~ A_1, \dots, A_n,
A_\infty\\
{\rm modulo} ~ {\rm diagonal} ~ {\rm conjugations} \end{matrix}\right\} \to
(q_1, \dots, q_g, p_1, \dots, p_g)
$$
becomes singular, or equivalently condition (\ref{nonzero1}) is violated, is never invariant under the monodromy preserving deformation.\footnote{In the theory of {\it iso-spectral}\/ deformations an analogous problem arises. In this case one needs to check that the dynamics on the Jacobian of the spectral curve is never tangent to the Theta-divisor.} This can be proved under the following two assumptions:

\noindent i) If $A_\infty$ has a resonance of order one then the corresponding logarithmic correction $R_1$  is not zero (see (\ref{lev5})).
 
\noindent ii) For at least one $j$, the entries of the $j$-th row of the matrices $A_1$, \dots, $A_n$ satisfy the following condition 
\begin{equation}\label{nonzero}
{\rm for }~{\rm every}~ i\neq j ~ {\rm there} ~{\rm exists} ~k~ {\rm such}~{\rm
that} ~{A_k}_{ji}\neq 0.
\end{equation}

\noindent Under these assumptions, by performing a small monodromy preserving deformation, condition (\ref{nonzero1}) is satisfied.

In fact suppose by contradiction that $\sum_{l=1}^n u_l
A_{l_{ji}}(u)\equiv 0$ for some $i,j$.  A simple differentiation using the Schlesinger equations gives
\begin{equation}\label{weakcondition}
{\pal\over \pal u_k} \sum_{l=1}^n u_l
A_{l_{ji}}=-\left(1+\lambda_j^{(\infty)}
-\lambda_i^{(\infty)}\right) A_{k_{ji}}=0.
\end{equation}
Now if $1+\lambda_j^{(\infty)} -\lambda_i^{(\infty)}=0$ then $A_\infty$ has a resonance of order one and since $R_{1_{ji}}= \sum_{l=1}^n u_l
{A_l}_{ji}$ must be zero, assumption i) is contradicted. Therefore $ {A_k}_{ji}=0$, but this contradicts assumption ii).

Clearly the two assumptions are satisfied in a large Zariski open set in the space of solutions of the Schlesinger equations. 
\end{rmk}


Let us rewrite the equation (\ref{eqscal}) 
in the matrix form 
\begin{equation}
\ddz\Psi={\cal B}({z})\Psi
\label{matB}\end{equation}
where 
\begin{equation}
\Psi=\left(\begin{array}{c}
y\\ y'\\ \dots \\ y^{(m-1)} \\
\end{array}\right), \qquad
{\cal B}({z})=\left(\begin{array}{ccccc}
0&1&0&\dots &0\\ 0&0&1&0&\dots\\
0&\dots&0&1&0\\\dots &\dots &\dots &\dots &\dots\\
d_{0}&d_{1}&\dots &\dots& d_{m-1}
\end{array}\right),
\label{matB1}\end{equation}
with $d_{l}(z) := {f_{l}({z})\over 
\Delta({z}){R}({z})^{m-l}}$, $R(z)=\prod_{k=1}^{n}({z}-u_k)$, 
$\Delta(z)=\prod_{i=1}^{g}({z}-q_i)$ and ${\rm deg}\, f_l(z)=(n-1)(m-l)+g.$
Recall that the system (\ref{matB1}) is obtained from the original Fuchsian system by a gauge transformation.

\begin{lm}
If the apparent singularities $q_i$, $i=1,\dots,g$ in the equation
(\ref{eqscal}) are distinct,
for each $i=1,\dots,g$, the matrix ${\mathcal B}({z})$ has one 
and only one eigenvalue $\rho_i({z})$ with a
simple pole at $q_i$. For each $i=1,\dots,g$, we define
$$
p_i = \left\{\begin{array}{l}
\res_{z=q_i}{\rho_i({z})\over{z}-q_i},\quad
\hbox{for}\quad 
q_i\neq u_k,\,\forall\,k=1,\dots,n,\\
\res_{z=q_i}{\rho_i({z})},\quad
\hbox{for}\quad  
q_i= u_k,\\
\end{array}\right.
$$
Then 
\begin{equation}
p_i={\tilde f_{m-2}(q_i)\over\tilde f_{m-1}(q_i)}
{{\tilde f}'_{m-1}(q_i)\over\Delta'(q_i)}-
{\tilde f}_{m-1}(q_i){\Delta''(q_i)\over\Delta'(q_i)}=\res_{z=q_i}
\left( d_{m-2}(z) + {1\over 2}d_{m-1}(z)^2\right),
\label{pi}\end{equation} 
where for $l=0,\dots,m-1$
$$
{\tilde f}_{l}({z}) = \left\{\begin{array}{l}
{f_l({z})\over{R}({z})^{m-l}},\quad
\hbox{\rm for}\quad 
q_i\neq u_k,\,\forall\,k=1,\dots,n,\\
{f_l({z})\over{{R}'({z})}^{m-l}},\quad
\hbox{\rm for}\quad  
q_i= u_k.\\
\end{array}\right.
$$
\end{lm}

Observe that as an immediate consequence of the second part of the equation (\ref{pi}), one obtains that the momenta $p_i$ coincide with those defined in (\ref{mom}).

\noindent{\bf Proof.}
The characteristic equation of ${\cal B}$ is
$$
\Delta({z})\rho({z})^m = \sum_{l=0}^{m-1}{f_l({z})
\over{R}({z})^{m-l}}\rho({z})^{l}. 
$$
Let us define $\tilde\rho({z})={R}({z})\rho({z})$.
Since the polynomials $f_l({z})$ are regular at ${z}=q_i$,
there is only one eigenvalue $\tilde\rho_i({z})$ that has a 
pole at $q_i$. If $q_i\neq u_k$ for all $k=1,\dots,n$, this pole is simple.
Let us expand $\tilde\rho_i$ at 
${z}=q_i$ and compare the left and right hand sides of the 
characteristic equation. We obtain
$$
\begin{array}{ll}
\tilde\rho_i({z}) = &{{f}_{m-1}(q_i)\over\Delta'(q_i)({z}-q_i)}+ 
{f_{m-2}(q_i)\over f_{m-1}(q_i)}+
{{f}'_{m-1}(q_i)\over\Delta'(q_i)}-
{f}_{m-1}(q_i){\Delta''(q_i)\over\Delta'(q_i)}
+{\mathcal O}({z}-q_i).\\
\end{array} 
$$
Therefore $\rho_i({z}){{\tilde f}_{m-1}(q_i)\over\Delta'(q_i)({z}-q_i)}+p_i+
{\mathcal O}({z}-q_i)$, where 
${\tilde f}_{l} = {f_l\over{R}(z)^{m-l}}$. This proves 
(\ref{pi}) for $q_i\neq u_k$ for all $k=1,\dots,n$.
Analogously if $q_i= u_k$ for one value of $k=0,\dots,\infty$, then
$\rho_i({z})$ has a double pole at $q_i$ and $\rho_i({z}){{\tilde f}_{m-1}(q_i)\over\Delta'(q_i)({z}-q_i)^2}+
{p_i\over({z}-q_i)}+{\mathcal O}(1)$ and again we obtain
(\ref{pi}) as we wanted to prove. 
The second part of formula (\ref{pi}) is immediately obtained
from the formula (\ref{papper}) for $a_1(z)=d_{m-1}(z)$. \epf

\begin{df}
We call the set $(q_1,\dots,q_g,p_1,\dots,p_g)$ the {\it isomonodromic
coordinates} of the Schlesinger equations.
\end{df}

\begin{thm}\label{kanon}
On a generic reduced symplectic leaf
${\cal O}_1\times\dots\times{\cal O}_{n}/Diag$ 
the quantities $(q_1,\dots,q_g,p_1,\dots,p_g)$ 
are canonical coordinates. The Schlesinger equations in these
coordinates are written in the canonical form
\begin{eqnarray}
&&
{\pal q_i\over \pal u_k} = {\pal {\mathcal H}_k \over \pal p_i}
\nn\\
&&
{\pal p_i\over \pal u_k} = -{\pal {\mathcal H}_k \over \pal q_i}\nn
\end{eqnarray}
where the Hamiltonians in canonical coordinates are
given by the formula
\begin{equation}
{\mathcal H}_k = -\res_{z=u_k}\left(d_{m-2}({z})+{1\over2}d_{m-1}({z})^2\right),
\label{hamnostra}\end{equation}
where $d_{m-2}({z})$ and $d_{m-1}({z})$ are defined in Theorem 
\ref{skal-red}.
\end{thm}

\begin{cor}\label{expham}
The Hamiltonians (\ref{hamnostra}) are given by
\begin{eqnarray}\nn
&&
{\mathcal H}_k=\left[\sum_{s=1}^g\frac{\delta_1^{(s)}-p_s}{u_k-q_s}+
\sum_{i\neq k}\frac{\beta_{2}^{(i)}R'(u_i)}{u_k-u_i}  -
\frac{\beta_2^{(k)}R''(u_k)}{2}-
\beta_{2}^{(\infty)} u_k^{n-2} -P_{n-3}(u_k)\right] {1\over R'(u_k)}-\nn
\\ 
&&
 -\left\{\sum_{s=1}^g {1\over u_k-q_s} + \sum_{i\neq k} { 
\sum_{j=1}^m \lambda_j^{(i)} - {m\, (m-1)\over 2}\over u_k-u_i}\right\}\left[
\sum_{j=1}^m \lambda_j^{(i)} - {m\, (m-1)\over 2}\right],
\nn\end{eqnarray}
where the coefficients of the polynomial $P_{n-3}(z)$ are rational functions of $p_1,\dots,p_g$, $q_1,\dots,q_g$ uniquely determined by (\ref{main-ur}).
\end{cor}

\begin{ex}\label{ex33}
In the $2\times 2$ case the polynomial $\Delta(z)$ coincides with the
$(1,2)$-matrix entry of $A(z)=\sum_k A_k /(z-u_k)$,
$$
\Delta(z) = R(z) A_{12}(z).
$$
So the isomonodromic coordinates $q_i$ coincide with the spectral coordinates
(see below).
Our $p_1,\dots,p_g$ are slightly different
from the usual momenta $\hat p_1,\dots,\hat p_g$ defined for Garnier 
systems (see \cite{IKSY}). In fact in our case we imposed
the trace of all matrices ${A}_k$ to be zero, while in \cite{IKSY}, the
determinant is zero. 
The relation between our coordinates and \cite{IKSY} is given by
$$
p_i= -\hat p_i +\sum_{k=1}^{n} {\lambda_1^{(k)}\over q_i-u_k}+ 
{q_i(2 \sum_{k=1}^{n} u_k-n q_i)+\sum_{1\leq k<l\leq n} u_k u_l
\over\prod_{k=1}^{n} (q_i-u_k)}
$$
Keeping track of this time-dependent canonical transformation, it is not 
difficult to verify that our Hamiltonian functions (\ref{hamnostra}) 
coincide with the one given in \cite{IKSY}.
\end{ex}

\noindent{\bf Proof of the Theorem.}
Since the system (\ref{matB1}) is gauge equivalent to the original Fuchsian system (\ref{N1}), it suffices to perform all computations with (\ref{matB1}).
First of all formula (\ref{hamnostra}) is obtained by straightforward computation applying formula (\ref{hamk}) to the matrix (\ref{matB1}). We want to show that  $q_i$, $p_i$  are canonical coordinates on
the reduced symplectic leaf (\ref{red-leaf}) and that in those coordinates the Hamiltonian is indeed (\ref{hamnostra}). To this aim observe that we can always put equations of the form (\ref{eqscal}) 
in the matrix form (\ref{matB}), (\ref{matB1}).  
The apparent singularities are poles of the eigenvectors of the matrix ${\cal B}({z})$. 

In this way the proof reduces to proving the following

\begin{lm}
If $q_i\neq u_1,\dots,u_{n},\infty$ for all $i=1,\dots,g$,
the symplectic structure (\ref{omegak}) on the space of monodromy 
data of the linear system of ODEs of the form (\ref{matB}),
(\ref{matB1}) is
$$
\omega_K = \sum_{i=1}^g{\rm d}p_i\wedge{\rm d}q_i.
$$
\end{lm}

\noindent{\bf Proof.}
Due to gauge invariance of the form $\omega_K$ we have to compute
$$
\omega_K=-{1\over2}\sum_{k=1}^{n} \res_{z=u_k}{\rm Tr}
(\delta{\cal B}\wedge\delta\Psi\,\Psi^{-1})
-{1\over2}\sum_{i=1}^g \res_{z=q_i}{\rm Tr}
(\delta{\cal B}\wedge\delta\Psi\,\Psi^{-1}).
$$
Observe that $\delta{\cal B}_{i l}$ is zero for all $i\neq m$ and 
$\delta{\cal B}_{m l}=\delta d_{l-1}$ so that $\omega_K$ depends
only on the $m$-th column of the matrix $\Psi^{-1}$, i.e. 
$$
\begin{array}{ll}
\omega_K =& -{1\over2}\sum_{k=1}^{n} \res_{z=u_k}\sum_{l,j=0}^m
\delta d_{l-1}\wedge\delta\Psi_{lj} \left(\Psi^{-1}\right)_{jm} -\\
&-{1\over2} \sum_{i=1}^g \res_{z=q_i}\sum_{l,j=0}^m
\delta d_{l-1}\wedge\delta\Psi_{lj} \left(\Psi^{-1}\right)_{jm}.\\
\end{array}
$$ 
 
Let us deal with the apparent singularities first. We have to compute 
the expansion of ${\cal B}({z})$ and $\Psi({z})$ at $q_i$.

Choose
$m$ linear independent solutions ${y}_1,\dots,{y}_m$ having expansions 
at ${z}=q_i$ of the form described in Lemma \ref{basisapp}
where the constants $\alpha_l$ are determined by the differential
equation (\ref{matB}) and are
\begin{eqnarray}\label{eqalpha}
&&
\alpha_{l+1}^{(i)} = -\res_{z=q_i}d_l({z}), 
\qquad l=0,\dots,m-2\nn\\
&&
\alpha_m^{(i)} = \res_{z=q_i}\left(2{d_{m-1}({z})\over{z}-q_i}+
d_{m-2}({z})\right).
\end{eqnarray}
Comparing these with (\ref{pi}) we have 
$p_i={\alpha_m^{(i)}-\alpha_{m-1}^{(i)}\over2}=\res_{z=q_i}\left({d_{m-1}({z})\over z-q_i}+d_{m-2}({z})\right)$. 

Observe that the fundamental matrix $\Psi$ of the scalar equation 
(\ref{eqscal}) has matrix elements given
by
$$
\Psi_{lj}={{\rm d}^{l-1}{y}_j\over{\rm d}{z}^{l-1}},
\qquad l,j=1,\dots,m.
$$
We can show that in the computation of the residue in (\ref{omegak})
at $q_i$, one can neglect ${\mathcal O}({z}-q_i)^{m}$ in 
${y}_1,\dots,{y}_m$ and in the coefficient  $d_{m-1}$ and  ${\mathcal O}(1)$ in 
$d_0,\dots,d_{m-2}$. This follows by straightforward computations based on
a list of observations:
\begin{enumerate}
\item $\left(\Psi^{-1}\right)_{jm}={\cal O}({z}-q_i)^{m-j}$ for
$j=1,\dots,m-1$ and $\left(\Psi^{-1}\right)_{mm}{\cal O}({z}-q_i)^{-1}$.
\item From (\ref{mom}) and (\ref{eqalpha}) we have that
$$
{\alpha_m^{(i)}+\alpha_{m-1}^{(i)}\over 2}=\sum_{j\neq i}
{1\over q_i-q_j} +\sum_{k=1}^{n} 
{\sum_{j=1}^{m}\lambda_j^{(k)}-{m(m-1)\over2}\over q_i-u_k}.
$$ 
\item For $l=1,\dots,m-1$, 
$$
\delta d_{l-1}= {-\alpha_l^{(i)}\delta q_i\over({z}-q_i)^2}-
{\delta\alpha_l^{(i)}\over({z}-q_i)}+{\cal O}(1),
$$
and
$$
\delta d_{m-1}= {\delta q_i\over({z}-q_i)^2}+{1\over2}
\left(\delta\alpha_m^{(i)}+\delta\alpha_{m-1}^{(i)}\right)+
{\cal O}(1)\delta q_i+{\cal O}({z}-q_i).
$$
\item For $1\leq j\leq m-1$ and $l\leq j-1$:
\begin{eqnarray}\nn
&&
\delta\Psi_{lj}= -{\delta q_i\over(j-l-1)!} ({z}-q_i)^{j-l-1}+
{\delta\alpha_j^{(i)}\over(m-l)!}({z}-q_i)^{m-l}-\nn\\
&&
-{\alpha_j^{(i)}\delta q_i\over(m-l-1)!} ({z}-q_i)^{m-l-1}+ 
\delta q_i
{\cal O}({z}-q_i)^{m-l+1}+{\cal O}({z}-q_i)^{m-l+2},
\nn\end{eqnarray}
for $1\leq j\leq m-1$ and $l\geq j$:
\begin{eqnarray}
&&
\delta\Psi_{lj}= {\delta\alpha_j^{(i)}\over(m-l)!}({z}-q_i)^{m-l}
-{\alpha_j^{(i)}\delta q_i\over(m-l-1)!} ({z}-q_i)^{m-l-1}+\nn \\
&&
+\delta q_i {\cal O}({z}-q_i)^{m-l+1}+
{\cal O}({z}-q_i)^{m-l+2},
\nn\end{eqnarray}
and for $j=m$
\begin{eqnarray}\nn
&&
\delta\Psi_{lm}= -{\delta q_i\over(m-l)!} ({z}-q_i)^{m-l}+
{\delta\alpha_m^{(i)}\over(m-l+2)!}({z}-q_i)^{m-l+2}-\nn\\
&&
-{\alpha_m^{(i)}\delta q_i\over(m-l-1)!} ({z}-q_i)^{m-l-1}+ 
\delta q_i
{\cal O}({z}-q_i)^{m-l+2}+{\cal O}({z}-q_i)^{m-l+3}.
\nn\end{eqnarray}
\end{enumerate}
From 1) and 3) we immediately see that only the terms with $j=m-1,m$
can contribute to the residue. From 1), 2), 3) and 4) we have 
\begin{eqnarray}\nn
&&
\sum_{i=1}^g \res_{z=q_i}
\sum_{l}\delta d_{l-1}\wedge\delta\Psi_{l,m-1}
\left(\Psi^{-1}\right)_{m-1,m}=\nn\\
&&
=\sum_{i=1}^g \res_{z=q_i}\delta d_{m-1}\wedge\delta\Psi_{m,m-1}
\left(\Psi^{-1}\right)_{m-1,m}{1\over2} \delta\alpha_{m-1}^{(i)}\wedge\delta q_i,
\nn\end{eqnarray}
and 
\begin{eqnarray}\nn
&&
\sum_{i=1}^g \res_{z=q_i}
\sum_{l}\delta d_{l-1}\wedge\delta\Psi_{lm}\left(\Psi^{-1}\right)_{mm} =
\nn\\
&&
=\sum_{i=1}^g \res_{z=q_i}
\left(\delta d_{m-1}\wedge\delta\Psi_{mm}
\left(\Psi^{-1}\right)_{mm}+\delta d_{m-2}\wedge\delta\Psi_{m-1,m}
\left(\Psi^{-1}\right)_{mm}\right) =\nn\\
&&
=-\delta\alpha_{m}^{(i)}\wedge\delta q_i+\frac{1}{2}\delta\alpha_{m-1}^{(i)}\wedge\delta q_i,
\nn\end{eqnarray}
thus we obtain
$$ 
-{1\over2} \sum_{i=1}^g \res_{z=q_i}\sum_{l,j}
\delta d_{l-1}\wedge\delta\Psi_{lj} \left(\Psi^{-1}\right)_{jm}=-\frac{1}{2}  \sum_{i=1}^g \delta\left(\alpha_{m-1}^{(i)}-\alpha_m^{(i)}\right)\wedge\delta q_i= 
\sum_{i=1}^g \delta p_i\wedge\delta q_i.
$$

Let us now show that
$$
\res_{z=u_k}\sum_{l,j}
\delta d_{l-1}\wedge\delta\Psi_{lj} \left(\Psi^{-1}\right)_{jm}=0,\qquad\forall\,k=1,\dots,n.
$$
Let us expand $d_{l-1}$ at $u_k$. We have
$d_{l-1} = {c_{l-1}\over({z}-u_k)^{m-l+1}} +
{\mathcal O}({z}-u_k)^{l-m}$
where $c_{l-1}$ are uniquely determined by the indicial equation, thus
by the exponents. As a consequence $\delta d_{l-1}= 
{\mathcal O}({z}-u_k)^{l-m}$. Analogously to estimate 
$\Psi_{lj} = {{\rm d}^{l-1}{y}_j\over{\rm d}{z}^{l-1}}$, we can again 
normalize the solutions ${y}_j$ at $u_k$ is such a way that
${y}_j=({z}-u_k)^{\lambda_j^{(k)}} + 
{\mathcal O}({z}-u_k)^{\lambda_j^{(k)}+1}$ so that  
$\delta{y}_j= {\mathcal O}({z}-u_k)^{\lambda_j^{(k)}-1}$. Thus
$\delta d_{l-1}\wedge\delta\Psi_{lj}={\mathcal O}({z}-u_k)^{\lambda_j^{(k)}-m}$.
Now $\left(\Psi^{-1}\right)_{jm}={\mathcal O}({z}-u_k)^{m-\lambda_j^{(k)}}$
so that near $u_k$ we have $\sum_{l,j}\delta d_{l-1}\wedge\delta\Psi_{lj}
\left(\Psi^{-1}\right)_{jm}={\mathcal O}(1)$ and the pole $u_k$ does
not contribute to the residue. \epf

The above lemma proves that  $q_i$, $p_i$  are canonical coordinates on
the reduced symplectic leaf (\ref{red-leaf}). We now want to prove that in those coordinates the Hamiltonians are indeed given by formula (\ref{hamnostra}).  To this aim we need to extend the phase space.


Let us consider the space of all matrices ${\cal B}$ of the form (\ref{matB1}) with coefficients $d_0,\dots,d_{m-1}$, 
\begin{eqnarray}
&&
d_{m-1}(z) = \sum_{s=1}^g {1\over z-q_s} + \sum_{i=1}^n {1\over z-u_i} \left[
\sum_{j=1}^m \lambda_j^{(i)} - {m\, (m-1)\over 2}\right]
\nn\\
&&
d_k(z) = \left[-\sum_{s=1}^g {c_{k+1}^{(s)}R(q_s)^{m-k-1}\over z-q_s}
+ (-1)^{m-k-1} \sum_{i=1}^n {\beta_{m-k}^{(i)}\over z-u_i} [R'(u_i)]^{m-k-1}
\right.+\nn\\
&&
\left.+\beta_{m-k}^{(\infty)} z^{(m-k)(n-1) - n } 
+P_{(m-k)(n-1) - n  -1}(z)\right] {1\over R(z)^{m-k-1}},\nn
\\
&&
k=0, \dots, m-2.
\nn
\end{eqnarray}
We remind that this means that the equation (\ref{eqscal}) has $n+1$ Fuchsian poles at $u_1,\dots,u_n,\infty$ with indices $\lambda_1^{(k)},\dots,\lambda_m^{(k)}$ for $k=1,\dots,n,\infty$ given by equations (\ref{betas}) and (\ref{betas1}) and it has simple poles at the points $q_1,\dots,q_g$. 

Near each simple pole $q_s$, the matrix ${\cal B}$ can be expanded as
$$
{\cal B}={{\cal B}_0^{(s)}\over z-q_s} + {\cal B}_1^{(s)}+{\mathcal O}(z-q_s)
$$
where 
$$
{\cal B}_0^{(s)}=\left(\begin{array}{cccc}
0&0&\dots &0\\ 
\dots &\dots  &\dots &\dots\\
0&0&\dots&0\\
-c_1^{(s)}&\dots&-c_{m-1}^{(s)}&1\\
\end{array}\right),\\ \quad 
{\cal B}_1^{(s)}=\left(\begin{array}{ccccc}
0&1&0&\dots &0\\ 0&0&1&0&\dots\\
0&\dots&0&1&0\\\dots &\dots &\dots &\dots &\dots\\
\delta_m^{(s)}&\delta_{m-1}^{(s)}&\dots &\dots& \delta_1^{(s)}
\end{array}\right),
$$
where $\delta_1^{(s)},\dots,\delta_m^{(s)}$ are given in equations (\ref{mom}).
If we do not impose equations (\ref{ur-ja}), that is if we do not assume the singularities $q_1,\dots,q_g$ to be apparent,  there exists a fundamental solution of the form
\begin{equation}\label{psiR}
\Psi=G^{(s)}(z)(z-q_s)^\Lambda(z-q_s)^{R^{(s)}}
\end{equation}
where $G^{(s)}(z) = G_0^{(s)}+G_1^{(s)}(z-q_s)+{\mathcal O}(z-q_s)^2$, 
$\Lambda={\rm diagonal}(0,\dots,0,1)$, $R^{(s)}$ is an off-diagonal matrix with all entries equal to zero, apart from the last row. The matrices $R^{(s)}$, $G^{(s)}_0$  and $G^{(s)}_1$ are determined by the following equations:
$${\mathcal B}^{(s)}_0 G^{(s)}_0=G^{(s)}_0\Lambda,\qquad
{\mathcal B}^{(s)}_0 G^{(s)}_1+{\mathcal B}^{(s)}_1 G^{(s)}_0=G^{(s)}_1+ G^{(s)}_0 R^{(s)}+
G^{(s)}_1\Lambda.
$$
This fact is a simple consequence of the gauge formula applied to the gauge $\Psi=G^{(s)}(z)\tilde\Psi$
that maps ${\cal B}$ to ${\Lambda\over z-q_s}+R^{(s)}$.

Observe that the $m-1$ equations (\ref{ur-ja}) imposing that the simple pole $q_s$ is apparent, coincide with the $m-1$ equations $R^{(s)}\equiv 0$. 

We consider now the symplectic structure on the extended space of isomonodromic deformations of systems of the form  (\ref{matB}),  (\ref{matB1}) where the entries of the matrices $R^{(s)}$, $s=1,\dots,g$ are not necessarily null. 


The following lemma concludes the proof of Theorem \ref{kanon}, by proving that the Hamiltonians in the isomonodromic coordinates are indeed given by (\ref{hamnostra}).  

\begin{lm}
On the extended space $S\times X_n$, where $S$ are the symplectic leaves and
$X_n=\complessi^n\setminus\{diagonals\}$ is the configuration space of $n$ points, the symplectic structure (\ref{omegak}) becomes
$$
\omega_K = \sum_{i=1}^g{\rm d}p_i\wedge{\rm d}q_i-
\sum_{k=1}^n {\rm d}{\mathcal H}_k\wedge{\rm d}u_k.
$$
\end{lm}

\noindent{\bf Proof.} 
There are two main differences with the previous proof. The first one is that now we have to take into account the variations $\delta u_k$, $k=1,\dots,n$, the second one is that now the entries of the matrices $R^{(s)}$ are not necessarily zero. Let first look at the term $\delta\Psi\wedge\Psi^{-1}$ near the point $q_s$. Using formula (\ref{psiR}) we obtain
\begin{eqnarray}\nn
&&
\delta\Psi \Psi^{-1}=\delta G^{(s)}(z)\left(G^{(s)}(z)\right)^{-1} - 
G^{(s)}(z){\Lambda\over z-q_s}\left(G^{(s)}(z)\right)^{-1}\delta q_s-\nn\\
&&
\qquad\qquad-G^{(s)}(z){R^{(s)}\over z-q_s} \left(G(z)^{(s)}\right)^{-1}\delta q_s\nn
\end{eqnarray}
because all resonances are of order one and $(z-q_s)^\Lambda R^{(s)} (z-q_s)^{-\Lambda }= R^{(s)}$.

Therefore, when computing the residue at $q_s$ of ${\rm Tr}(\delta{\cal B}\wedge\delta\Psi\,\Psi^{-1})$
we just need to add the contribution of the term $-{\rm Tr}\left(\delta{\cal B}\wedge G(z)^{(s)}{R^{(s)}\over z-q_s} \left(G(z)^{(s)}\right)^{-1}\right)\delta q_s$ to ${\rm d}p_s\wedge{\rm d} q_s$. We are now going to prove that this extra contribution is zero.

In fact
$$
G(z){R^{(s)}\over z-q_s} \left(G(z)^{(s)}\right)^{-1}={G_0^{(s)}R^{(s)}(G^{(s)}_0)^{-1}\over z-q_s}-\left[G_0^{(s)}R^{(s)}(G^{(s)}_0)^{-1}, G_1^{(s)} (G_0^{(s)})^{-1}\right]+{\mathcal O}(z-q_s).
$$
Only the last column contributes to the trace. It is not difficult to see that the last column of $G_0^{(s)}R^{(s)}G^{(s)}_0$ is zero and that the only non--zero element of the last column of 
$\left[G_0^{(s)}R^{(s)}(G^{(s)}_0)^{-1}, G_1^{(s)} (G_0^{(s)})^{-1}\right]$ is the last one. Since 
$\delta d_m={\delta q_s\over z-q_s}+{\mathcal O}(1)$, the residue is zero.

Let us now show that
$$
\res_{z=u_k}\sum_{l,j}
\delta d_{l-1}\wedge\delta\Psi_{lj} \left(\Psi^{-1}\right)_{jm}=2\delta {\mathcal H}_k\wedge \delta u_k,\qquad\forall\,k=1,\dots,n.
$$
The contribution of the matrices $R^{(s)}$ does not play any role here because we are expanding at $u_k$. On the other side, this time we need to take into account the variations $\delta u_1,\dots,\delta u_n$.

Let us expand $d_{l-1}$ at $u_k$. We have
\begin{eqnarray}\label{expduk}
d_{m-1}=\frac{\beta_1^{(k)}-\frac{m(m-1)}{2}}{z-u_k} + D_{m-1}^{(k)} + {\mathcal O}(z-u_k)\\
d_{l-1}=\frac{(-1)^{m-l}\beta_{m-l+1}^{(k)}}{(z-u_k)^{m-l+1}} + \frac{D_{l-1}^{(k)}}{(z-u_k)^{m-l}} +
{\mathcal O}(z-u_k)^{l+1-m}, \nn
\end{eqnarray}
where (cfr. (\ref{papper}) with $d_{l-1}(z)=a_{m-l+1}(z)$)
\begin{eqnarray}\nn
&&
D_{m-1}^{(k)} = \sum_{s=1}^g {1\over u_k-q_s} + \sum_{i\neq k} {1\over u_k-u_i} \left[
\sum_{j=1}^m \lambda_j^{(i)} - {m\, (m-1)\over 2}\right]
\nn\\
&&
D_{l-1}^{(k)} =  \left[   (-1)^{m-l}
\sum_{i\neq k}{\beta_{m-l+1}^{(i)}\over u_k-u_i} [R'(u_i)]^{m-l}\right.+
\beta_{m-l+1}^{(\infty)} u_k^{(n-1)(m-l+1)-n} +
\nn\\
&&
\left.
\qquad+P_{(m-l+1)(n-1)-n-1}(u_k)\right] {1\over R'(u_k)^{m-l}},\qquad
l=0, \dots, m-2.
\nn
\end{eqnarray}
We need to introduce some notation:
$$
[\lambda]_0:=1,\qquad [\lambda]_1:= \lambda,\qquad
[\lambda]_n:= \lambda(\lambda-1)\dots(\lambda-n+1),\, \forall \,n=2,3,\dots.
$$
The indicial equations (\ref{betas}) read
\begin{equation}\label{betas3}
[\lambda^{(k)}]_m=\left(\beta_1^{(k)}-\frac{m(m-1)}{2}\right)[\lambda^{(k)}]_{m-1}+
\sum_{l=1}^{m-1}(-1)^{m-l}\beta_{m-l+1}^{(k)} [\lambda^{(k)}]_{l-1}
\end{equation}
To start with, we perform our computation in the case when the exponents
$\lambda_1^{(k)} \dots, \lambda_m^{(k)}$ of the pole $u_k$ are  non-resonant. Thank to (\ref{fund-ur}), there exists a basis of solutions $y_1,\dots,y_m$ of the form
$$
y_i(z)=(z-u_k)^{\lambda_i^{(k)}}(1+\eta^{(k)}_i (z-u_k)+{\mathcal O}(z-u_k)^2),\qquad i=1,\dots,m.
$$
Therefore we have
$$
\Psi_{li}(z)= \frac{{\rm d}^{l-1}y_i}{{\rm d}z^{l-1}}
= [\lambda_i^{(k)}]_{l-1} (z-u_k)^{\lambda^{(k)}_i-l}+ [\lambda_i^{(k)}+1]_{l-1} \eta^{(k)}_i
(z-u_k)^{\lambda^{(k)}_i-l+1}+{\mathcal O}(z-u_k)^{\lambda^{(k)}_i-l+2},
$$
where the constants $\eta^{(k)}_1,\dots,\eta^{(k)}_m$ are determined by the following equations
\begin{eqnarray}\label{eqeta}\nn
&&
\left\{[\lambda^{(k)}_i+1]_m-
\left(\beta_1^{(k)}-\frac{m(m-1)}{2}\right)[\lambda^{(k)}_i+1]_{m-1}-\right.\qquad\nn\\
&&
\left.\qquad-\sum_{l=1}^{m-1}(-1)^{m-l}\beta_{m-l+1}^{(k)} [\lambda^{(k)}_i+1]_{l-1}
\right\}\eta^{(k)}_i=\sum_{l=1}^{m-1}D_{l-1}^{(k)} [\lambda^{(k)}_i]_{l-1}.
\end{eqnarray}
As above, only the last column of the inverse matrix $\Psi^{-1}$ enters in the computation of the symplectic structure $\omega$
$$
\left(\Psi^{-1}\right)_{im}=\frac{(z-u_k)^{m-1-\lambda^{(k)}_i}}{\prod_{j\neq i}
(\lambda^{(k)}_i-\lambda^{(k)}_j)}+{\mathcal O}(z-u_k)^{m-\lambda^{(k)}_i}.
$$
Proceeding in a similar way as in the first part of this proof, we arrive at the formula
\begin{eqnarray}\label{symuk}
&&
\res_{z=u_k}\sum_{l,i}
\delta d_{l-1}\wedge\delta\Psi_{li} \left(\Psi^{-1}\right)_{im}
= \sum_{i=1}^m \frac{-1}{\prod_{j\neq i} (\lambda^{(k)}_i-\lambda^{(k)}_j)}
\left\{\sum_{l=1}^{m-1}[\lambda_i^{(k)}]_{l}\, \delta D_{l-1}^{(k)}\wedge\delta u_k+\right.\nn\\  
&&
\qquad +\left(\beta_1^{(k)}-\frac{m(m-1)}{2}\right)[\lambda^{(k)}_i+1]_{m-1} \,
\delta \eta^{(k)}_i\wedge\delta u_k+\\
&&
\qquad \left.+\left(\sum_{l=1}^{m-1}(-1)^{m-l}(m-l+1)\beta_{m-l+1}^{(k)} [\lambda^{(k)}+1]_{l-1}\right)
\delta \eta^{(k)}_i\wedge\delta u_k\right\}.\nn
\end{eqnarray}
Now using the indicial equation (\ref{betas3}), we get
\begin{eqnarray}\label{neceta}
\nn
&&
\left(\beta_1^{(k)}-\frac{m(m-1)}{2}\right)[\lambda^{(k)}_i+1]_{m-1}+
\sum_{l=1}^{m-1}(-1)^{m-l}(m-l+1)\beta_{m-l+1}^{(k)} [\lambda^{(k)}+1]_{l-1} =\nn\\
&&
(\lambda_i^{(k)}-m+1)\left\{ [\lambda^{(k)}+1]_{m}-\left(\beta_1^{(k)}-\frac{m(m-1)}{2}\right)
 [\lambda^{(k)}+1]_{m-1}-\right.\\
 &&\left.-\sum_{l=1}^{m-1}(-1)^{m-l}
 \beta_{m-l+1}^{(k)} [\lambda^{(k)}+1]_{l-1}
\right\}.\nn
\end{eqnarray}
Observe that the right-hand-side of equation (\ref{neceta}) is $(\lambda^{(k)}_i-m+1)$ times the coefficient of $\eta_i^{(k)}$ in (\ref{eqeta}). Using this in equation (\ref{symuk}), we get
\begin{eqnarray}\label{symuk1}
&&
\res_{z=u_k}\sum_{l,i}
\delta d_{l-1}\wedge\delta\Psi_{li} \left(\Psi^{-1}\right)_{im}
=\sum_{i=1}^m \frac{-1}{\prod_{j\neq i} (\lambda^{(k)}_i-\lambda^{(k)}_j)}
\left\{2  [\lambda^{(k)}_i]_{m}\, \delta D^{(k)}_{m-1}\wedge \delta u_k +\right.\nn\\
&&
\left.
\qquad+ \sum_{l=1}^{m-1}[\lambda_i^{(k)}]_{l-1}(2\lambda^{(k)}_i+2-m-l) \, \delta D^{(k)}_{l-1}\wedge \delta u_k \right\}.
\end{eqnarray}
To conclude we observe that
\begin{eqnarray}\nn
&&
\sum_{i=1}^m [\lambda_i]_l \frac{2 \lambda_i-m+2-l} 
{\prod_{j\neq i} (\lambda^{(k)}_i-\lambda^{(k)}_j)}
=\sum_{i=1}^n
{\rm res}_{\lambda=\lambda_i} [\lambda]_l \,\frac{2\lambda-m-l+2}{\prod_{j=1}^n
(\lambda-\lambda_j)}=\nn\\
&&
=-{\rm res}_{\lambda=\infty} [\lambda]_l \,\frac{2\lambda-m-l+2}{\prod_{j=1}^n
(\lambda-\lambda_j)}=\left\{\begin{array}{l} 0,\quad\hbox{for } l=0,1,\dots,m-3,\\
2,\quad\hbox{for } l=m-2.\\
\end{array}\right.
\nn\end{eqnarray}
Analogously
$$
\sum_{i=1}^m \frac{[\lambda_i^{(k)}]_m} 
{\prod_{j\neq i} (\lambda^{(k)}_i-\lambda^{(k)}_j)}
= \beta_1^{(k)}-\frac{m(m-1)}{2}.
$$
Using these in (\ref{symuk1}), we get finally
\begin{eqnarray}\nn
&&
\res_{z=u_k}\sum_{l,i}
\delta d_{l-1}\wedge\delta\Psi_{li} \left(\Psi^{-1}\right)_{im}
=-2 \left(\beta_1^{(k)}-\frac{m(m-1)}{2}\right) \,\delta D_{m-1}^{(k)}\wedge\delta u_k-\nn\\
&&
\quad 2\, \delta D_{m-2}^{(k)}\wedge\delta u_k = 2 \, {\rm d} {\mathcal H}_k\wedge {\rm d}u_k,
\nn\end{eqnarray}
as we wanted to prove.

To conclude, we observe that if some of the exponents
$\lambda_1^{(k)} \dots, \lambda_m^{(k)}$ of the pole $u_k$ are resonant, then by Theorem 2.1 there exists a fundamental solution
$$
\Psi=G^{(k)}(z)(z-u_k)^{\Lambda^{(k)}}(z-u_k)^{R^{(k)}}
$$
where $G^{(k)}(z) = \sum_{j=0}^\infty G_j^{(k)}(z-u_k)$, and $R^{(k)}=\sum R^{(k)}_j$ is a finite sum of off-diagonal matrices such that
$$
(z-u_k)^{\Lambda^{(k)}}  R^{(k)}(z-u_k)^{-\Lambda^{(k)}}={R^{(k)}_0}+{R^{(k)}_1(z-u_k)}+\dots.
$$
By applying enough iterates of Lemma 4.20  we can increase the order of the resonances arbitrarily, i.e. we can always assume that 
$$
(z-u_k)^{\Lambda^{(k)}}  R^{(k)}(z-u_k)^{-\Lambda^{(k)}}={R^{(k)}_p}(z-u_k)^p+
R^{(k)}_{p+1}(z-u_k)^{p+1}+\dots.,
$$
with $p$ large enough. Then the extra term in $\delta\Psi\Psi^{-1}$ due to $R^{(k)}$ is given by
$$-G^{(k)}(z)(R^{(k)}_p(z-u_k)^{p-1}+{\mathcal O}(z-u_k)^p) \left(G(z)^{(k)}\right)^{-1}\delta u_k$$
which does not contribute to the residue.
\epf


\subsection{An example.}

As we already know, for $m=2$ the isomonodromic coordinates coincide with the spectral ones. Starting from $m=3$ they are different. 

In this subsection we give an explicit parametrization of special Fuchsian equations (\ref{scal})
in the first non-trivial case $m=3$ and $n=3$, in terms of our isomonodromic coordinates $q_1,\dots,q_g,p_1,\dots,p_g$ and compute the Hamiltonians of the Schlesinger equations in the canonical coordinates.
Note that $g=4$ for $m=3$ and $n=3$.

Starting from a Fuchsian system
$$
\frac{dY}{dz} = A(z) Y, \quad A(z) =\frac{A_1}{z-u_1} + \frac{A_2}{z-u_2}  + \frac{A_3}{z-u_3} 
$$
where $A_i$ are $3\times 3$ matrices with the eigenvalues $\lambda_1^{(i)}$, $\lambda_2^{(i)}$, $\lambda_3^{(i)}$, $i=1, \, 2, \, 3$ satisfying
$$
-(A_1+A_2+A_3) =A_\infty =\diag(\lambda_1^{(\infty)}, \lambda_2^{(\infty)}, \lambda_3^{(\infty)})
$$
we arrive at a third order Fuchisan equation with eight regular singularities with the following Riemann scheme:
$$
{\cal P} \left\{ \begin{array}{cccccccc} \infty & u_1 & u_2 & u_3 & q_1 & q_2 & q_3 & q_4 \\
\lambda_1^{(\infty)} & \lambda_1^{(1)} & \lambda_1^{(2)} & \lambda_1^{(3)} & 0 & 0 & 0 & 0\\
\lambda_2^{(\infty)}+1 & \lambda_2^{(1)} & \lambda_2^{(2)} & \lambda_2^{(3)} & 1 & 1 & 1 & 1\\
\lambda_3^{(\infty)}+1 & \lambda_3^{(1)} & \lambda_3^{(2)} & \lambda_3^{(3)} & 3 & 3 & 3 & 3
\end{array}\right\}
$$
satisfying the additional constraint of absence of logarithmic terms at the points $q_1$, \dots, $q_4$. From the previous considerations it follows that the Fuchsian equation must have the form
\eqa\label{papperitz3}
&&
y'''= \left[ \sum_{s=1}^4 \frac1{z-q_s} -3\sum_{i=1}^3 \frac1{z-u_i}\right] \, y'' +
\\
&&
+\left[\sum_{s=1}^4 \frac{ c_2^{(s)} R(q_s)}{z-q_s} +\sum_{i=1}^3 \frac{\beta_2^{(i)} R'(u_i)}{z-u_i} +\beta_2^{(\infty)} z +h\right] \, \frac{y'}{R(z)}
\nn\\
&&
+\left[ -\sum_{s=1}^4 \frac{c_1^{(s)} R^2(q_s)}{z-q_s} +\sum_{i=1}^3 \frac{\beta_3^{(i)} {R'}^2(u_i)}{z-u_i}
+\beta_3^{(\infty)} z^3 + a\, z^2 + b\, z + c\right] \,\frac{y}{ R^2(z)}.
\nn
\eeqa
Let us spell out the notations. The polynomial $R(z)$ is given by
$$
R(z)=(z-u_1)(z-u_2)(z-u_3).
$$
The coefficients $\beta_k^{(i)}$, $\beta_k^{(\infty)}$ are given by the following formulae
\eqa\label{spur1}
&&
\beta_2^{(i)} = \lambda_1^{(i)} \lambda_2^{(i)} +\lambda_1^{(i)} \lambda_3^{(i)} +\lambda_2^{(i)} \lambda_3^{(i)} -5, \quad \beta_3^{(i)}=\lambda_1^{(i)} \lambda_2^{(i)} \lambda_3^{(i)},  \quad i=1, \, 2, \, 3
\nn\\
&&
\\
&&
\beta_2^{(\infty)} =[\lambda_1^{(\infty)}]^2 +[\lambda_2^{(\infty)}]^2+[\lambda_3^{(\infty)}]^2 -2 \lambda_1^{(\infty)} -5, \quad \beta_3^{(\infty)} = -\lambda_1^{(\infty)} ( \lambda_2^{(\infty)}+1)( \lambda_3^{(\infty)}+1).
\nn
\eeqa
In this example we assume all the matrices $A_i$ to be traceless:
\beq\label{spur0}
\tr A_i =0, \quad i=1, \, 2, \, 3.
\eeq
We also put 
$$
c_2^{(s)}=-p_s-3\, \frac{R'(q_s)}{R(q_s)}+\frac12 \frac{\Delta''(q_s)}{\Delta'(q_s)}
$$
as in equation (\ref{mom}). We introduce the following quantities
\beq\label{spur2}
p_s^{[k]}=p_s + k\,  \frac{R'(q_s)}{R(q_s)}, \quad k=0, \, 1, \, 2, \, 3, \quad s=1, \dots, 4
\eeq
and
\beq\label{spur3}
\tilde p_s^{[k]} = p_s^{[k]} -\frac12 \frac{\Delta''(q_s)}{\Delta'(q_s)}
\eeq
where, as above the monic polynomial $\Delta(z)$ is defined by
\beq\label{spur4}
\Delta(z)=(z-q_1) \dots (z-q_4) \equiv z^4 - \sigma_1 z^3 + \sigma_2 z^2 - \sigma_3 z +\sigma_4.
\eeq
Here $\sigma_1$, \dots, $\sigma_4$ are just the elementary symmetric functions of $q_1$, \dots, $q_4$.  
In these notations, $c_2^{(s)}=-\tilde p^{[3]}_s$.

The coefficients $h$, $a$, $b$, $c$ and $c_1^{(s)}=\alpha_1^{(s)}$, $s=1, \dots, 4$ are to be expressed in terms of
the canonical coordinates $q_1$, \dots, $q_4$, $p_1$, \dots, $p_4$ and $u_1$, $u_2$, $u_3$ from the assumption of absence of logarithmic terms at the apparent singularities. This assumption yields a linear $8 \times 8$ system for the above unknowns.
Eliminating the unknowns $c_1^{(s)}$ one arrives at the following system:
\beq\label{sistema3}
a\, q_s^2 + b \, q_s + c + h\cdot \sum_t M_{st} (p^{[2]}, q) R(q_t) =w_s, \quad s=1, \dots, 4
\eeq
where the $4\times 4$ matrix $M(p,q) =\left( M_{st}(p,q)\right)$ is defined by
$$
M_{st}(p,q)=\left\{ \begin{array}{cc} p_s, & t=s\\ & \\ \frac1{q_t-q_s} , & t\neq s\end{array}\right. 
$$
and
\eqa\label{w-s}
&&
w_s=\sum_{t,r} M_{st}(p^{[2]}, q)R(q_t) M_{tr}(p^{[1]},q) R(q_r) \tilde p_r^{[0]}
\nn\\
&&
-\sum_t M_{st} (p^{[2]}, q) R(q_t) f_2(q_t) - f_3(q_s), \quad s=1, \dots, 4
\eeqa
and
\eqa
&&
f_2(z) =-\left[ \sum_{i=1}^3 \beta_2^{(i)} \frac{R'(u_i)}{z-u_i} + \beta_2^{(\infty)} z\right]
\nn\\
&&
f_3(z) = \sum_{i=1}^3 \beta_3^{(i)} \frac{{R'}^2(u_i)}{z-u_i} + \beta_3^{(\infty)} z^3.
\nn
\eeqa
Denote
\beq\label{spur4}
D=D(p,q,u) =\sum_{s=1}^4 \tilde p_s^{[2]} \frac{R(q_s)}{\Delta'(q_s)}
\eeq
the determinant of the linear system (\ref{sistema3}).
Then
\begin{eqnarray}\nn
&&
a=-\sum_{s=1}^4  \frac{w_s}{D}\left[
\frac12\sum_{j,k,l} {\rm sign}(s,j,k,l) p_j^{[2]}R(q_j)\frac{q_k-q_l}{W(q)}
+\frac{R(q_s)}{{\Delta'}^2(q_s)}+\right.\nn \\
&&
\quad\left.
+\frac1{\Delta'(q_s)}{\sum_{j=1}^4(q_s-q_j)^2\frac{R(q_j)}{{\Delta'}^2(q_j)}(\sigma_1-q_s-3 q_j)}\right]
\nn
\end{eqnarray}
\begin{eqnarray}\nn
&&
b=\sum_{s=1}^4 \frac{w_s }{D}\left[\frac12\sum_{j,k,l} {\rm sign}(s,j,k,l) p^{[2]}_jR(q_j)
\frac{q_k^2-q_l^2}{W(q)}+\right.\nn \\
&&
\quad\left.+\frac{R(q_s)(\sigma_1-2 q_s)}{{\Delta'}^2(q_s)}
+\frac1{\Delta'(q_s)}{\sum_{j=1}^4(q_s-q_j)^2\frac{R(q_j)}{\Delta'(q_j)^2}
(\sigma_1^2-\sigma_2-\sigma_1 q_s-2 \sigma_1 q_j+2 q_s q_j)}
\right]\nn \\
&&
\nn
\end{eqnarray}
\begin{eqnarray}\nn
&&
c=-\sum_{s=1}^4  \frac{w_s}{D}\left[\frac12{\sum_{j,k,l} {\rm sign}(s,j,k,l) p_j^{[2]}R(q_j)
q_k q_l\frac{q_k-q_l}{W(q)}}+\right.\nn \\
&&
\quad
+\frac1{\Delta'(q_s)}{\sum_{j=1}^4(q_s-q_j)^2\frac{R(q_j)}{{\Delta'}^2(q_j)} q_j
\left(\sigma_1^2-4\sigma_2-3 q_s^2+2\sigma_1 q_s+3 \frac{ \sigma_4}{q_sq_j}\right)}+
 \nn \\
&&
\quad\left.+\frac{R(q_s)}{{\Delta'}^2(q_s)}(3 q_s^2-2\sigma_1+2 \sigma_2)
\right]
\nn
\end{eqnarray}
\beq\label{spurh}
h=\frac1{D} \sum_{s=1}^4 \frac{w_s}{\Delta'(q_s)}
\eeq
where 
\beq\label{discr}
W(q) = \prod_{i<j}(q_i-q_j),
\eeq
${\rm sign}(s,j,k,l)$ is the sign of the permutation $\left(\begin{array}{cccc} 1 & 2 & 3 & 4 \\ s & j & k & l\end{array}\right)$.

Using (\ref{hamnostra}) one obtains the following expression for the Hamiltonians of Schlesinger equations $S_{(3,3)}$
\eqa\label{gamil}
&&
{\cal H}_i = -\frac1{R'(u_i)} \left[ \beta_2^{(\infty)} u_i + h -\sum_{s=1}^4 \tilde p_s^{[3]} \frac{R(q_s)}{u_i-q_s}
+\sum_{j\neq i} \frac{\beta_2^{(j)}}{u_i-u_j} R'(u_j)\right]
\nn\\
&&
+\frac12 \beta_2^{(i)} \frac{R''(u_i)}{{R'}(u_i)} + 3\, \frac{\Delta'(u_i)}{\Delta(u_i)} -\frac92 \, \frac{R''(u_i)}{R'(u_i)},
\quad i = 1, \, 2, \, 3
\eeqa
where the rational function $h=h(p,q,u)$ was defined in (\ref{spurh}). Clearly of these three Hamiltonians only one is independent: the solutions depend only on the combination $(u_3-u_1)/(u_2-u_1)$.


\setcounter{equation}{0}
\setcounter{theorem}{0}

\section{Comparison of spectral and isomonodromic coordinates}

\subsection{Spectral coordinates.}

We recall the construction of the algebro-geometric 
Darboux coordinates on the generic reduced symplectic leaves (\ref{red-leaf})
following the scheme of 
\cite{veselov, adams, DuDi, gekhtman}. 
We call these algebro-geometric Darboux coordinates {\it spectral 
coordinates}. 


The spectral coordinates are defined as follows. 
Let us assume that all the matrices ${ A}_i$ have
pairwise distinct nonzero eigenvalues $\lambda^{(i)}_1$, \dots, 
$\lambda^{(i)}_m$, and that the diagonal matrix
$$
{A}_\infty:=-({A}_1+\dots + {A}_{n})={\rm diag}\, (\lambda^{(\infty)}_1, \dots, \lambda^{(\infty)}_m)
$$
has distinct nonzero diagonal entries.
Consider the characteristic polynomial of the matrix 
$A(z)$ of 
the form
\begin{equation}\label{matrix}
A(z)= \sum_{i=1}^{n} {{A}_i\over z-u_i}
\end{equation}
with constant matrices ${A}_1$, \dots, ${A}_{n}$
satisfying the following properties. 
Denote
\begin{equation}\label{char}
{\cal R}(z,w) =\det (w-A(z))=w^m+\al_1(z) w^{m-1} + \dots +\al_m(z)
\end{equation}
the characteristic polynomial of the matrix $A(z)$. Denote by
\begin{equation}\label{discr}
D(z) := \left[  \prod_{i=1}^{n} (z-u_i)\right]^{m(m-1)}\prod_{i\neq
j}(w_i(z)-w_j(z))
\end{equation}
the discriminant of the polynomial ${R}(z,w)$. In this formula 
$w_1(z)$,\dots,$w_m(z)$ are roots of the equation ${R}(z,w)=0$. The resulting
expression is a polynomial in the coefficients $\al_1(z)$, \dots, $\al_m(z)$.
Under the
assumption that the matrix ${A}_\infty$ has simple spectrum the degree 
of the discriminant is equal to
\begin{equation}\label{degN}
N=m(m-1)(n-1).
\end{equation}

\noindent {\bf Assumption 1.} The $N$ roots of the discriminant are 
simple and pairwise distinct. Also we require that
\begin{equation}\label{dui}
D(u_i)\neq 0, ~~i=1, \dots, n.
\end{equation}

Due to this assumption the {\it spectral curve}
\begin{equation}\label{curve}
{\cal R}(z,w)=0
\end{equation}
of the matrix $A(z)$ is smooth outside the lines $z=u_1$, \dots, 
$z=u_{n}$, $z=\infty$. These lines intersect the spectral curve in 
singular points of multiplicity $m$.

Let us introduce the row vectors ${b}_0$, ${b}_1(z)$, $b_2(z)$, \dots by
\begin{equation}\label{bk}
b_0=(1, 0, \dots, 0), ~~b_k(z) = b_0 A^{k-1}(z), ~k>0.
\end{equation}
Denote $B(z)$ the $m\times m$ matrix with the rows $b_0$, $b_1(z)$, \dots,
$b_{m-1}(z)$
\begin{equation}\label{b-z}
B(z) = \left( \begin{matrix}b_0 \cr b_1(z) \cr \cdot \cr \cdot \cr \cdot \cr
b_{m-1}(z)\cr\end{matrix} \right).
\end{equation}

Put
\begin{equation}\label{Delta}
{\Delta_0}(z) = \left[ \prod _{i=1}^{n} (z-u_i)\right]^{m(m-1)\over 2} 
\det B(z).
\end{equation}

{\bf Assumption 2.} All the roots $\gamma_1$, \dots, $\gamma_g$ of the 
polynomial ${\Delta_0}(z)$ are pairwise distinct and the degree of 
${\Delta_0}(z)$ is equal to
the maximal value
\begin{equation}\label{degDelta}
g=\deg {\Delta_0}(z) ={1\over 2} m\, (m\,{n}-m-{n}-1)+1.
\end{equation}
We also assume that the roots $\gamma_1$, \dots, $\gamma_g$ do not coincide
with the poles $z=u_i$ of the Fuchsian system neither with the zeroes of the
discriminant $D(z)$. Under this assumption there exists, for any 
$i=1, \dots, g$, a unique, up to normalization, eigenvector
$\psi^{i}$ of the matrix $B(\gamma_i)$ with zero eigenvalue. The first 
component
of the eigenvector vanishes. It is also an eigenvector of the matrix
$A(\gamma_i)$ with some eigenvalue $\mu_i$,
$$
B(\gamma_i) \psi^i =0, ~~A(\gamma_i) \psi^i= \mu_i \psi^i, ~~\psi^i(\psi_1^i, \psi_2^i, \dots,
\psi_m^i)^T, ~~\psi_1^i=0, ~~i=1, \dots, g.
$$
Observe that the genus of the Riemann surface (\ref{curve}) is equal to $g$.
The spectral curve (\ref{curve}) together with the divisor 
\begin{equation}\label{divisor}
{\cal D} = \sum_{i=1}^g (\gamma_i, \mu_i)
\end{equation}
determines the matrix $A(z)$ uniquely up to a conjugation by a constant 
diagonal matrix.
Moreover, the matrices $A(z)$ satisfying the above assumptions form a Zariski
open subset in the space of all matrices of the form (\ref{matrix}). All these
facts are rather standard for the theory of algebraically completely
integrable systems. We give a sketch of proofs of these statements in the
Appendix.

\begin{df}
We call the set $(\gamma_1,\dots,\gamma_g,\mu_1,\dots,\mu_g)$ the 
{\it spectral coordinates}\/ on (\ref{red-leaf})
\end{df}

\begin{ex} For the $m=3$ case the polynomial $\Delta(z)$ determining the
isomonodromic coordinates reads
$$
\Delta(z) = R(z) \left[ A_{12} A_{13} \left( A_{22} - A_{33}\right)
-A_{12}^2 A_{23} +A_{13}^2 A_{32} -A_{12} A_{13}' +A_{12}' A_{13}\right]
$$
where, as usual $R(z)=\prod_i (z-u_i)$ and $A_{ij}=A_{ij}(z)$ are the entries
of the $3\times 3$ matrix $A(z) = \sum_i A_i/(z-u_i)$. The positions of the
spectral coordinates are determined by the polynomial
$$
\Delta_0(z) = R(z) \left[ A_{12} A_{13} \left( A_{22} - A_{33}\right)
-A_{12}^2 A_{23} +A_{13}^2 A_{32}\right].
$$
We see that, unlike the case $m=2$ (see above Example \ref{ex33}) 
starting from
$m=3$ the spectral and isomonodromic coordinates do not coincide.
\end{ex}


\subsection{Spectral coordinates versus isomonodromic coordinates}

In this subsection we prove that in certain semi--classical limit
the isomonodromic coordinates
$q_1$, \dots, $q_g$, $p_1$, \dots, $p_g$  converge to the 
algebro--geometric Darboux coordinates.

Let us consider the following family of Fuchsian systems depending on a
small parameter $\epsilon$
\begin{equation}\label{eps}
\epsilon \, {d\Phi\over d{z}} =A({z})\, \Phi, ~~ \Phi=(\phi_1({z}), 
\dots, \phi_m({z}))^T.
\end{equation}

\begin{thm} \label{teorem2} Under the assumptions of the Theorem \ref{skal-red} the apparent singularities of the
scalar reduction of the Fuchsian system admit the following expansion
\begin{equation}\label{qeps}
q_k=\gamma_k + O(\epsilon), ~~ \epsilon \to 0, ~~ k=1, \dots, g.
\end{equation}
Moreover, the scalar reduction can be written in the following form
\begin{equation}\label{scalar}
\epsilon^m y^{(m)} + \epsilon^{m-1} a_1(z, \epsilon)y^{(m-1)}
+\dots +a_m(z,\epsilon) y=0
\end{equation}
where the functions $a_1(z, \epsilon)$, \dots, $a_m(z, \epsilon)$ are analytic
in $(z, \epsilon)$ for $z\neq u_i$, $z\neq \gamma_j$, $|\epsilon |<<1$ and
\begin{equation}\label{scalar1}
a_k(z,\epsilon) = \alpha_k(z) - \epsilon\, \beta_k(z) +O(\epsilon^2), ~~k=1,
\dots, m
\end{equation}
where
\begin{equation}\label{scalar2}
\beta_k(z) ={1\over z-\gamma_j} \left[ \mu_j^{k-1} + \alpha_1(\gamma_j)
\mu_j^{k-2}+ \dots + \al_{k-1}(\gamma_j)\right] +O(1), ~~ z\to \gamma_j.
\end{equation}
In particular, the parameters $p_j$ defined in (\ref{pi}) have the following
expansion
\begin{equation}\label{peps}
p_j = \epsilon^{-1} [\mu_j+\alpha_1(\gamma_j)] +O(1), ~~ \epsilon \to 0, ~~ j=1, \dots, g.
\end{equation}
\end{thm}

Here $\gamma_k$, $\mu_k$ are the spectral coordinates of the matrix $A(z)$.

\begin{lm} \label{lemma1} The following formula holds true for any $k\geq 0$
\begin{equation}\label{k-th}
\epsilon^k y^{(k)} = [b_k(z) +\epsilon \, \tilde b_k(z) +O(\epsilon^2)]\, 
\Phi
\end{equation}
where the row vectors $b_k(z)$ were defined in (\ref{bk}) and the row vectors
$\tilde b_k(z)$ are defined by the following recursive procedure
\begin{equation}\label{bk-til}
\tilde b_0=0, ~~ \tilde b_{k+1}(z) = \tilde b_k(z) A(z) + b_k'(z), ~ k\geq 0.
\end{equation}
\end{lm}

Here and below it will be understood that the product of the row vector by the
column vector is a scalar.

\noindent{\bf Proof.} For $k=0$ (\ref{k-th}) is obvious. 
Since the first row of $A({z})$ is $b_1({z})$, the first equation 
of the Fuchsian system can be recast into the form
$$
\epsilon\, y' = b_1({z}) \,\Phi .
$$
This proves (\ref{k-th}) for $k=1$. Let us now assume (\ref{k-th}) for $k$ 
and prove it for $k+1$. Differentiating both sides of (\ref{k-th}) in 
${z}$ and multiplying by $\epsilon$ yields
$$
\begin{array}{ll}
\epsilon^{k+1} \, y^{(k+1)} =& 
\epsilon[b_k'({z})+\epsilon \tilde b_k'({z})+O(\epsilon^2)]\,
\Phi({z})+
[b_k({z})+\epsilon \, 
\tilde b_k({z}) +O(\epsilon^2)]\,A({z}) \Phi({z})=\\
&=[b_{k+1}({z})+\epsilon \,(\tilde b_k({z}) A({z})+ 
b_k'({z}))+O(\epsilon^2)]\,\Phi({z}).\\
\end{array}
$$
The proof of the Lemma is completed by induction. \epf

\begin{cor} Define $m\times m$ matrix valued function $\tilde B(z)$ with the
rows $\tilde b_0(z)$, $\tilde b_1(z)$, \dots, $\tilde b_{m-1}(z)$. Then the
scalar reduction of the Fuchsian system reads
\begin{equation}\label{red1}
\epsilon^m y^{(m)} =[b_m(z) + \epsilon \tilde b_m(z) + O(\epsilon^2)]
\left[ B(z) + \epsilon \tilde B(z) + O(\epsilon^2)\right]^{-1} \hat y
\end{equation}
where
\begin{equation}\label{yhat}
\hat y: = (y, \epsilon \, y', \dots, \epsilon^{m-1} y^{(m-1)})^T.
\end{equation}
\end{cor}

\noindent{\bf Proof. } From (\ref{k-th}) we obtain 
$$
\left(\begin{array}{c}
y\\ \epsilon \, y'\\ \dots\\ \epsilon^{m-1} y^{(m-1)}\\
\end{array}\right)=(B+\epsilon\tilde B+O(\epsilon^2))\, \Phi.
$$
This proves (\ref{red1}). \epf
\vskip 0.1 cm

From the Corollary it readily follows the claim of the Theorem about
expansions (\ref{qeps}) of the apparent singularities. It also follows
analyticity of the coefficients 
\begin{equation}\label{red2}
\begin{array}{l}
(a_m(z,\epsilon), a_{m-1}(z,\epsilon), \dots, a_1(z,\epsilon))=\\
=[b_m(z) + \epsilon \tilde b_m(z) + O(\epsilon^2)]
\left[ B(z) + \epsilon \tilde B(z) + O(\epsilon^2)\right]^{-1}=\\
=b_m(z) B^{-1}(z) + \epsilon\left[ \tilde b_m(z) B^{-1}(z) - b_m(z)
B^{-1}(z) \tilde B(z)\, B^{-1}(z)\right] +O(\epsilon^2)\\
\end{array}
\end{equation}
for small $\epsilon$ and for $z$ away from the poles $z=u_i$ and $z=\gamma_j$.

Let us now simplify the r.h.s. of the formula (\ref{red2}).

\begin{lm} \label{lemma2} The leading term in the r.h.s. of (\ref{red2})
reads
\begin{equation}\label{lead}
b_m(z)\, B^{-1}(z)=-(\al_m(z), \al_{m-1}(z), \dots, \al_1(z)).
\end{equation}
\end{lm}

\noindent{\bf Proof.}   Using Cayley - Hamilton theorem we obtain
$$
\begin{array}{ll}
b_m(z)\, B^{-1}(z)&=b_0 A^m(z) B^{-1}(z) =\\ 
&=-b_0 \left( \al_m(z)+\al_{m-1}(z) A(z)
+ \dots + \al_1(z) A^{m-1}(z)\right) B^{-1}(z)=\\
&=-\left[ \al_m(z) b_0 +\al_{m-1}(z) b_1(z) + \dots + 
\al_1(z) b_{m-1}(z)\right]
B^{-1}(z) =\\
&=-(\al_m(z), \al_{m-1}(z), \dots, \al_1(z)).\\
\end{array}
$$
In the last line we use the definition of the inverse matrix $B^{-1}(z)$. The
lemma is proved. \epf

We will now simplify the linear in $\epsilon$ term. We need the following simple

\begin{lm} \label{lemma3} Let us introduce matrix $T(z)$ by 
\begin{equation}\label{matt}
T(z) = \left(\begin{array}{ccccc}
0 & 1 & 0 & \dots & 0 \\
0 & 0 & 1 & \dots & 0 \\
  & \dots &  & \dots & \\
0 & 0  &  & \dots & 1 \\ 
  -\al_m(z) & -\al_{m-1}(z) & & \dots & -\al_1(z) \\
\end{array}\right).
\end{equation}
Then
\begin{equation}\label{id}
T(z) B(z) = B(z)A(z) .
\end{equation}
\end{lm}

\noindent{\bf Proof.}   Use the definition of the matrix $B(z)$ and 
Cayley - Hamilton theorem. \epf

\begin{lm} \label{lemma4} The following identity holds true
\begin{eqnarray}\label{long1}
&&
\tilde b_m(z) B^{-1}(z) - b_m(z) B^{-1}(z) \tilde B(z) B^{-1}(z)
\nn\\
&& 
=\sum_{k=1}^m e_k B'(z) B^{-1}(z) \left[ T^{m-k}(z) + \al_1(z) T^{m-k-1}(z)
+\dots + \al_{m-k}(z)\right].
\end{eqnarray}
Here $e_1$, \dots, $e_m$ are row vectors of the standard basis in
${{\mathbb C}^m}^*$, $(e_k)_i=\delta_{i\, k}$.
\end{lm}

\noindent{\bf Proof.}   Using arguments of Lemma \ref{lemma2} we replace 
$b_m(z) B^{-1}(z)$ by
$-(\al_m(z), \dots, \al_1(z))$. Next, using induction we derive the following
formula for the row vectors (\ref{bk-til})
$$
\tilde b_k(z) = e_2 B'(z) A^{k-2}(z) + e_3 B'(z) A^{k-3}(z)+ \dots + e_k B'(z),
~~ k\geq 2.
$$
Using identity (\ref{id}) the last formula can be recast into the form
$$
\tilde b_k(z) = \left[ e_2 B'(z) B^{-1}(z) T^{k-2}(z) + e_3 B'(z) B^{-1}(z)
T^{k-3}(z)\right.
$$
$$
\left. + \dots + e_k B'(z) B^{-1}(z) \right] B(z), ~~ k\geq 2.
$$
This implies the formula (\ref{long1}) with the summation in the r.h.s. starting
from $k=2$. Since the first row of the matrix $B'(z)$ identically vanishes,
adding the term with $k=1$ does not change the sum. \epf

We now want to compute the Laurent expansion of the coefficients of the scalar
reduction (\ref{scalar}) at the points $z=\gamma_j = q_j +O(\epsilon)$. 
Let $\gamma$ be one of the zeroes of ${\Delta_0}(z)$, $\mu$ the eigenvalue
of the matrix $A(\gamma)$ such that the corresponding eigenvector $\psi=(\psi_1,
\psi_2, \dots, \psi_m)^T$
$$
A(\gamma) \psi = \mu\, \psi
$$
satisfies
$$
\psi_1=0.
$$
According to our assumptions the eigenvector is defined uniquely up to a scalar
factor. As we know, $\psi$ is also an eigenvector of the matrix $B(\gamma)$ 
with
zero eigenvalue,
$$
B(\gamma)\psi=0.
$$
Moreover, there exists an analytic function $\lambda(z)$ defined for
$|z-\gamma|<<1$ being an eigenvalue of $B(z)$ s.t. $\lambda(z)$ has a simple 
zero
at $z=\gamma$, and $\lambda(z)$ does not coincide with other eigenvalues of
$B(z)$. Denote $\psi(z)= (\psi_1(z), \dots, \psi_m(z))^T$ the analytic vector
valued function s.t.
$$
B(z) \psi(z) = \lambda(z) \psi(z), ~~ |z-\gamma|<<1,
$$
$$
\lambda(\gamma)=0, ~~ \psi(\gamma)=\psi.
$$
We also introduce left eigenvector $\psi^*(z) = (\psi_1^*(z), \dots,
\psi_m^*(z))$,
\begin{equation}\label{left}
\psi^*(z)B(z) = \lambda(z) \psi^*(z), ~~|z-\gamma|<<1.
\end{equation}
Denote
$$
\psi^*:= \psi^*(\gamma).
$$

\begin{lm} \label{lemma5} $\psi^*$ is a left eigenvector of $T(\gamma)$ with
the eigenvalue $\mu$,
\begin{equation}\label{star0}
\psi^* T(\gamma) = \mu \psi^*.
\end{equation}
In particular, it
can be chosen in the form
\begin{equation}\label{star}
\psi^*_k= \mu^{m-k} + \al_1(\gamma) \mu^{m-k-1} + \dots + \al_{m-k}(\gamma), ~~
k=1, \dots, m.
\end{equation}
\end{lm}

\noindent{\bf Proof.}   From (\ref{left}) for $z=\gamma$ it follows
\begin{equation}\label{star1}
\psi^* B(\gamma)=0.
\end{equation}
Let us prove that $\psi^* T(\gamma)$ is again a left eigenvector of $B(\gamma)$
with zero eigenvalue. Indeed, using (\ref{id}) we obtain
$$
\psi^* T(\gamma) B(\gamma) = \psi^* B(\gamma) A(\gamma) = 0.
$$

The eigenvector of $T(\gamma)$ with the eigenvalue $\mu$ can be written
in the form (\ref{star}). Let us prove that this eigenvector satisfies 
(\ref{star1}). 

According to our assumptions all the eigenvalues of the matrix $A(\gamma)$ are
pairwise distinct. Of course, they coincide with the eigenvalues of the matrix
$T(\gamma)$. Therefore it suffices to prove that
$$
\psi^* B(\gamma) \psi'=0
$$
for an arbitrary eigenvector $\psi'$ of the matrix $A(\gamma)$,
$$
A(\gamma) \psi' = \mu' \psi'.
$$
Indeed, from (\ref{star}) we obtain
$$
\psi^* B(\gamma) \psi' = b_0 \psi' \, \hat R(\mu, \mu',\gamma)
$$
where
$$
\hat  R(\mu, \mu',\gamma) = {{\cal R}(z,w)-{\cal R}(z,w')\over w-w'}|_{z=\gamma, ~
w=\mu,~w'=\mu'}
$$
for $\mu'\neq \mu$ and
$$
\hat{\cal R}(\mu, \mu, \gamma) ={\pal{\cal R}(z,w)\over 
\pal w}|_{z=\gamma, ~ w=\mu}.
$$
It is clear that $\hat  R(\mu, \mu', \gamma)=0$ for $\mu'\neq \mu$. So
$\psi^* B(\gamma) \psi'=0$. For $\mu'=\mu$ we have 
$\hat  R(\mu, \mu, \gamma)\neq 0$
(since $\gamma$ is not a zero of the discriminant $D(z)$) but $b_0\, \psi=0$
since $\psi_1=0$. The Lemma is proved. \epf

We will now compute the leading term of the Laurent expansion of the 
logarithmic
derivative $B'(z) B^{-1}(z)$ at $z\to \gamma$.

\begin{lm} \label{lemma6} For $z\to \gamma$
\begin{equation}\label{logder}
B'(z) B^{-1}(z) = {1\over z-q} {B'(\gamma)\psi \otimes \psi^*\over \psi^*
B'(\gamma)\psi} +O(1).
\end{equation}
\end{lm}

\noindent{\bf Proof.}   Let 
\begin{equation}\label{proj}
\Pi_1(z)={\psi(z) \otimes \psi^*(z)\over \psi^*(z) \psi(z)}
\end{equation}
be the projector of ${\mathbb C}^m$ onto the direction of the eigenvector
$\psi(z)$ parallel to the $(m-1)$-dimensional subspace spanned by other
eigenvectors. Denote $\Pi_2(z) = {\rm id} - \Pi_1(z)$ the complementary
projector and put
$$
B_2(z) := B(z) \Pi_2(z).
$$
We have
$$
B(z) = \lambda(z) \Pi_1(z) + B_2(z).
$$
All these matrix valued functions are analytic for $z$ sufficiently close to
$\gamma$ and since $B(z)$ has a unique zero eigenvalue, 
$$
{\rm rank}\left( B_2(\gamma)\right)=m-1,
$$
and the image of $B_2({z})$ is transverse to that of $\Pi_1$.
So
$$
B'(z) B^{-1}(z) = (\log\lambda(z))' \Pi_1(z) + \lambda^{-1}(z) B_2'(z) \Pi_1(z)
+ {\rm regular ~ terms}.
$$
Since $B_2(z) \Pi_1(z) \equiv 0$, we obtain
$$
B_2'(z) \Pi_1(z) = - B_2(z) \Pi_1'(z) =-{B_2(z) \psi'(z)\otimes \psi^*(z)\over
\psi^*(z) \psi(z)} + \dots
$$
$$
 =  -{B(z) \psi'(z)\otimes \psi^*(z)\over
\psi^*(z) \psi(z)} + \dots
$$
In the last equation dots denote terms analytic at $z=\gamma$. We obtain
$$
B'(z) B^{-1}(z) = {1\over z-\gamma} \left[ {\psi \otimes \psi^*\over \psi^*
\psi} - {1\over \lambda'(\gamma) }{B(\gamma) \psi'(\gamma) \otimes \psi^*\over
\psi^* \psi}\right] + O(1).
$$
Using the well-known formula of the ``perturbation theory''
$$
\lambda'(\gamma) ={\psi^* B'(\gamma)\psi\over \psi^* \psi}
$$
(observe that the formula implies $\psi^* B'(\gamma)\psi\neq 0$) and the
identity
$$
B(\gamma) \psi'(\gamma) +B'(\gamma)\psi = \lambda'(\gamma) \psi
$$
we complete the proof of the Lemma. \epf

\noindent{\bf End of the proof of the Theorem \ref{teorem2}}. We are to compute the sum
$$
\sum_{k=1}^m e_k B'(z) B^{-1}(z) \left[ T^{m-k}(z) + \al_1(z) T^{m-k-1}(z) +
\dots + \al_{m-k}(z)\right]
$$
$$
= {1\over z-\gamma} \sum_{k=1}^m e_k {B'(\gamma)\psi\otimes \psi^*\over \psi^*
B'(\gamma)\psi}
\left[ T^{m-k}(\gamma) + \al_1(\gamma) T^{m-k-1}(\gamma) +
\dots + \al_{m-k}(\gamma)\right] +O(1)
$$
$$
= {1\over z-\gamma} \sum_{k=1}^m e_k {B'(\gamma)\psi\otimes \psi^*\over \psi^*
B'(\gamma)\psi}
\left[ \mu^{m-k} + \al_1(\gamma) \mu^{m-k-1} +
\dots + \al_{m-k}(\gamma)\right] +O(1)
$$
$$
={1\over z-\gamma} \sum_{k=1}^m {\psi_k^* \left( B'(\gamma)\psi\right)_k\over
\psi^* B'(\gamma)\psi} \, \psi^* +O(1) ={1\over z-\gamma} \psi^* +O(1)
$$
where the left eigenvector $\psi^*$ is chosen in the form (\ref{star}). The
Theorem is proved. \epf


\subsection{Canonical transformations.}
Let us consider the isomonodromic coordinates $(q_1,\dots,q_g,p_1,
\dots,p_g)$ obtained from the scalar reduction w.r.t. the first row 
of (\ref{N1}).

\begin{prop}\label{prsimmk}
For every $k=2,\dots,n$, the following 
transformation
\begin{equation}
S_k:\quad \left\{\begin{array}{l} 
\tilde q_i=u_1+ u_k-q_i, \quad i=1,\dots,g,\\
 \tilde p_i=-p_i,\quad i=1,\dots,g,\\
 \tilde u_l = u_1+u_k-u_l,\quad l=1,\dots,n, \\
\tilde \lambda^{(k)}_j= 
\lambda^{(1)}_j,\quad j=1,\dots,m,\\
\tilde\lambda^{(1)}_j= 
\lambda^{(k)}_j,\quad j=1,\dots,m,\\
\tilde {\mathcal H}_l=-{\mathcal H}_l,\quad l=1,\dots,n,
\end{array}\right.
\label{simmk}\end{equation}
is a birational canonical transformation of the Schlesinger systems.
This transformation acts on the monodromy matrices 
as follows
\begin{eqnarray}\label{ss-inf3}
&&
\tilde  M_1=M_1^{-1}\dots M_{k-1}^{-1} M_k M_{k-1}\dots M_1, \nn\\
&&
\tilde M_j=M_{j},\quad j\neq 1,k,\\
&&
\tilde M_k=M_{k-1} \dots M_2 M_1 M_2^{-1} \dots M_{k-1}^{-1},\quad i=k+1,\dots,n.\nn
\end{eqnarray}
\end{prop}

\noindent{\bf Proof. }  The transformation (\ref{simmk}) is obviously birational. To show 
that (\ref{simmk}) is a canonical transformation of the Schlesinger systems, we just
observe that it is obtained by the conformal transformation $\zeta=u_1+u_k-{z}$
of the scalar reduction (\ref{eqscal}). In fact (\ref{eqscal}) 
is transformed to
$$
{{\rm d}^m y\over{\rm d}\zeta^m}=\sum_{l=0}^{m-1}
\tilde d_l(\zeta){{\rm d}^l y\over{\rm d}\zeta^l},
$$
where $\tilde d_l(\zeta)=(-1)^{l-m}d_l(u_1+u_k-\zeta){\tilde f_{l}(\zeta)\over\tilde\Delta(\zeta)\tilde R(\zeta)^{m-l}}$
with $\tilde R(\zeta)= \prod_{k=1}^{n}(\zeta-\tilde u_k)$, 
$\Delta(\zeta)= \prod_{i=1}^{g}(\zeta-\tilde q_i)$ and
$\tilde f_{l}(\zeta)= (-1)^{(m-l)(n-1)+g} f_{l}(u_1+u_k-\zeta)$. 
To obtain $\tilde p_i$ we use formula
(\ref{pi}):
$$
\begin{array}{ll}
\tilde p_i=&\res_{\zeta=\tilde q_i}\left(\tilde d_{m-2}(\zeta)+
{\tilde d_{m-1}(\zeta)\over\zeta-\tilde q_i} \right)=\\
&=\res_{\zeta=\tilde q_i}\left(d_{m-2}(u_1+u_k-\zeta)-
{d_{m-1}(u_1+u_k-\zeta)\over\zeta-\tilde q_i}\right)=-p_i.\\
\end{array}
$$
To obtain the formulae for the exponents, just observe that the conformal 
transformation $\zeta=u_1+u_k-{z}$ permutes $u_1$ with $u_k$. 

Let us now prove the formula (\ref{ss-inf3}). The involution 
$i_k:u_1\leftrightarrow u_k$ changes the basis in the fundamental group 
$\pi_1(\overline\complessi\backslash\{u_1,\dots,u_n,\infty\})$. In fact,
as explained in Section 2.2, the cuts $\pi_1,\dots,\pi_n$ along which we 
take our basis $l_1,\dots,l_n$, are ordered according to the 
order of the poles. Applying the transformation $i_k$ we then arrive 
at the new basis of loops
\begin{eqnarray}
&&
l_1'=l_1l_2\cdots l_{k-1}l_k l_{k-1}^{-1}\cdots l_2^{-1}l_1^{-1},\nn\\
&&
l_j'=l_{j},\quad j\neq 1,k,\nn\\
&&
l_k'=l_{k-1}^{-1}\cdots l_2^{-1} l_1 l_2\cdots l_{k-1}.\nn
\end{eqnarray}
from these formulae we immediately obtain (\ref{ss-inf3}).
\epf

\begin{prop}
The following transformation
\begin{equation}
S_\infty:\quad \left\{\begin{array}{l} 
\tilde q_i={1\over q_i-u_1}, \quad i=1,\dots,g,\\
\tilde p_i= -p_iq_i^2-{2m^2-1\over m}q_i,\quad i=1,\dots,g,\\
\tilde u_l= {1\over u_l-u_1},\quad l=2,\dots,n, \\
u_{1}\mapsto\infty,\\
\infty\mapsto u_{1},\\
\tilde\lambda^{(\infty)}_1= 
\lambda^{(1)}_1+{m-1\over m},\\
\tilde\lambda^{(\infty)}_j= 
\lambda^{(1)}_j-{1\over m},\quad j=2,\dots,m,\\
\tilde\lambda^{(1)}_1=\lambda^{(\infty)}_1-
{m-1\over m},\\
\tilde\lambda^{(1)}_j=\lambda^{(\infty)}_j+
{1\over m},\quad j=2,\dots,m,\\
\tilde{\mathcal H}_1=\mathcal H_1,\\
\tilde{\mathcal H}_l= -\mathcal H_l (u_l-u_1)^2 + (u_l-u_1) (d^0_{m-1}(u_l-u_1))^2 -\\
\qquad-(u_l-u_1) {(m-1)(m^2-m-1)\over m}d^0_{m-1}(u_l-u_1),\quad l=2,\dots,n\\
\end{array}\right.
\label{simminf}\end{equation}
where 
$$
d^0_{m-1}(u_k) = \sum_{s=1}^g {1\over u_k-q_s}  - {m\, (m-1)\over 2} 
\sum_{l\neq k} {1\over u_k-u_l}.
$$
is a birational canonical transformation of the Schlesinger systems.
This transformation acts the monodromy matrices $M_1,\dots,M_n,M_\infty$
as follows
\begin{eqnarray}\label{ss-inf1}
&&
\tilde M_\infty=e^{-{2\pi i \over m}}  C_1 M_\infty^{-1}  M_1  M_\infty C_1^{-1},\nn\\
&&
\tilde M_1= e^{2\pi i \over m} C_1 M_\infty C_1^{-1},\\
&&
\tilde M_j=C_1^{-1}M_j C_1,\quad j\neq 1,\infty,\nn
\label{last}\end{eqnarray}
where $C_1$ is the connection matrix of $M_1$.
\label{prsimminf}\end{prop}

\noindent{\bf Proof. } The fact that the above transformation is birational 
is trivial. To show 
that it is a canonical transformation of the Schlesinger systems, we just
observe that it is obtained by a conformal transformation 
$\zeta={1\over{z-u_1}}$ and a gauge transformation 
$y=g(\zeta)\tilde y$, $g(\zeta)=\zeta^{m-1\over m}$ of the scalar 
reduction (\ref{eqscal}). In fact (\ref{eqscal}) is transformed to
$$
{{\rm d}^m\tilde y\over{\rm d}\zeta^m}=-\sum_{p=1}^{m-1}
\left(\begin{array}{c}m\\ p\\ \end{array}\right)
{1\over g(\zeta)}{{\rm d}^p g\over{\rm d}\zeta^p}
{{\rm d}^{m-p}\tilde y\over{\rm d}\zeta^{m-p}}
+\sum_{s=0}^{m-1}\hat d_s \sum_{p=0}^{s}
\left(\begin{array}{c}s\\ p\\ \end{array}\right)
{1\over g(\zeta)}{{\rm d}^p g\over{\rm d}\zeta^p}
{{\rm d}^{s-p}\tilde y\over{\rm d}\zeta^{s-p}},
$$
where
$$
\begin{array}{l}
\hat d_0=(-1)^m \zeta^{-2m}d_0({1\over\zeta}+u_1),\\
\hat d_s= (-1)^{m+1} c_{s+1}^{m+1}\zeta^{s-m}+
(-1)^m \sum_{l=s}^{m-1}\zeta^{l+s-2m}c_{s+1}^{l+1}d_l({1\over\zeta}+u_1),\\
\end{array}
$$
and $c_i^j:=(-1)^{j-1}(j-i)!\left(\begin{array}{c}j-2\\i-2\\ 
\end{array}\right) \left(\begin{array}{c}j-1\\i-1\\ 
\end{array}\right)$. Using the above formula and (\ref{pi}) it is 
a straightforward computation to obtain the formulae for $\tilde q_i$, $\tilde p_i$ and 
$\tilde{\mathcal H}_l$ in  (\ref{simminf}). 
The transformation law of the exponents is obtained in two stages: first the
conformal transformation $\zeta={1\over{z}}+u_1 $ maps 
\begin{eqnarray}
&&\lambda^{(\infty)}_1\to \lambda^{(1)}_1,\quad
\lambda^{(\infty)}_i+1\to \lambda^{(1)}_i,\quad i=2,\dots,m,\nn\\
&&
\lambda^{(1)}_1\to\lambda^{(\infty)}_1,\quad
\lambda^{(1)}_i \to\lambda^{(\infty)}_i+1,\quad i=2,\dots,m,\nn
\end{eqnarray}
then the gauge transformation 
$y=g(\zeta)\tilde y$, $g(\zeta)=\zeta^{m-1\over m}$ adds ${m-1\over m}$ to all 
exponents at infinity and subtracts the same quantity to all exponents at $0$. 
To show (\ref{ss-inf1}) we proceed as in the previous proof: the involution 
$i_\infty:u_1\leftrightarrow \infty$ changes the base point of the fundamental group (this is obtained by conjugating all monodromy matrices with the connection matrix $C_1$ of $M_1$), and it changes the  basis of loops as in the previous proof with $k$  replaced by $\infty$ and $k-1$ by $n$. 
This implies immediately (\ref{ss-inf1}). \epf

\begin{rmk}
Obviously we can obtain analogous birational canonical transformations
acting on the isomonodromic coordinates obtained from the scalar 
reduction w.r.t. any row of (\ref{N1}). 
\end{rmk}

\begin{rmk}
Apart from the above symmetries, there are other birational canonical 
transformations. In fact let us denote by $(q_1^{(j)},\dots,q_g^{(j)},
p_1^{(j)},\dots,p_g^{(j)})$ the isomonodromic coordinates obtained 
from the scalar reduction w.r.t. the $j$-th row. The transformation
that maps $(q_1^{(j)},\dots,q_g^{(j)},p_1^{(j)},\dots,p_g^{(j)})$
to the isomonodromic coordinates obtained from the scalar reduction 
w.r.t. the $i$-th row$(q_1^{(i)},\dots,q_g^{(i)},p_1^{(i)},\dots,p_g^{(i)})$
is by construction a birational canonical transformation. These 
transformations are the analogues of Okamoto's
$w_3$ for the Painlev\'e sixth equation (see \cite{Ok2}).
\end{rmk}

\noindent{\bf Acknowledgments}
The authors are very grateful to A. Bolibruch for many helpful 
conversations. This work is
partially supported by European Science Foundation Programme ``Methods of
Integrable Systems, Geometry, Applied Mathematics" (MISGAM), Marie Curie RTN ``European Network in Geometry, Mathematical Physics and Applications"  (ENIGMA), 
and by Italian Ministry of Universities and Researches (MIUR) research grant PRIN 2004
``Geometric methods in the theory of nonlinear waves and their applications". The researches of M.M. are also supported by EPSRC,
SISSA, ETH, IAS, and IRMA (Strasbourg).

\setcounter{equation}{0}
\setcounter{theorem}{0}


\def\thetheorem{A.\arabic{theorem}}
\def\theprop{A.\arabic{prop}}
\def\thelemma{A.\arabic{lemma}}
\def\thecor{A.\arabic{cor}}
\def\theexam{A.\arabic{exam}}
\def\theremark{A.\arabic{remark}}
\def\theequation{A.\arabic{equation}}

\appendix
\makeatletter
\renewcommand{\@seccntformat}[1]{{Appendix:}\hspace{-2.3cm}}
\makeatother
\renewcommand{\thesection}{Appendix:}
\section{\quad\qquad \ \ Algebro-geometric Darboux coordinates}


Here we outline the construction of the so-called algebro-geometric Darboux coordinates. Our
presentation follows \cite{DuDi}. However, the idea of constructing canonical coordinates for integrable systems by using the projections of the points in the divisor of a suitable normalized line--bundle  on the
spectral curve appeared for the first time in a paper by H. Flaschka and D.W. McLaughlin
\cite{flaschka} were the special cases of the Toda system and KdV equation were dealt with. Later
S.P.Novikov and A.P.Veselov \cite{veselov} generalized the construction to any hyperelliptic spectral curve and introduced a 
general class of finite-and infinite-dimensional Poisson brackets. 
The Flaschka--McLaughlin construction was then generalized to generic rational Lax pairs by M.R. Adams, J. Harnad and J. Hurtubise \cite{adams}.
A quantum version of this method was
initiated by E.Sklyanin \cite{sklyanin}. 

A construction of the polynomial ${\Delta_0}(z)$ equivalent to ours was given
in \cite{scott, gekhtman}. Our Theorem \ref{interpol} that enables to
constructing {\it rational} Darboux coordinates on the reduced symplectic leaves seems
to be new (cf. however the recent paper \cite{babel} where a similar approach to
constructing the Darboux coordinates was developed).

Let us rewrite the characteristic polynomial (\ref{char}) of the matrix
\begin{equation}\label{comden1}
A(z) = \sum_{i=1}^{n} {{ A}_i\over z-u_i} = {\hat A(z)\over 
\prod_{i=1}^{n} (z-u_i)}
\end{equation}
\begin{equation}\label{comden2}
\hat A(z) = -{ A}_\infty z^{{n}-1} + O(z^{{n}-2}), \quad
z\to\infty
\end{equation}
in the form
\begin{equation}\label{comden3}
\hat {\cal R}(z, \hat w) = \det \left( \hat w - \det \hat A(z) \right) = 
\left(
\prod_{i=1}^{n} (z-u_i)\right) ^m {\cal R}(z, w), ~~ \hat w = w \, 
\prod_{i=1}^{n} (z-u_i).
\end{equation}
Expanding the determinant one obtains a polynomial
\begin{equation}\label{comden4}
\hat {\cal R}(z, \hat w)=\hat w^m + \hat\alpha_1(z) \hat w^{m-1}+\dots +
\hat\alpha_m(z)=\sum_{i+({n}-1) j\leq m({n}-1)} a_{ij} z^i \hat w^j
\end{equation}
where
\begin{equation}\label{comden5}
a_{ij}= a_{ij}({ A}_1, \dots, { A}_{n}; 
u_1, \dots, u_{n})
\end{equation}
are some polynomials in the entries of the matrices 
${ A}_k$ and in $u_l$,
$$
a_{0m}=1, ~~ \hat\alpha_s(z) = 
\sum_{i=0}^{s\, ({n}-1)} a_{i\,, m-s}z^i, ~~s=0, \,
1, \dots, m.
$$
It is well known that, for a Zariski open subset in the projective
space  ${\bf P}^{M-1}$,
\begin{equation}\label{comden6}
M= ({n}-1){m(m+1)\over 2} + m+1
\end{equation}
with the homogeneous coordinates 
$$
a_{ij}, ~~i+({n}-1) j \leq m({n}-1)
$$
the affine algebraic curve
\begin{equation}\label{comden7}
\hat{\cal  R}(z, \hat w)=\sum_{i+({n}-1) j\leq m({n}-1)} a_{ij} z^i \hat w^j=0
\end{equation}
is smooth. Indeed, it suffices to check smoothness of one of the curves of the
above family, e.g. of
$$
\hat w^m = z^{m ({n}-1)} -1.
$$
Under the smoothness assumption the standard compactification of (\ref{comden7}) gives a compact
Riemann surface $\Gamma$ of the genus
\begin{equation}\label{genus}
g = ({n}-1){m(m-1)\over 2} -m+1
\end{equation}
(cfr. the formula (\ref{degDelta})). The infinite part of $\Gamma$ is a divisor
$D_\infty$ of the degree $m$. If the normal Jordan form of the matrix ${\cal
A}_\infty$ contains $k$ Jordan blocks of the multiplicities $m_1$, \dots, $m_k$
then the divisor $D_\infty$ has the form
\begin{equation}\label{d-inf}
D_\infty = m_1 \infty_1 + \dots + m_k \infty_k.
\end{equation}
Here $\infty_1, \dots, \infty_k \in \Gamma$ are the points added at infinity.
In particular, if the spectrum of ${A}_\infty$ is simple then the divisor
$D_\infty$ is a sum of $m$ distinct points $\infty_1$, \dots, $\infty_m$:
\begin{equation}\label{inf}
\infty_k := \{ z\to \infty, ~ \hat w\to \infty, ~ 
{\hat w\over z^{{n}-1}} \to
-\lambda^{(\infty)}_k\} , ~~ k=1, \dots, m.
\end{equation}
Here $\lambda_1^{(\infty)}$, \dots, $\lambda_m^{(\infty)}$ are the eigenvalues 
of the
matrix ${ A}_\infty$. 

Let us give an intrinsic characterization
of algebraic curves of the form (\ref{comden7}).
      
The coordinate functions $z$ and $\hat w$ have poles at the divisors $D_\infty$
and $({n}-1)D_\infty$ respectively. Conversely, the following simple statement
can be proved by using standard arguments based on the Riemann - Roch theorem

\begin{lm}\label{lemma1.1} Let $\Gamma$ be a Riemann surface of the 
genus (\ref{genus}). Let
$D_\infty$ be a divisor of the degree $m$ on $\Gamma$ such that 
$$
\dim H^0(\Gamma, {\cal O}(D_\infty))=2, ~~ 
\dim H^0(\Gamma, {\cal O}(({n}-1)D_\infty))={n}+1.
$$
Then the Riemann surface can be represented in the form (\ref{comden7}).
\end{lm}

\noindent {\bf Proof.} The first of the assumptions implies existence 
of a non-constant meromorphic function
$z$ with poles at the points of the divisor $D_\infty$. The second one yields
existence of another function $\hat w$ with poles at $(n-1)D_\infty$ that 
cannot
be represented as a polynomial in $z$. Let us now consider the space
$H^0(\Gamma,{\cal O}\left( m(n-1)D_\infty\right))$. The $M$ monomials
\begin{equation}\label{span}
z^i \hat w^j, \quad i+(n-1) j \leq m(n-1)
\end{equation}
belong to this space. Let us prove that these monomials are linearly dependent.
To this end let us compute the dimension
$$
\dim H^0(\Gamma,{\cal O}\left( m(n-1)D_\infty\right)).
$$
First of all, the degree of the divisor ${\cal D}:= m(n-1)D_\infty$
equals 
$$
\deg {\cal D} = m^2 (n-1) =2g-2+m(n+1) >2g-2.
$$
So the Riemann - Roch theorem gives
$$
\dim H^0(\Gamma,{\cal O}\left( m(n-1)D_\infty\right))=\deg {\cal D}-g+1
=M-1.
$$
This proves linear dependence of $M$ functions of the form (\ref{span}). The
Lemma is proved.\epf

Observe that, for ${n} >1$ the eigenvalues 
$\lambda_1^{(\infty)}$,
\dots, $\lambda_m^{(\infty)}$ are determined  by the Riemann surface uniquely
up to permutations and common affine transformations
$$
\lambda_i^{(\infty)} \mapsto a\lambda_i^{(\infty)} +b, ~~i=1, \dots, m.
$$
We will say that our Riemann surface $\Gamma$ is $D_\infty$-{\it generic} if 
the eigenvalues are pairwise distinct,
$$
\lambda_i^{(\infty)} \neq \lambda_j^{(\infty)} , ~~i\neq j.
$$

Instead of using the coefficients $a_{ij}$ as the homogeneous coordinates 
in the space of 
algebraic curves (\ref{comden7}) we will construct another system of
coordinates on a Zariski open subspace in ${\bf P}^{M-1}$. Let us choose 
${n}+g$
pairwise distinct numbers $u_1$, \dots, $u_{n}$, $\gamma_1$, \dots, 
$\gamma_g$ and also
${n}+1$
$m$-tuples of pairwise distinct numbers $\lambda_1^{(i)}$, \dots, 
$\lambda_m^{(i)}$ for every $i=1, \dots, n, \infty$ 
$$
\lambda_r^{(i)}\neq \lambda_s^{(i)}, ~~s\neq r, ~~i=1, \, 2, \dots, m, \, 
\infty
$$
satisfying the constraint
\begin{equation}\label{f-cons}
\sum_{i=1}^{n} \sum_{r=1}^m \lambda_r^{(i)} + 
\sum_{r=1}^m \lambda_r^{(\infty)}=0.
\end{equation}
Finally, choose arbitrary $g$ numbers $\mu_1$, \dots, $\mu_g$.
Denote
\begin{equation}\label{muhat}
\hat \lambda_r^{(i)} := \prod_{j\neq i} (u_i-u_j) \lambda_r^{(i)} , ~~r=1,
\dots, m, ~~i=1, \dots, {n}
\end{equation}
\begin{equation}\label{phat}
\hat \mu_s = \prod_{i=1}^{n} (\gamma_s - u_i)\mu_s, ~~s=1, \dots, g.
\end{equation}

\begin{thm} \label{interpol} For generic values of the parameters 
\begin{equation}\label{param}
u_1, \dots, u_{n}, \quad 
\lambda_1^{(i)}, \dots, \lambda_m^{(i)}, \quad i=1, \dots, n, \infty, 
\quad \gamma_1, \dots, \gamma_g,\, \mu_1, \dots, \mu_g
\end{equation} 
satisfying the constraint (\ref{f-cons})
there exists a unique 
curve $\hat {\cal R}(z, \hat w)=0$ of the form
(\ref{comden7}) with $a_{0m}=1$ passing through the points 
\begin{equation}\label{pointu}
(u_i, \hat\lambda_r^{(i)}), ~r=1, \dots, m, ~i=1, \dots, {n}, 
\end{equation}
\begin{equation}\label{pointinf}
z, \, \hat w \to \infty, ~ {\hat w\over z^{{n}-1}}\to - 
\lambda_r^{(\infty)}, ~
r=1, \dots, m
\end{equation}
\begin{equation}\label{pointq}
(\gamma_s, \hat \mu_s), ~~ s=1, \dots, g.
\end{equation}
\end{thm}

\noindent{\bf Proof.} Let us denote 
$$
\sigma_k^{(i)}:= \sigma_k(\lambda_1^{(i)}, \dots, \lambda_m^{(i)}), \quad k=1, \dots, m
$$
the value of the $k$-th elementary symmetric function of 
$\lambda_1^{(i)}, \dots, \lambda_m^{(i)}$, $i=1, \dots, n, \infty$.
The equation of the algebraic curve must have the form
\begin{eqnarray}\label{ansatz}
&&
\hat w^m -R(z) \sum_{i=1}^{n} {\sigma_1^{(i)}\over z-u_i} \hat w^{m-1}
\\
&&
+R(z) \sum_{k=2}^m \left[ (-1)^k \sum_{i=1}^{n} 
{\sigma_k^{(i)} \over z-u_i}[R'(u_i)]^{k-1} + \sigma_k^{(\infty)} z^{k\,
n-n-k}+ p_{k\, n - n - k -1}(z)\right] \hat w^{m-k} =0.\nn
\end{eqnarray}
Here, as above
$$
R(z) := \prod_{i=1}^n (z-u_i),
$$
the polynomials $p_{n-3}(z)$, $p_{2n-4}(z)$, \dots, $p_{m\, n-m-n-1}(z)$ labeled by their degrees are to be determined later.
Such a curve will pass through the points (\ref{pointu}), (\ref{pointinf}). To have it passing also through (\ref{pointq})
the following system of equations must be satisfied
\begin{equation}\label{ansatz1}
\sum_{k=2}^m p_{k\, n -n -k -1}(\gamma_s) \hat \mu_s^{m-k} + Q(\gamma_s, \hat\mu_s) = 0, \quad s=1, \dots, g
\end{equation}
where
\begin{eqnarray}
&&
Q(z,\hat w) := {\hat w^m\over R(z)} 
- \sum_{i=1}^{n} {\sigma_1^{(i)}\over z-u_i} \hat w^{m-1}
+ \sum_{k=2}^m \left[ 
(-1)^k \sum_{i=1}^{n} {\sigma_k^{(i)} \over z-u_i}[R'(u_i)]^{k-1}\right.
\nonumber\\
&&
\left. + \sigma_k^{(\infty)} z^{k\,
n-n-k}\right] \hat w^{m-k}.\nn
\end{eqnarray}
This is a linear system for the $g$ coefficients of the polynomials $p_{n-3}(z)$, $p_{2n-4}(z)$, \dots, $p_{m\, n-m-n-1}(z)$. Let us prove
that the determinant of this linear system does not vanish identically. Indeed, this determinant is equal to the following polynomial
in $\gamma_1$, \dots, $\gamma_g$, $\hat\mu_1$, \dots, $\hat\mu_g$
\begin{eqnarray}\label{det}
&&
W_{m,n}(\gamma_1, \dots, \gamma_g, \hat\mu_1, \dots, \hat\mu_g) 
\nonumber\\
&&
:=\sum_\pi (-1)^{|\pi|} (\hat\mu_{i_1}\dots \hat\mu_{i_{n-3}})^{m-2} 
(\hat\mu_{j_1} \dots \hat\mu_{j_{2n-4}})^{m-3}
(\hat\mu_{k_1} \dots \hat\mu_{k_{3n-5}})^{m-4} \dots
\nonumber\\
&&
\times V(\gamma_{i_1}\dots \gamma_{i_{n-3}})V(\gamma_{j_1} \dots 
\gamma_{j_{2n-4}})V(\gamma_{k_1} \dots \gamma_{k_{3n-5}})
V(\gamma_{l_1}, \dots \gamma_{l_{m\, n - m - n -1}}).
\end{eqnarray}
Here the summation is over the partitions
\begin{equation}\label{parti}
\pi:\,  \{ 1, 2, \dots, g\} = \{i_1, \dots, i_{n-3}\} \sqcup \{ j_1, \dots, j_{2n-4}\} \sqcup
\{ k_1, \dots, k_{3n-5}\} \sqcup \dots \sqcup \{ l_1, \dots, 
l_{m\, n - m - n -1}\}
\end{equation}
$|\pi|$ stands for the parity of the permutation $\pi\in S_g$, 
$$
V(x_1, \dots, x_k) := \prod_{1\leq i<j\leq k} (x_i-x_j)
$$
is the Vandermonde determinant. It is clear that this polynomial is not an identical zero. The Theorem is proved.\epf

We want now to show that the same data used in Theorem \ref{interpol} 
determine the matrix valued polynomial $A(z)$ in the determinant representation
(\ref{comden3}). 
We will now prove that any generic curve of the form (\ref{comden7}) can be
represented in the determinant form (\ref{comden3}). Actually, this can be done
in many ways; we will describe the parameters of all determinant representations
of $\Gamma$.

Let $\Gamma$ be the spectral curve (\ref{comden4}) of a matrix $\hat A(z)$.
Assuming smoothness of the spectral curve, we will associate with the
determinant representation a degree $g$ divisor $D$ on $\Gamma$. Let us first consider
the eigenvector line bundle ${\cal L}$ on $\Gamma$
\begin{equation}\label{lb}
\hat A(z) \psi = \hat w \psi,
\end{equation}
$\psi= (\psi_1, \dots, \psi_m)^T$. Define the divisor
\begin{equation}\label{div}
\hat D:= \{ \psi_1=0\}.
\end{equation}
In other words, $\hat D$ is the divisor of poles of the meromorphic functions
$\psi_2/\psi_1$, $\psi_3/\psi_1$, \dots, $\psi_m/\psi_1$.

\begin{lm} The degree of the divisor $\hat D$ is equal to
$$
\deg \hat D = g+m-1.
$$
\end{lm}

See \cite{griffiths} for a simple proof.

Denote
\begin{equation}\label{div1}
D = \hat D\cap \Gamma \setminus D_\infty.
\end{equation}

\begin{lm} The point $(z_0, \hat w_0)\in D$
only if 
$$
{\Delta_0}(z_0)=0.
$$
Here the polynomial ${\Delta_0}(z)$ was defined in (\ref{Delta}).
Conversely, for any root $z_0$ of the polynomial ${\Delta_0}(z)$ 
there exists
a point $(z_0, \hat w_0)\in D$. 
\end{lm}

\noindent{\bf Proof.} For the convenience of the reader we will give here 
the proof. Rewriting the equation of the divisor in the form
$$
<b_0, \psi>=0
$$
(here $<~,~>$ stands for the natural pairing between row- and column-vectors,
the row-vector $b_0$ was defined in (\ref{bk})) we also derive that, for $k=1,
\dots, m-1$
$$
<b_k(z_0), \psi> = <b_0, A^k(z_0) \psi> = w^k_0 <b_0, \psi> =0, 
~~w_0 ={\hat w_0\over
{R}(z_0)}
$$
The determinant of this linear homogeneous system must be equal to 0. 
This gives
${\Delta_0}(z_0)=0$. 

Conversely, let $z_0$ be a zero of ${\Delta_0}$ . Using the identity 
(\ref{id}) we derive
that the subspace ${\mathrm{Ker}}\,B(z_0)$ of the $m$-dimensional space is
invariant w.r.t. the linear operator $A(z_0)$. Here the matrix $B(z)$ was
defined in the line after the formula (\ref{bk}). Therefore there exists an
eigenvector $\psi \in {\mathrm{Ker}}\,B(z_0)$, 
$$
A(z_0) \psi =w_0 \psi, ~~ <b_0, \psi>=\psi_1=0.
$$
By definition the point $(z_0, w_0)$ belongs to the divisor $D\in \Gamma$. The
Lemma is proved. \epf

Let us now compute the degree of the polynomial ${\Delta_0}(z)$. Let
$$
\hat A(z) = -A_\infty z^{{n}-1} +C z^{{n}-2} + O(z^{{n}-3}).
$$
Explicitly, the matrix $C=(C_{ij})$ reads 
\begin{equation}\label{matc}
C = A_\infty \bar u+ \sum_{i=1}^{n} u_i A_i, ~~\bar u = 
\sum_{i=1}^{n} u_i.
\end{equation}

\begin{lm} The polynomial ${\Delta_0}(z)$ has the form
$$
{\Delta_0}(z) = (-1)^{m-1} C_{12} C_{13} \dots C_{1m} \prod_{2\leq i<j}
(\lambda_i^{(\infty)} -\lambda_j^{(\infty)}) \, z^g +O(z^{g-1})
$$
where $g$ is given by the formula (\ref{genus}).
\end{lm}

\noindent{\bf Proof.}   For the $j$-th coordinate of the row-vector
$$
\hat b_k : = b_0 \hat A^{k-1}(z)
$$
(cf. (\ref{bk}) one obtains
$$
(\hat b_k)_j = 
(-1)^k \left[ \left( \lambda_1^{(\infty)}\right)^k \delta_{1\, j}
z^{k({n}-1)} - {(\lambda_1^{(\infty)})^k-(\lambda_j^{(\infty)})^k\over
\lambda_1^{(\infty)}-\lambda_j^{(\infty)}} C_{1\,j} z^{k({n}-1)-1} + 
\dots\right]
$$
where dots stand for the terms of lower order in $z$. Computing the determinant
of this matrix we obtain the proof of the needed formula.\epf

\begin{cor} If the eigenvalues of the matrix $A_\infty$ are pairwise distinct
and all the elements of the first row of the matrix (\ref{matc}) are not equal
to zero then the degree of the divisor $D$ is equal to $g$. The remaining 
points of the divisor $\hat D$ are at infinity
$$
\hat D - D = \infty_2+ \infty_3 + \dots + \infty_m.
$$
\end{cor}

The statement of the Corollary is a formalization of the following asymptotic
behavior of the eigenvectors $\psi=(\psi_1, \dots, \psi_m)^T$ 
of the matrix $\hat A(z)$ at $z\to\infty$:
$$
\psi_k = \delta_{kj} + O\left( {1\over z}\right), ~~ z\to\infty, 
~~{\hat w\over
z^{{n}-1}} \to -\lambda_j^{(\infty)}.
$$
Such normalized eigenvector will have $g+m-1$ poles on $\Gamma \setminus
D_\infty$. 
Under the above assumptions the first component $\psi_1$ has simple zeroes at
the points $\infty_2$, \dots, $\infty_m$.

\begin{rmk} 
We observe that, as it was shown in the proof of Lemma \ref{lmb0}, the element
$C_{ij}$ of the matrix (\ref{matc}) is identically equal to zero for an
isomonodromic deformation ${  A}_k{  A}_k(u_1, \dots, u_{n})$, if and only if
$$
(1-\lambda^{(\infty)}_i+\lambda^{(\infty)}_j){  A}_{k_{ij}}
=-{  A}_{\infty_{ij}},
$$
that is either $1-\lambda^{(\infty)}_i+\lambda^{(\infty)}_j=0$   or
$A_{k_{ij}} =0$ for all $k$.
\end{rmk}

Denote 
\begin{equation}\label{canonic}
\gamma_s = \gamma_s({\bf A}), ~~\mu_s= \mu_s({\bf A}), ~~s=1, \dots, g
\end{equation}
the coordinates of the points of the divisor $D$,
\begin{equation}\label{canonic1}
D=(\gamma_1, \hat \mu_1) + \dots +(\gamma_g, \hat \mu_g), ~~ \hat \mu_s = \mu_s \prod_{i=1}^{n}
(\gamma_s-u_i).
\end{equation}
We want to show that, given the generic values of the functions $\gamma_s({\bf A})$,
$\mu_s({\bf A})$, $s=1, \dots, g$ together with the the numbers $u_1$, \dots, $u_{n}$
and the pairwise distinct eigenvalues $\lambda_1^{(i)}$, \dots,
$\lambda_m^{(i)}$, $i=1, \dots, {n}, \infty$ satisfying the constraint
(\ref{f-cons}) one can uniquely determine
the conjugacy class of the ${n}$-tuple of matrices 
${\bf A}=(A_1, \dots, A_{n})$
modulo diagonal conjugations and permutations. To this end we will prove,
essentially using the technique of \cite{mfgo} the converse
statement that shows that, for a Zariski open subset in the space ${\bf
P}^{M-1}$ of algebraic curves of the form (\ref{comden7}) the curve can be
represented in the determinant form. We will also describe parameters of such 
determinant representations of a given curve.

Let $D$ be a divisor of the degree $g$ on the Riemann surface $\Gamma$. We will
say that the divisor is $D_\infty$-{\it non-special} if 
\begin{equation}\label{nonspec}
\dim H^0(\Gamma, {\cal O}(D+\infty_i -\infty_1)) =1, ~~i=1, \dots, m.
\end{equation}

\begin{thm} \label{ind-th} Any smooth affine curve $\hat {\cal R}(z, \hat w)=0$
of the form (\ref{comden7}) can be
represented in the determinant form (\ref{comden4}) for a matrix $\hat A(z)$
of the form (\ref{comden2}). For a $D_\infty$-generic curve 
such representations, considered modulo diagonal
conjugations
$$
\hat A(z) \mapsto K^{-1} \hat A(z)\, K, ~~ K={\mathrm{diag}} (k_1, \dots, k_m)
$$
and permutations of coordinates
$$
\hat A(z) \mapsto P^{-1} \hat A(z) \,P, ~~ P\in S_{m-1} 
$$
preserving the vector $(1, 0, \dots, 0)$ are in one-to-one correspondence
with the degree $g$ $D_\infty$-non-special divisors $D$ on $\Gamma$.
\end{thm}

\noindent{\bf Proof.}  Let us order the infinite points of $\Gamma$ and 
choose nonzero sections
\begin{equation}\label{sect}
\psi_k \in H^0 (\Gamma, O(D+\infty_k-\infty_1)), ~~k=2, \dots, m.
\end{equation}
They are determined uniquely up to constant factors. Introduce a vector valued
meromorphic function on $\Gamma$ putting
$$
\psi=(1, \psi_2, \dots, \psi_m)^T
$$
Introduce $m \times m$ matrix $\Psi(z)=(\Psi_{k\,j}(z))$ of Laurent series 
in $1/z$ expanding the functions $\psi_k$ near the infinite points 
$\infty_j\in\Gamma$.
Let $W (z)= {\mathrm{diag}}(W_1(z), \dots, W_m(z))$ be the diagonal matrix
obtained by taking the Laurent series of the function $\hat w$ on $\Gamma$
near $\infty_1$, \dots, $\infty_m$. Define a matrix of polynomials
\begin{equation}\label{matr}
\hat A(z) : = \left( \Psi(z) W(z)\Psi^{-1}(z)\right) _+.
\end{equation}
Here $(~)_+$ means the polynomial part in $z$ of the expansion. By construction
$$
\hat A(z) = -z^{m({n}-1)}{\mathrm{diag}}(\lambda_1^{(\infty)}, \dots,
\lambda_m^{(\infty)}) + O(z^{m(n-1)-1}).
$$
Let us prove that the vector function $\psi$ on $\Gamma$ satisfies
$$
\hat A(z) \psi = \hat w\, \psi.
$$
This can be done using the standard arguments of the Krichever's scheme
\cite{krichever}. Indeed, by the construction all the components of the 
difference 
$$
\hat w \psi - \hat A(z) \psi
$$
are analytic at the infinite points $\infty_1$, \dots, $\infty_m$. Therefore they
have poles only at the points of the divisor $D$. Due to nonspeciality of $D$
the difference is equal to zero. 

We have proved that $\Gamma$ coincides with the spectral curve of the polynomial
matrix $\hat A(z)$. By construction the divisor $D$ we started with coincides
with the one defined above.
It remains to observe that, choosing another basic sections
$$
\tilde\psi_k \in H^0 (\Gamma, O(D+\infty_k-\infty_1)), 
$$
$$
\tilde \psi_k = c_k \psi_k, ~~k=2, \dots, m
$$
yields the diagonal conjugation of the polynomial matrix $\hat A(z)$,
$$
\hat A(z) \mapsto C\, \hat A(z)\, C^{-1}, ~~ C={\mathrm{diag}}(1, c_2, \dots,
c_m).
$$
Moreover, changing the order of the infinite points preserving $\infty_1$
implies a permutation.
The theorem is proved. \epf

\begin{cor} The map
\begin{equation}\label{mappa}
\left[ A(z)\right] \mapsto {\mathcal Spec}
\end{equation}
is a birational isomorphism of the space of classes of equivalence of
rational matrix-valued functions of the form
$$
A(z) = \sum_{i=1}^n {A_i\over z-u_i}, \quad A_\infty = - (A_1 + \dots + A_n)
$$
with diagonal $A_\infty$ considered modulo diagonal conjugations and the space
of spectral data with the coordinates 
\begin{equation}\label{spec}
\left( u_1, \dots, u_n, {\mathrm {Spec}}\, A_1,
\dots, {\mathrm {Spec}}\, A_n, {\mathrm {Spec}}\, A_\infty, \gamma_1, \mu_1,\dots,
\gamma_g,  \mu_g\right)\in {\mathcal Spec}.
\end{equation} 
In particular, $(\gamma_1, \mu_1,\dots,
\gamma_g, \mu_g)$ are coordinates on the reduced symplectic leaves of the Poisson
bracket (\ref{ham2}).
\end{cor}

We will now prove that $\gamma_i=\gamma_i({\bf A})$, $\mu_i=\mu_i({\bf A})$, 
$i=1, \dots, g$
are {\it canonical coordinates} on the reduced symplectic leaves of the Poisson bracket.
It will be convenient to represent the Poisson bracket (\ref{ham2}) in the
following wellknown
$r$-matrix form (see \cite{Fad-Tak} regarding the definitions and notations).

\begin{lm}
The Poisson bracket (\ref{ham2}) can be represented in the form
\begin{equation}
\left\{{A}({z_1})\otim_, {A}(z_2)\right\}
= \left[  {A}({z_1})\otim\ID+\ID\otim
{A} (z_2), r({z_1} - z_2)\right] ,   \label{qpb}
\end{equation}
where $r({z})$ is a {\it classical $r$-matrix}, i.e.
a solution of the linearized Yang -- Baxter equation, given by
$$
r_{i k}^{j l}({z})={\delta_i^l\delta_k^j\over{z}}.
$$
Equivalently, (\ref{qpb}) reads
\begin{equation}\label{pb0}
\left\{ A^i_j(z_1), A^k_l(z_2)\right\} ={1\over z_1-z_2}
\left[ \delta^i_l \left(A^k_j(z_1)-A^k_j(z_2)\right) -\left(
A^i_l(z_1)-A^i_l(z_2)\right) \, \delta_j^k\right].
\end{equation}
\label{thm1.2}
\end{lm}

\noindent{\bf Proof.} By the definition
$$
\left\{{A}({z_1})\otim_,{A}(z_2) \right\}^{i k}_{j l}:\sum_p {{{A}_p}_j^k \delta_l^i -{{A}_p}_l^i \delta_j^k 
\over({z_1}-u_p)(z_2-u_p)}.
$$
So
\begin{eqnarray}
&&
\sum_p {{{A}_p}_j^k \delta_l^i -{{A}_p}_l^i \delta_j^k 
\over({z_1}-u_p)(z_2-u_p)}
={1\over({z_2}-z_1)}\sum_p\left({1\over{z_1}-u_p}-
{1\over z_2-u_p}\right)
\left({{A}_p}_j^k \delta_l^i -
{{A}_p}_l^i \delta_j^k \right)
\nonumber\\
&&
={1\over({z_2}-z_1)}\left( {A}^k_j({z_1})\delta_l^i-
{A}^i_l({z_1})\delta_j^k + {A}^i_l(z_2)\delta_j^k-
{A}^k_j(z_2)\delta_l^i\right)
\nonumber\\
&&
=\left[  {A}({z_1})\otim\ID+\ID\otim
{A} (z_2),r({z_1} - z_2)\right]^{i k}_{j l}.
\nn\end{eqnarray}
This concludes the proof. \epf
\vskip 0.3 cm

\begin{thm} The functions $\lambda_i^{(k)}=\lambda_i^{(k)}({\bf A})$, $i=1,
\dots, m$, 
$k=1, \dots, {n}, \infty$, $\gamma_i=\gamma_i({\bf A})$, $\mu_i=\mu_i({\bf A})$, $g=1, \dots, g$
on the space of $m\times m$ matrices $(A_1, \dots, A_{n})=:{\bf A}$
have the following canonical Poisson brackets w.r.t. the structure (\ref{pb0})
$$
\{ \gamma_i, \mu_j\} = \delta_{ij},
$$
all other Poisson brackets vanish.
\end{thm}

\noindent{\bf Proof.}   We already know that the eigenvalues $\lambda_i^{(k)}$ of the matrices 
$A_k$ are Casimirs of the Poisson bracket. It remains to compute the Poisson
brackets of the functions $\mu_i({\bf A})$ and $\gamma_j({\bf A})$.

Let us introduce the following notations. Let $z$ be not
a ramification point for the Riemann surface $\Gamma$. Let us fix some ordering
of the sheets of the Riemann surface. Denote
$$
|a>, ~~a=1, \dots, m
$$
the basis of eigenvectors of the matrix $ A=A(z)$,
\begin{equation}\label{ket}
 A |a> = w_a |a>
\end{equation}
normalized by the condition
\begin{equation}\label{norm1}
<b_0 | a>=1, ~~ b_0 = (1, 0, \dots, 0).
\end{equation}
Here 
$$
w_a=w_a(z), ~~ a=1, \dots, m
$$
are the roots of the characteristic equation $\det(w-A(z))=0$.
Denote $<a|$ the dual basis of row-vectors
\begin{equation}\label{dual}
<a|b> = \delta_{ab}.
\end{equation}
Due to (\ref{norm1}) one has
\begin{equation}\label{b0}
b_0 = \sum_{a=1}^m <a|.
\end{equation}

\begin{lm} The following formulae for the Poisson brackets hold true
\begin{equation}\label{pb1}
\{ \log {\Delta_0}(z), \log{\Delta_0} (z')\}
=0
\end{equation}
\begin{equation}\label{pb3}
\{ w_a(z), w_b(z')\} = 0.
\end{equation}
\begin{equation}\label{pb2}
\{\log{\Delta_0}(z), w_c(z')\} = {1\over (z-z')} 
\sum_{a\neq b} <a| c'> <c'|b>
\end{equation}
\end{lm}

In this formulae the primes mean that the corresponding function is computed
at the point $z'$, i.e.
$$
A(z') |c'> = w'_c |c'> ~~{\rm where}~ w'_c= w_c(z').
$$

\noindent{\bf Proof.} The following wellknown variational formulae will 
be useful in the computations
of the Poisson brackets
\begin{equation}\label{deltaw}
\delta w_a = <a|\delta  A | a>
\end{equation}
\begin{equation}\label{deltaa1}
<b|\delta a> = {<b|\delta  A| a>\over w_a - w_b}, ~~ b\neq a
\end{equation}
\begin{equation}\label{deltaa2}
<a |\delta a> =\sum_{b\neq a} {<b|\delta  A| a>\over w_b - w_a}.
\end{equation}
Here $|\delta a>$ is the variation of the eigenvector $|a>$. In the derivation
of the last formula we have used the normalization (\ref{norm1}). 

Denote $\Psi(z)$ the matrix with the columns $|1>$, \dots, $|m>$. The rows of
the inverse matrix $\Psi^{-1}$ coincide with the bra-vectors $<1|$, \dots,
$<m|$.
From
(\ref{b0}) it easily follows that the matrix $B(z)$ is equal to the product
of the Vandermonde matrix of the pairwise distinct
numbers $w_1$, \dots, $w_m$ by $\Psi^{-1}(z)$.
So
$$
\det B(z) = \prod_{i<j} (w_i-w_j) {\det}^{-1}\Psi(z).
$$
Using the Liouville formula
$$
\delta\log(\det\Psi) ={\mathrm{tr}}\, \Psi^{-1}\delta\Psi = \sum_{a=1}^m <a|\delta a>
$$
yields
\begin{equation}\label{deltadelta}
\delta \log{\Delta_0}(z) = {1\over 2}\sum_{a\neq b} {<a| \delta  A | a>
-<b|\delta A|b>
-<a| \delta A
| b>\over w_a - w_b}.
\end{equation}
It is understood that the values of the variables $u_i$ are fixed during the
variation. From the above formulae for $\delta{\Delta_0}(z)$, $\delta w_c(z')$ we
derive, by a straightforward calculation, the brackets (\ref{pb1}) -
{\ref{pb2}). The Lemma is proved. \epf 

\noindent{\bf Proof of the Theorem.}
From (\ref{pb1}) and from the representation
$$
\delta\log{\Delta_0}(z) = - \sum_{i=1}^m {\delta \gamma_i\over z-\gamma_i} + {\rm
regular}~{\rm terms}
$$
it easily follows that
$$
\{ \gamma_i, \gamma_j\} = 0.
$$
The commutation rule 
$$
\{\mu_i, \mu_j\}=0
$$
follows from (\ref{pb3}). Let us compute the brackets $\{ \gamma_i, \mu_j\}$. 
Due to Theorem \ref{ind-th} we may assume that the projections $z=\gamma_i$
of the points
of the divisor $D$ onto the $z$-plane are all pairwise distinct, they are
distinct from $u_j$ and from the ramification points of the Riemann surface
(cf. Assumption 2 above). Consider
first the case $j\neq i$.
Assume
that the numeration of the sheets of the Riemann surface at the neighborhoods
of the points $z=\gamma_i$ and $z'=\gamma_j$ is done in such a way that the pole of the
eigenvector $\psi$ of the matrix $A(z)$ belongs to the sheet labeled by $c$.
That means that the ket-vector $|c>$ has a simple pole at $z\to \gamma_i$,
\begin{equation}\label{pole1}
|c> = {|\tilde c>\over z-\gamma_i} + O(1), ~~ z\to \gamma_i.
\end{equation}
All other ket- and bra-vectors are analytic in $z$ near this point and 
the corresponding bra-vector $<c\,|$ has a simple zero
\begin{equation}\label{pole2}
<c\,| = (z-\gamma_i) <\tilde c| +O\left((z-\gamma_i)^2\right).
\end{equation}

From the already proven commutation rule of the
coordinates $\gamma_i$ it follows that
$$
\begin{array}{ll}
\{ \gamma_i, \mu_j\}& = -\lim_{z\to \gamma_i} \lim_{z'\to \gamma_j} 
(z-\gamma_i)\{ \log{\Delta_0}(z), w_c(z')\}=\\
&= - {1\over 2 (\gamma_i-\gamma_j)} \lim_{z\to \gamma_i} 
\lim_{z'\to \gamma_j} (z-\gamma_i)\sum_{a\neq b} <a|c'><c'|b>.\\
\end{array}
$$
Due to (\ref{pole1}), (\ref{pole2}) the r.h.s. is analytic at the point $z=\gamma_j$.
The singularity at $z=\gamma_i$ can come only from the terms with $b=c$. So, the
singular part in the sum equals
$$
\sum_{a\neq b} <a|c'><c'|b>= \sum_{a\neq c} <a|c'><c'|c> +{\rm regular}.
$$
Using the $z$-independent normalization (\ref{norm1}) we rewrite the singular
term in the form
$$
\sum_{a\neq c} <a|c'><c'|c>= -<c|c'><c'|c>.
$$
Using again (\ref{pole1}), (\ref{pole2}) we establish analyticity of the last
expression also at $z=\gamma_i$. Therefore $\{ \gamma_i, \mu_j\} =0$ for $i\neq j$.

To compute $\{ \gamma_i, \mu_i \}$ we will first calculate the limit
$$
\lim_{z'\to z} \{ \log{\Delta_0}(z), w_c(z')\}.
$$
Observe that, for $z'=z$ the numerator of the formula (\ref{pb2}) vanish:
$$
\left( \sum_{a\neq b} <a|c'><c'|b>\right)_{z'=z} =  \sum_{a\neq b} <a|c><c|b>
=\sum_{a\neq b}\delta_{ac}\delta_{cb}=0.
$$
So the needed limit is equal to the derivative
$$
\lim_{z'\to z} \{ \log{\Delta_0}(z), w_c(z')\}=-\left({d\over dz'} 
\sum_{a\neq b} <a|c'><c'|b>\right)_{z'=z}.
$$
Let us denote
$$
|\dot c> := {d\over dz'} |c'>\large |_{z'=z}, ~~ <\dot c| := {d\over d z'}
<c'|\, \large |_{z'=z}.
$$
We obtain, using again the $z'$-independent normalization (\ref{norm1})
$$
{d\over dz'} \left( \sum_{a\neq b} <a|c'><c'|b>\right)_{z'=z} 
=\sum_{a\neq c} <a| \dot c>  +\sum_{b\neq c}<\dot c|b>=-<c|\dot c> 
+\sum_{b\neq c}<\dot c|b>.
$$
Therefore
$$
\{ \gamma_i, \mu_i\} = \lim_{z\to \gamma_i} (z-\gamma_i) \left[ 
-<c|\dot c> 
+\sum_{b\neq c}<\dot c|b>\right].
$$
The last term in the brackets is analytic at the point $z=\gamma_i$. For the first
one we obtain, using (\ref{pole1}), (\ref{pole2})
$$
<c|\dot c> =-{<\tilde c| \tilde c>\over z-\gamma_i} +\, {\rm regular} ~{\rm terms}.
$$
The last step is to use the normalization 
$$
<c|c>\equiv 1
$$
to derive that
$$
<\tilde c| \tilde c>=1.
$$
Hence
$$
\{ \gamma_i , \mu_i\} =1.
$$
The Theorem is proved. \epf

\bibliography{bibgenerale.bib}
\bibliographystyle{plain}

\end{document}